\documentclass[12pt,leqno]{article}

\usepackage{amsmath,amsfonts,amssymb,epsfig,theorem}
\usepackage[usenames]{color}

\usepackage[all]{xy}
\usepackage{psfrag}
\usepackage{multirow}

\DeclareFontFamily{U}{rsf}{}
\DeclareFontShape{U}{rsf}{m}{n}{
  <5> <6> rsfs5 <7> <8> <9> rsfs7 <10->  rsfs10}{}
\DeclareMathAlphabet{\mathscr}{U}{rsf}{m}{n}
\newcommand{\mycal}[1]{\mathscr{#1}}
\newcommand{\op}[1]{\operatorname{#1}}

\newcommand{\bS}{\boldsymbol{S}}
\newcommand{\bD}{\boldsymbol{D}}
\newcommand{\bO}{\boldsymbol{O}}
\newcommand{\bo}{\boldsymbol{o}}
\newcommand{\bj}{\boldsymbol{j}}
\newcommand{\stS}{\text{\sf\bfseries S}}
\newcommand{\nc}{{\sf\bfseries nc}}
\newcommand{\ncSpec}{\op{\text{{\sf\bfseries ncSpec}}}}
\newcommand{\iso}{\op{\text{\sf\bfseries iso}}}
\newcommand{\bbf}{\boldsymbol{f}}
\newcommand{\QncHS}{(\op{\mathbb{Q}-\text{\nc{\sf{HS}}}})}
\newcommand{\QHS}{(\op{\mathbb{Q}-HS})}
\newcommand{\be}{\boldsymbol{\mathcal{E}}}
\newcommand{\bT}{\boldsymbol{\mathcal{T}}}
\newcommand{\bc}{\boldsymbol{c}}
\newcommand{\bM}{\boldsymbol{M}}
\newcommand{\perv}{\text{{\sf Perv}}}
\newcommand{\bi}{\boldsymbol{i}}
\newcommand{\bp}{\boldsymbol{p}}

\newcommand{\bw}{{\sf w}}
\newcommand{\bstar}{\boldsymbol{*}}
\newcommand{\sGr}{\text{\sf\bfseries Gr}}
\newcommand{\sDel}{\text{\sf\bfseries Del}}
\newcommand{\MF}{\text{\sf\bfseries MF}}
\newcommand{\bDelta}{\boldsymbol{\Delta}}
\newcommand{\bexp}{\text{\bfseries exp}}
\newcommand{\mir}{\boldsymbol{|}}

\newcommand{\bPi}{\boldsymbol{\Pi}}
\newcommand{\MM}{\op{\mathbb{M}\text{{\sf od}}}}
\newcommand{\bOmega}{\boldsymbol{\Omega}}
\newcommand{\ddiv}{\op{{\sf div}}}
\newcommand{\bmu}{\boldsymbol{\mu}}

\newcommand{\Perf}{\op{{\sf Perf}}}
\newcommand{\FS}{\op{{\sf FS}}}
\newcommand{\Fuk}{\op{{\sf Fuk}}}

\newcommand{\bxi}{\boldsymbol{\Xi}}
\newcommand{\blX}{\boldsymbol{\mathcal{X}}}
\newcommand{\bbgamma}{\boldsymbol{\gamma}}

\theoremstyle{plain}

\newtheorem{theo}{Theorem}[section]
\newtheorem{lemma}[theo]{Lemma}

\newtheorem{conn}[theo]{Conjecture}
\newtheorem{defi}[theo]{Definition}
\newtheorem{prop}[theo]{Proposition}

\newtheorem{claim}[theo]{Claim}
\newtheorem{prop-defi}[theo]{Proposition-Definition}
\newtheorem{lemma-defi}[theo]{Lemma-Definition}
{\theorembodyfont{\rmfamily} \newtheorem{rem}[theo]{Remark}}
{\theorembodyfont{\rmfamily} \newtheorem{ex}[theo]{Example}}
{\theorembodyfont{\rmfamily} \newtheorem{que}[theo]{Question}}

\newlength{\miniwidth}

\numberwithin{equation}{subsection}

\begin{document}

\title{{\bfseries Hodge theoretic aspects \\ of mirror symmetry}}
\author{L.Katzarkov \and M.Kontsevich \and T.Pantev}
\date{ }
\maketitle

\begin{abstract} We discuss the Hodge theory of algebraic non-commutative
spaces and analyze how this theory interacts with the Calabi-Yau
condition and with mirror symmetry.  We develop an abstract theory of
non-commutative Hodge structures, investigate existence and
variations, and propose explicit construction and classification
techniques. We study the main examples of non-commutative Hodge
structures coming from a symplectic or a complex geometry possibly
twisted by a potential. We study the interactions of the new Hodge
theoretic invariants with mirror symmetry and derive non-commutative
analogues of some of the more
interesting consequences of Hodge theory such as unobstructedness and
the construction of canonical coordinates on moduli.
\end{abstract}

\tableofcontents

\section{Introduction} \label{sec:introduction}

This paper is a first in a series aiming to develop a general
procedure associating a 2-dimensional cohomological field theory in
the sense \cite{kontsevich-manin} (CohFT in short) to a certain
structure in derived algebraic geometry.  More precisely, for any
Calabi-Yau $A_{\infty}$-category satisfying appropriate finiteness
conditions (smoothness and compactness), and such that a
noncommutative analog of the Hodge $\Rightarrow$ de Rham spectral
sequence collapses, we associate an infinite-dimensional family of
CohFTs.  The additional parameters needed to specify the CohFT are of
a purely cohomological nature. Conjecturally, our procedure applied to
the Fukaya category should give (higher genus) Gromov-Witten
invariants of the underlying symplectic manifold. 

This program was first outlined by the second author in several talks
given in 2003-2004, and some aspects of it were later developed in
depth by
K.Costello
\cite{costello-CY,costello-GW,costello-gauge,costello-renormalization}.
The whole picture turns out to be very intricate, and in the process
of writing we realized that we have to address a large variety of
problems. In this installment we do not discuss the general plan of
our approach but rather focus on those features of $A_{\infty}$ or dg
categories that can be captured by Hodge theoretic constructions.  We
propose a formalism that starts with Homological Mirror Symmetry and
extrapolates a geometric picture for the requisite categories that
makes them amenable to study via old and new Hodge theory.  Our hope
is that this geometric treatment will provide new invariants and will
expand the scope of possible applications in symplectic geometry and
algebraic geometry.

Mirror symmetry was introduced in physics as a special duality between
two $N=2$ super conformal field theories.  Traditionally a $N=2$ super
conformal field theory is constructed as a quantization of a
non-linear $\sigma$-model with target a compact Calabi-Yau manifold
equipped with a Ricci flat K\"{a}hler metric and a closed $2$-form -
the so called $B$-field. Two Calabi-Yau manifolds $X$
and $Y$ form a {\em mirror pair} $X\mir Y$ if the associated $N=2$
super conformal field theories are mirror dual to each other
\cite{cox-katz}.

Homological Mirror Symmetry was introduced in 1994 by the second
author for the case of Calabi-Yau manifolds but today the realm of its
applicability is much broader.  In particular many of our
considerations in the present work are governed by an analogue of
Homological Mirror Symmetry for geometries with potentials.  We
study the effect of such mirror symmetry on the associated categories
of $D$-branes and especially on the associated non-commutative Hodge
structures on homological invariants, i.e. on the Hochschild and
cyclic homology and cohomology of such categories. We study mirror
pairs consisting of a compact manifold on one side, and of a
Landau-Ginzburg model with a proper potential on a non-compact
manifold having $c_{1} = 0$ on the other. We formulate the mirror
symmetry conjecture on the Hodge theoretic level in both directions.
That is, we compare the non-commutative Hodge structures associated with a
compact complex manifold and a mirror holomorphic Landau-Ginzburg
model, and also the non-commutative Hodge structures associated with a
compact complex manifold with a chosen smooth anticanonical divisor
and with the mirror symplectic Landau-Ginzburg model. This picture is
clearly non-symmetric and has to be generalized. In order to completely
understand the Hodge theoretic aspect of mirror symmetry, one will
have to allow for non-trivial potentials on both sides of the duality
and include the novel toric dualities between formal Landau-Ginzburg
models of Clarke \cite{patrick} and the new smooth variations of
non-commutative Hodge structures of Calabi-Yau  type that we attach to
anticanonical $\mathbb{Q}$-divisors in section
\ref{ssec:anticanonical}. We plan to return to such a generalization
in a future work.

Due to its foundational nature the paper comes out somewhat long winded
and technical for which we apologize. It is organized in three major
parts:

The first part introduces and develops the abstract theory of
non-commutative (\nc )Hodge structures.  This theory is a variant of
the formalism of semi-infinite Hodge structures that was introduced by
Barannikov  \cite{b01,b02,b02b}. We discuss the general
theory of \nc-Hodge structures in the abstract and analyze the various
ways in which the Betti, de Rham and Hodge filtration data can be
specified. In particular we compare \nc \ and ordinary Hodge theory and
explain how \nc-Hodge theory fits within the setup of categorical
non-commutative geometry. We also pay special attention to the
\nc-aspect of Hodge theory and its interaction with the classification
of irregular connections on the line via topological 
data. One of the most useful technical results in this part is the
gluing Theorem~\ref{theo:gluenc} which allows us to assemble \nc-Hodge
structures out of some simple  geometric ingredients. This theorem is
used later in the paper for constructing \nc-Hodge
structures attached to geometries with a
potential.

The second part explains how symplectic and complex geometry give rise
to \nc-Hodge structures and how these structures can be viewed as
interesting invariants of Gromov-Witten theory, projective geometry
and the theory of algebraic cycles. In particular we analyze the Betti
part of the \nc-Hodge theory of a projective space (viewed as a
symplectic manifold) and use this analysis to propose a general
conjecture for the integral structure on the cohomology of the Fukaya
category of a general compact symplectic manifold. The formula for the
integral structure uses only genus zero Gromov-Witten invariants and
characteristic classes of the tangent bundle. Our conjecture is in
complete agreement with the recent work of Iritani
\cite{iritani-integral} who made a similar proposal based on mirror
symmetry for toric Fano orbifolds. We also discuss in detail the
origin of the Stokes data for holomorphic geometries with potentials
and investigate the possible categorical incarnations of this data.

In the third part we study \nc-Hodge structures and their variations
under the Calabi-Yau condition. We extend and generalize the standard
treatment of the deformation theory of Calabi-Yau spaces in order to
get a theory which works equally well in the \nc-context and to be
able to properly define the canonical coordinates in Homological
Mirror Symmetry. We approach the deformation-obstruction problem both
algebraically and by Hodge theoretic means and we obtain
unobstructedness results, generalized pre Frobenius structures and
some interesting geometric properties of period domains for \nc-Hodge
structure. We also study global and infinitesimal deformations and
describe different constructions of Betti and de Rham \nc-Hodge data
for ordinary geometry, relative geometry, geometry with potentials and
abstract \nc-geometry.

\

\bigskip

\noindent
{\bfseries Acknowledgments:} Throughout the preparation of this work
we have benefited from discussions with many people who generously
shared their thoughts and insights with us. Special thanks are due to
A D. Auroux, M. Abouzaid, .Bondal, R.Donagi, V.Golyshev, M.Gross,
A.Losev, D.Orlov, C.Simpson, Y.Soibelman, Y.Tschinkel, A.Todorov, and
B.To\"{e}n for expert help, encouragement and advice. We would also
like to thank the University of Miami for providing the productive
research environment in which most of this work was done. During
various stages of this work we have enjoyed the hospitality of several
outstanding research institutions. We thank the IAS, the IHES, the
Centre Interfacultaire Bernoulli at the EPFL, and the ESI for the
excellent working conditions they have provided. The first and third
author would especially like to thank the organizers of the conference
``From tQFT to $tt^{*}$ and integrability'' at the University of
Augsburg, for giving them an opportunity to speak and for the
invitation to contribute to the proceedings volume of the conference.

During the preparation of this work Ludmil Katzarkov was partially
supported by the Focused Research Grant DMS-0652633 and a research
grant DMS-0600800 from the National Science Foundation, and a FWF
grant  P20778. Tony Pantev was partially supported by NSF Focused Research
Grant DMS-0139799, NSF Research Training Group Grant DMS-0636606, and
NSF grant DMS-0700446.

\section{Non-commutative Hodge structures} \label{sec:nc}

In this section we will discuss the notion of a pure non-commutative
(\nc) Hodge structure. The \nc-Hodge structures are analogues of the
classical notion of a pure Hodge structure on a complex vector
space. Both the \nc-Hodge structures discussed presently and Simpson's
non-abelian Hodge structures \cite{simpson-hf} generalize classical
Hodge theory. In Simpson's theory, one allows for non-linearity in the
substrate of the Hodge structure: the non-abelian Hodge structures of
\cite{simpson-hf} are given by imposing Hodge and weight filtrations
on non-linear topological invariants of a K\"{a}hler space, e.g. on
cohomology with non-abelian coefficients, or on the homotopy type. In
contrast the \nc-Hodge structures discussed in this paper consist of a
novel filtration-type data (the twistor structure of
\cite{simpson-mts,hertling-crelle,sabbah-twistor}) which are still
specified on a vector space, e.g. on the periodic cyclic homology of
an algebra.

Similarly to ordinary
Hodge theory \nc-Hodge structures arise naturally on the de Rham
cohomology of non-commutative spaces of categorical origin.

\subsection{\bfseries Hodge structures} \label{sec-hodge-structures}

We will give several different descriptions of a \nc-Hodge structure
in terms of local data. We begin with the notion of a rational and
unpolarized \nc-Hodge structures, ignoring for the time being the
existence of polarizations and integral lattices. 

\

\bigskip

\subsubsection {\bfseries Notation} \ The \nc-Hodge structures will be
described in terms of geometric data on the punctured complex line, so
we fix once and for all a coordinate $u$ on $\mathbb{C}$ and the
compactification $\mathbb{C} \subset \mathbb{P}^{1}$.  We will write
$\mathbb{C}[[u]]$ for the algebra of formal power series in $u$, and
$\mathbb{C}((u))$ for the field of formal Laurent series in
$u$. Similarly, we will write $\mathbb{C}\{u\}$ for the algebra of
power series in $u$ having positive radius of convergence, and
$\mathbb{C}\{u\}[u^{-1}]$ for the field of meromorphic Laurent series
in $u$ with a pole at most at $u = 0$.

\

\bigskip

\subsubsection {\bfseries Meromorphic connections on the complex line} \ We
will need some standard notions and facts related to meromorphic
connections on Riemann surfaces. We briefly recall those next. More
details can be found in e.g. \cite[chapter~II]{sabbah-frobenius}.

Let $\mycal{M}$ be a finite dimensional vector space over
$\mathbb{C}\{u\}[u^{-1}]$, and let $\nabla$ be a {\em\bfseries
  meromorphic connection} on $\mycal{M}$. Explicitly $\nabla$ is given
by a $\mathbb{C}$-linear map $\nabla_{\frac{d}{du}} : \mycal{M} \to
\mycal{M}$ which satisfies the Leibniz rule for multiplication by
elements in $\mathbb{C}\{u\}[u^{-1}]$. A {\em\bfseries holomorphic
  extension of
  $\mycal{M}$} is a  free finitely generated
$\mathbb{C}\{u\}$-submodule $\mycal{H} \subset \mycal{M}$, such that
\[
\mycal{H}\otimes_{\mathbb{C}\{u\}} \mathbb{C}\{u\}[u^{-1}] =
  \mycal{M}.
\]
Traditionally (see e.g. \cite[section~0.8]{sabbah-frobenius})
a holomorphic extension is called a {\em lattice}. We will avoid this
classical terminology since it may lead to confusion with the integral
lattice structures that we need.

As usual the data $(\mycal{M},\nabla)$ or $(\mycal{H},\nabla)$ should
be viewed as local models for geometric data on a Riemann surface:
$(\mycal{M},\nabla)$ is the local model of a germ of a meromorphic
bundle with connection on a Riemann surface, and $(\mycal{H},\nabla)$
is the local model of a holomorphic bundle with meromorphic connection
on a Riemann surface.

Suppose $(\mycal{M},\nabla)$ is a meromorphic bundle with connection
over $\mathbb{C}\{u\}[u^{-1}]$ and let \linebreak $\mycal{H} \subset
\mycal{M}$ be a holomorphic extension.  We say that $\mycal{H}$ is
      {\bfseries logarithmic} with respect to $\nabla$ if
      $\nabla(\mycal{H}) \subset \mycal{H}\frac{du}{u}$. We say that
      $(\mycal{M},\nabla)$ has at most a {\bfseries regular
        singularity at $0$} if we can find a holomorphic
      extension $\mycal{H} \subset \mycal{M}$ which is logarithmic
      with respect to $\nabla$.

\

\begin{rem} \label{rem-algebraize-lattice}
The Riemann-Hilbert correspondence implies (see
e.g. \cite[II.3.7]{sabbah-frobenius})  that the functor 
of taking the germ at $0 \in \mathbb{P}^{1}$:
\[
\xymatrix@1{
\left(
\text{
\begin{minipage}[c]{2.5in}
finite rank algebraic vector bundles with connections on
$\mathbb{A}^{1}-\{0\}$ with a regular
singularity at $\infty$
\end{minipage}
}
\right) 
\ar[r]^-{{\mathfrak{G}_{0}}}
&
\left(
\text{
\begin{minipage}[c]{2.2in}
finite dimensional $\mathbb{C}\{u\}[u^{-1}]$-vector spaces 
with meromorphic connections
\end{minipage}
}
\right)
}
\]
 is an equivalence of
categories. For future reference we choose once and for all an inverse
$\mathfrak{B}_{0}$ of $\mathfrak{G}_{0}$. We will call
$\mathfrak{B}_{0}$ the {\bfseries algebraization functor} or the
{\bfseries Birkhoff extension functor}. 

Suppose $\mycal{H}$ is a {\em holomorphic bundle} over
$\mathbb{C}\{u\}$ equipped  with a meromorphic connection
$\nabla$. Let $\mycal{M} =
\mycal{H}\otimes_{\mathbb{C}\{u\}} \mathbb{C}\{u\}[u^{-1}]$ and let
$(M,\nabla) = \mathfrak{B}_{0}(\mycal{M},\nabla)$ be the corresponding
Birkhoff extension. The algebraic bundle $M$ on $\mathbb{A}^{1} -
\{0\}$ admits a natural extension to a holomorphic bundle $H$ on
$\mathbb{A}^{1}$: on a small punctured disc centered at $0 \in
\mathbb{A}^{1}$, the bundle $M$ is analytically isomorphic to
$\mycal{M}$, and so $\mycal{H}$ gives us an extension to
the full disc. In particular $\mathfrak{G}_{0}$ and $\mathfrak{B}_{0}$
can be promoted to a pair of inverse equivalences
\[
\xymatrix@1{
\left(
\text{
\begin{minipage}[c]{2.5in}
finite rank
algebraic vector bundles on $\mathbb{A}^{1}$  equipped with a
meromorphic  connection with poles at
most at $0$ 
and $\infty$, and a regular singularity at $\infty$
\end{minipage}
}
\right) 
\ar@<1ex>[r]^-{{\mathfrak{G}_{0}}}
&
\left(
\text{
\begin{minipage}[c]{2.2in}
finite rank free $\mathbb{C}\{u\}$-modules 
equipped with a meromorphic connection
\end{minipage}
}
\right)
\ar@<1ex>[l]^-{\mathfrak{B}^{0}}
}
\]
which we will denote again by $\mathfrak{G}_{0}$ and
$\mathfrak{B}_{0}$. 
\end{rem}

\

\bigskip

\subsubsection {\bfseries Stokes data} \ Let $(\mycal{H},\nabla)$ be a
holomorphic bundle with meromorphic connection over
$\mathbb{C}\{u\}$. We will need the Deligne-Malgrange description of
the associated meromorphic connection $(\mycal{M},\nabla)$ 
via a filtered sheaf on the
circle. We briefly 
recall this description next. More details can be found in
\cite{malgrange-classification} and \cite{babbitt-varadarajan-memoir}.
Let $(M,\nabla) := \mathfrak{B}_{0}((\mycal{M},\nabla))$ be the
Birkhoff extension of $(\mycal{M},\nabla)$ to
$\mathbb{P}^{1}$. Consider the circle $\bS^{1} :=
\mathbb{C}^{\times}/\mathbb{R}^{\times}_{+}$. The sheaf of local
$\nabla$-horizontal sections of $M^{\op{an}}$ on $\mathbb{C}^{\times}$ is a
locally constant sheaf on $\mathbb{C}^{\times}$, which by
contractability of $\mathbb{R}^{\times}_{+}$ induces a locally
constant sheaf $\stS$ of $\mathbb{C}$-vector spaces on $\bS^{1}$.

The sheaf $\stS$ is equipped with a natural local filtration
by subsheaves $\{\stS_{\leq \omega}\}_{\omega \in \sDel}$, where
\begin{itemize}
\item[{\bfseries (i)}] $\sDel$ is the complex local system on
  $\bS^{1}$ for which for every open $U \subset \bS^{1}$ the space of sections 
$\sDel(U)$ is defined to be the space of all holomorphic one forms
  $\omega$ on the sector
\[
\op{Sec}(U) := \left\{ \left. re^{i\varphi} \right| r >0, \varphi \in
U\right\} 
\]
which are of the form
\[
\omega = \left( \sum_{\substack{a \in \mathbb{Q} \\ a < -1 }}
c_{a}u^{a}\right)du, 
\]
where at most finitely many $c_{a} \neq 0$ and the branches $u^{a}$
are chosen arbitrarily. 

Note that the germs of sections of $\sDel$ are naturally ordered: if
$\omega', \omega'' \in \sDel(U)$, $\varphi \in U$, and if 
\[
\omega' - \omega'' = c_{a} u^{a} + \left(\text{\begin{minipage}[c]{0.7in}
  higher order terms\end{minipage}}\right),
\]
then one sets
\[
\omega' <_{\varphi} \omega'' \quad \Leftrightarrow \quad 
\op{Re} \left( \frac{c_{a} e^{i\varphi (a+1)}}{a+1}\right) < 0.
\]
\item[{\bfseries (ii)}] For every $\varphi \in \bS^{1}$ and every
  $\omega \in \sDel_{\varphi}$ the stalk 
\[
\left(\stS_{\leq \omega}\right)_{\varphi} \subset \stS_{\varphi} 
\]
is defined to be the subspace
\[
\left(\stS_{\leq \omega}\right)_{\varphi}  := \left\{ s \in
\stS_{\varphi} = 
\Gamma\left(\mathbb{R}^{\times}_{+}e^{i\varphi},M^{\op{an}}\right)^{\nabla}
\;\; 
\left| \;\; \text{\begin{minipage}[c]{2in} $e^{-\int\omega}s$ has
    moderate growth in the direction $\varphi$, i.e. 
\[
\left\|e^{-\int\omega}s\right\|_{|\mathbb{R}^{\times}_{+}e^{i\varphi}}
= \bO\left(r^{-N}\right),
\]
when $r \to 0$, $N\gg 0$.
\end{minipage}}
\right. \right\}
\]
Here $\|\bullet\|$ is the Hermitian norm of a section of $M$ computed
in some (any) meromorphic trivialization of $M^{\op{an}}$ near $u = 0$. 
\end{itemize}

\

\begin{defi} \label{defi-Stokes-filtration} 
The filtration we just defined is the  {\bfseries Deligne-Malgrange-Stokes
  filtration}, and the $\sDel$-filtered sheaf $\stS$ is called the
{\bfseries Stokes structure} associated to $(\mycal{M},\nabla)$.  
\end{defi}

\

\begin{rem} \label{rem-Stokes-property} The Deligne-Malgrange-Stokes
  filtration satisfies the following property.  
First of all, there exists a covariantly local system of finite sets
$\sDel_{(\mycal{M},\nabla)}\subset \sDel$ canonically associated with
$(\mycal{M},\nabla)$ such that the filtration by all of $\sDel$ is
determined by a filtration by all  $\sDel_{(\mycal{M},\nabla)}(U)$ and all
consecutive factors are non-trivial at all points of $\bS^{1}$ except
finitely many (called the directions of the  Stokes rays). Outside the
 Stokes rays the set $\sDel_{\mycal{M},\nabla}(\phi)$ is totally
ordered. It is easy to see that when we cross a Stokes ray then the
order changes by flipping the order on several disjoint intervals
(e.g. $\{1,2,3,4,5,6\}\to \{2,1,3,6,5,4\}$).  Moreover, on the
subquotients corresponding to these intervals, two filtrations coming
from the left and from the right of the anti Stokes ray are opposed to each
other.  This implies that the graded pieces with respect to the
Deligne-Malgrange-Stokes filtration are locally constant sytems of
vector spaces on the total space of stalks of the sheaf
$\sDel_{(\mycal{M},\nabla)}$ (which is a disjoint union of finite
coverings of $\bS^{1}$).
\end{rem}

\begin{rem} \label{rem:Deligne-Malgrange}
A fundamental theorem of Deligne and Malgrange
\cite[Theorem~4.2]{malgrange-classification},
\cite[Theorem~4.7.3]{babbitt-varadarajan-memoir} asserts that the
functor $(\mycal{M},\nabla) \mapsto \left(\stS, \left\{\stS_{\leq
\omega} \right\}_{\omega \in \sDel}\right)$ is an equivalence between
the category of meromorphic connections over $\mathbb{C}\{u\}[u^{-1}]$
and the category of $\sDel$-filtered local systems on $\bS^{1}$
satisfying the conditions described in
Remark~\ref{rem-Stokes-property}.  We will use this equivalence to
define the Betti part of a \nc-Hodge structure.
\end{rem}

\

\bigskip

\subsubsection {\bfseries The definition of a \nc-Hodge structure} \
After these preliminaries we are now ready to define \nc-Hodge
structures.

\

\begin{defi} \label{defi-nchodge} A  {\bfseries rational pure
    \nc-Hodge structure} 
  consists of the data
  $(H,\mycal{E}_{B},\iso)$, where 
\begin{itemize}
\item $H$ is a  $\mathbb{Z}/2$-graded
algebraic vector bundle on $\mathbb{A}^{1}$. 
\item $\mycal{E}_{B}$ is a local system of finite dimensional
  $\mathbb{Z}/2$-graded $\mathbb{Q}$-vector spaces on \linebreak
  $\mathbb{A}^{1}-\{0\}$. 
\item $\iso$ is an {\em analytic}  isomorphism  of holomorphic vector
  bundles on  
$\mathbb{A}^{1} - \{ 0 \}$: 
\[
\iso : \mycal{E}_{B}\otimes \mathcal{O}_{\mathbb{A}^{1} - \{ 0 \}}
  \stackrel{\cong}{\longrightarrow}
  H_{|\mathbb{A}^{1} - \{ 0 \}}. 
\]
\end{itemize}
\

\noindent
{\bfseries Note:} The isomorphism $\iso$ induces a natural flat
holomorphic connection $\nabla$  on  $ H_{|\mathbb{A}^{1} - \{ 0
  \}}$. 

\

\noindent
These data have to satisfy the following axioms:
\begin{description}
\item[{\sf\bfseries (\nc-filtration axiom)}] $\nabla$ is a meromorphic
  connection on $H$ with a pole of order 
  $\leq 2$ at $u 
  = 0$ and a regular singularity at $\infty$. More precisely, there exist:
\begin{itemize}
\item  a  holomorphic frame of $H$ near $u = 0$ in which 
\[
\nabla = d + \left(\sum_{k \geq -2} A_{k} u^{k}\right)du
\]
with $A_{k} \in \op{Mat}_{r\times r}(\mathbb{C})$, $r =
\op{rank} \ H$.
\item a  meromorphic frame of $H$ near $u = \infty$ in which 
\[
\nabla = d + \left(\sum_{k \geq -1} B_{k} u^{-k}\right)d(u^{-1})
\]
and $B_{k} \in \op{Mat}_{r\times r}(\mathbb{C})$.
\end{itemize}
\item[{\sf\bfseries ($\mathbb{Q}$-structure axiom)}] The
  $\mathbb{Q}$-structure $\mycal{E}_{B}$ on $(H,\nabla)$ is compatible
  with Stokes data. More precisely, let $\left(\stS, \left\{\stS_{\leq \omega}
\right\}_{\omega \in \sDel}\right)$ be the Stokes structure
corresponding to the germ \linebreak $(\mycal{H},\nabla) :=
\mathfrak{G}_{0}(H,\nabla)$, and let  $\stS_{B} \subset \stS$ be the
  $\mathbb{Q}$-structure on $\stS$ induced from
  $\mycal{E}_{B}$ via the isomorphism $\iso$. We require that the
Deligne-Malgrange-Stokes 
  filtration on $\stS$ is defined over $\mathbb{Q}$, i.e.
\[
\left(\stS_{\leq \omega} \cap \stS_{B}\right) \otimes_{\mathbb{Q}}
\mathbb{C} = \stS_{\leq \omega}
\]
for all  local sections $\omega \in \sDel$.
\item[{\sf\bfseries (opposedness axiom)}] The $\mathbb{Q}$-structure
  $\stS_{B}$ induces a real structure on $\stS$ and hence a complex
  conjugation $\tau : \stS \to \stS$. Let $\widehat{H}$ be the
  holomorphic bundle on $\mathbb{P}^{1}$ obtained as the gluing of
  $H^{\op{alg}}_{|\{ |u| \leq 1 \}}$ and
  $\left(\gamma^{*}\overline{H^{\op{alg}}}\right)_{|\{ |u| \geq 1 \}}$
  via $\tau$, where where
$\gamma : \mathbb{P}^{1} \to \mathbb{P}^{1}$ is the real structure on
$\mathbb{P}^{1}$ which fixes the unit circle, i.e. $\gamma(u) :=
1/\bar{u}$. Then we require that $\widehat{H}$ be holomorphically
  trivial, i.e. $\widehat{H} \cong
  \mathcal{O}_{\mathbb{P}^{1}}^{\oplus r}$.
\end{description}

\noindent
A morphism $\bbf :
(H_{1},\mycal{E}_{B,1},\iso_{1}) \to
(H_{2},\mycal{E}_{B,2},\iso_{2})$ of \nc-Hodge
structures is a pair \linebreak $\bbf = (f,f_{B})$, where $f : H_{1} \to
H_{2}$, is an algebraic map of vector bundles which
intertwines the connections, and $f_{B} : \mycal{E}_{B,1} \to
\mycal{E}_{B,2}$ is a map of $\mathbb{Q}$-local systems, such that
$f \circ \iso_{1} \linebreak = \iso_{2}\circ (f_{B}\otimes
\op{id}_{\mathcal{O}})$.  
We will write $\QncHS$ for the category of pure
\nc-Hodge structures. 
\end{defi} 

\

\begin{rem} The meromorphic connection $(M,\nabla)$ where
  $M=H\otimes_{\mathbb{C}[u]}\mathbb{C}[u,u^{-1}]$ can be thought of
  as the de Rham data of the \nc-Hodge structure, the local system
  $\stS_{B}$ of rational vector spaces over $\bS^{1}$ endowed with the
  rational Stokes filtration (see $\mathbb{Q}$-structure axiom) can be
  thought of as the Betti data, and the holomorphic extension $H$ of
  $M$ can be thought of as the analogue of the Hodge filtration.
\end{rem}

\

\bigskip

\subsubsection {\bfseries Variations of \nc-Hodge structures}
\label{sssec:variations} \ One can also
define variations of \nc-Hodge structures:

\begin{defi} \label{defi:variations} Let $S$ be a complex manifold. A
{\bfseries variation  of pure \nc-Hodge structures} over $S$ is a triple
$(H,\mycal{E}_{B},\iso)$, where

\begin{itemize}
\item  $H$ is a holomorphic  $\mathbb{Z}/2$-graded vector bundle
  on $\mathbb{A}^{1} 
  \times S$ which is algebraic in the $\mathbb{A}^{1}$-direction.  
\item $\mycal{E}_{B}$ is a local system of $\mathbb{Z}/2$-graded
  $\mathbb{Q}$-vector spaces on $(\mathbb{A}^{1}-\{0\})\times S$.
\item $\iso$ is an analytic isomorphism of holomorphic vector bundles 
\[
\iso : \mycal{E}_{B}\otimes \mathcal{O}_{(\mathbb{A}^{1}-\{0\})\times
  S} \stackrel{\cong}{\to} H_{|(\mathbb{A}^{1}-\{0\})\times
  S}.
\]
\end{itemize}
Let $\nabla$ be the induced meromorphic connection on $H$.
The data $(H,\mycal{E}_{B},\iso)$ should satisfy: 
\begin{description}
\item[{\sf\bfseries (\nc-filtration axiom)}] The connection $\nabla$
  has a regular singularity along $\{\infty\}\times S$ and
  Poincar\'{e} rank $\leq 1$ along $\{0\}\times S$,  i.e.  
\[
u^2\cdot \nabla_{\frac{\partial}{\partial u}} :  H  \to  H
\] 
is a holomorphic differential operator on $H$ of order $\leq 1$.
\item[{\sf\bfseries (Griffiths transversality axiom)}] For every
  locally defined vector  field 
$\xi \in T_{S}$   we have that 
\[
u\cdot \nabla_{\xi}: H \to H,
\]
is a holomorphic differential operator on $H$ of order $\leq 1$.
\item[{\sf\bfseries ($\mathbb{Q}$-structure axiom)}] The
Stokes structure on the local system $\stS$ on $\bS^{1}\times S$ is
well defined, i.e.  the steps in the Deligne-Malgrange-Stokes filtration on
$\stS$ are sheaves on $\bS^{1}\times S$. Furthermore the
$\mathbb{Q}$-structure $\mycal{E}_{B}$ is compatible with the Stokes
data as in Definition~\ref{defi-nchodge}.
\item[{\sf\bfseries (opposedness axiom)}] The relative version of the
  gluing construction for \nc-Hodge structures gives a globally
  defined complex vector bundle $\widehat{H}$ on $\mathbb{P}^{1}\times
  S$, which is holomorphically trivial in the $\mathbb{P}^{1}$
  direction. Moreover, with respect to the extension $\widehat{H}$ the
  connection $\nabla$ is meromorphic with Poincar\'{e} rank one along
  $(\{0\}\times S)\cup(\{\infty \}\times S)$.
\end{description} 
\end{defi}

\

\bigskip

\subsubsection {\bfseries Relation to other definitions} \ Various special
cases and partial versions of our notion of a \nc-Hodge structure have
been studied before in slightly different but related setups. We list
a few of the relevant notions and references without going into
detailed comparisons:

\begin{itemize}
\item A version of ($\mathbb{Z}$-graded) \nc-Hodge structures 
  appears in the fundamental work of K.Saito (see
  \cite{ksaito-primitive,ksaito-duality,ksaito-survey} and references
  therein) on the Hodge theoretic invariants of quasi-homogeneous
  hypersurface singularities under the name {\em\bfseries weight
    system}.
\item A version of the notion of a variation of (complex) \nc-Hodge
  structure appears in the work of Cecotti-Vafa in Conformal Field
  Theory \cite{cv1,cv2,cv3,bcov} under the name {\em\bfseries
    $tt^{*}$-geometry}.
\item Various versions of the notion of a (complex,polarized)
  \nc-Hodge structure appear in algebraic geometry and non-abelian
  Hodge theory in the works of Simpson \cite{simpson-hf,simpson-mts}
  and T.Mochizuki \cite{tmochizuki-puretwistor,
    tmochizuki-KH,tmochizuki-asymptotic.i,tmochizuki-asymptotic.ii}
  under the names of {\em\bfseries (tame or wild) harmonic bundle} or
  {\em\bfseries pure twistor structure}, and in the work of Sabbah
  \cite{sabbah-twistor} under the name {\em\bfseries integrable pure
    twistor $D$-module}.
\item The analytic germ of a (complex) variation of \nc-Hodge
  structures appears in mirror symmetry in the work of Barannikov
  \cite{b01,b02,b02b} and Barannikov and the second author 
 \cite{bk98} under the name
       {\em\bfseries semi-infinite Hodge structure}. The
       integral structures on semi-infinite Hodge structure were
       recently introduced and studied in the work of Iritani
       \cite{iritani-integral}. 
\item A version of the notion of a (real) \nc-Hodge structure appears
  in singularity theory in the work of Hertling
  \cite{hertling-crelle,hertling-survey} and Hertling-Sevenchek
  \cite{hertling-sevenchek} under the name {\em\bfseries TER
    structure}. Hertling and Sevenchek also consider polarized and
  mixed \nc-Hodge structures. Those appear under the names
  {\em\bfseries TERP structure} and {\em\bfseries mixed TERP
    structure} respectively. In particular in
  \cite{hertling-sevenchek} Hertling and Sevenchek study the
  degenerations of of TERP structures and prove a version of Schmid's
  nilpotent orbit theorem which gives rise to the notion of a limiting
  mixed TERP structure. Degenerations of variants of \nc-Hodge
  structures, as well as limiting mixed \nc-Hodge structures appear
  also in the works of Sabbah \cite{sabbah-degeneration} and S.Szabo
  \cite{szabo-nahm}.
\end{itemize}

\

\bigskip

\subsubsection {\bfseries Relation to usual Hodge
  theory} \label{sssec-usual-hodge} \ Recall (see
e.g. \cite{deligne-h2}) that a {\em\bfseries pure rational Hodge structure}
of weight $w$ is a triple 
$(V,F^{\bullet}V,V_{\mathbb{Q}})$ where:
\begin{itemize}
\item $V$ is a
complex vector space,
\item $V_{\mathbb{Q}} \subset V$ is a
$\mathbb{Q}$-subspace such that $V =
V_{\mathbb{Q}}\otimes_{\mathbb{Q}} \mathbb{C}$, and 
\item $F^{\bullet}V$ is a {\em\bfseries Hodge filtration} of weight
$w$ on $V$, i.e $F^{\bullet}V$ is a decreasing finite exhaustive
filtration by complex subspaces which satisfies $F^{p}V \oplus
\overline{F^{w+1-p}V} = V$, where the complex conjugation on $V$ is
the one given by the real structure $V_{\mathbb{R}} =
V_{\mathbb{Q}}\otimes \mathbb{R} \subset V$.
\end{itemize}
A {\em\bfseries pure Hodge structure} is a direct sum of pure Hodge
structures of various weights, and a morphism of pure Hodge structures
is a linear map of complex vector spaces, which maps the rational
structures into each other and is strictly compatible with the
filtrations. We will write $\QHS$ for the category of pure rational
Hodge structures. It is well known \cite{deligne-h2} that $\QHS$ is an
abelian $\mathbb{Q}$-linear tensor category. For every $w \in
\mathbb{Z}$ we have a $\otimes$-invertible object in
$\QHS$ of pure weight $2w$: the {\em\bfseries Tate Hodge structure
  $\mathbb{Q}(w)$}   given by $\mathbb{Q}(w) :=
(\mathbb{C},F^{\bullet},\mathbb{Q})$, where $F^{i} = \mathbb{C}$ for
$i \leq w$ and $F^{i} = \{0\}$ for $i > w$.

It turns out that pure Hodge structures can be viewed as \nc-Hodge structures. 
This is achieved through a version of the Rees module construction
(see e.g. \cite{simpson-hf})  which converts a filtered vector space into a
bundle over the affine line $\mathbb{A}^{1}$. Specifically, given a pure Hodge
structure $(V,F^{\bullet}V,V_{\mathbb{Q}})$ of weight $w$ we consider
the rank one meromorphic bundle with connection 
\[
\bT_{\frac{w}{2}} := \left(\mathbb{C}\{u\}[u^{-1}], d - \frac{w}{2}\cdot
\frac{du}{u} \right)
\]
and we set
\begin{itemize}
\item $\mycal{H} := \mycal{H}^{w \mod 2}:= \sum_{i} u^{-i}F^{i}V\{u\}$
  viewed as a $\mathbb{C}\{u\}$-submodule in
  $\mathbb{C}\{u\}[u^{-1}]\otimes_{\mathbb{C}}V$. Clearly, this submodule is
  preserved by the operator $\nabla_{u\frac{d}{du}}$ for the 
connection $\nabla := \left(d - \frac{w}{2}\cdot
  \frac{du}{u}\right)\otimes \op{id}_{V}$,
  i.e. $\left(\mycal{H},\nabla\right)$ is a logarithmic holomorphic
  extension of the meromorphic bundle with connection
  $\bT_{\frac{w}{2}}\otimes_{\mathbb{C}}V$.

\

\noindent
{\bfseries Note:} Consider the algebraization $(H,\nabla) =
\mathfrak{B}_{0}\left(\mycal{H},\nabla\right)$ 
of $\left(\mycal{H},\nabla\right)$. The fiber
$H_{1} := \linebreak 
H/(u-1)H$ of $H$ at $1 \in
\mathbb{A}^{1}$ is canonically identified with $V$. By definition the
connection $\nabla$ on $H$ has monodromy
$(-1)^{w}\op{id}_{V}$ and so preserves any rational subspace in
$V$.

\item $\mycal{E}_{B} := \mycal{E}_{B}^{w \mod 2}$ - the
  $\mathbb{Q}$-local system on 
  $\mathbb{A}^{1} - \{0\}$ defined as the subsheaf 
$\mycal{E}_{B} \subset H$
consisting of
    sections whose value at $1$ is in $V_{\mathbb{Q}} \subset V =
    H/(u-1)H$. In other words $\mycal{E}_{B}$ is
    the locally constant sheaf on $\mathbb{A}^{1}-\{0\}$ with
  fiber $V_{\mathbb{Q}}$ and monodromy $(-1)^{w}\op{id}_{V_{\mathbb{Q}}}$.
\item $\iso$ is the isomorphism of complex local systems,
  corresponding to the embedding $\mycal{E}_{B} \subset H$. 
\end{itemize}

\

\begin{rem} \label{rem-bases} On every simply connected open (in the
analytic topology) subset \linebreak 
$U \subset \mathbb{A}^{1}-\{0\}$  the bundle with
  connection  $\bT_{\frac{w}{2}}$ has a horizontal section
  $u^{w/2}$. In particular on such opens we have $H_{|U} =
  \sum_{i} u^{w/2}u^{-i}F^{i}[u]$. 
\end{rem}

\

\noindent
The data $\left(H,\mycal{E}_{B},\iso\right)$ satisfy tautologically
the {\sf\bfseries ($\mathbb{Q}$-structure axiom)} and the
{\sf\bfseries (opposedness axiom)} from
Definition~\ref{defi-nchodge}. Indeed, the {\sf\bfseries
($\mathbb{Q}$-structure axiom)} is satisfied since by definition
$\nabla$ has a regular singularity at $0$ and so $\stS_{\leq \omega} =
\stS$ or $0$ for all $\omega$. The {\sf\bfseries (opposedness axiom)}
is satisfied as it is equivalent in the case of regular singularities
to the oposedness property in the definition of the usual Hodge
structures.

Thus, the assignment $(V,F^{\bullet}V,V_{\mathbb{Q}}) \to
\left(H,\mycal{E}_{B},\iso \right)$ 
gives a functor 
\[
\mathfrak{n} : \QHS \to \QncHS
\] 
which by definition
factors through the orbit category (see e.g. \cite{keller-orbit} for
the definition of an orbit category)
\[
\pi : \QHS \to \QHS/(\bullet
\otimes \mathbb{Q}(1)),
\] 
i.e we have $\mathfrak{N} = \mathfrak{n}\circ
\pi$ for a functor 
\[
\mathfrak{N} :
\QHS/(\bullet \otimes \mathbb{Q}(1)) \to \QncHS.
\]

\

\noindent
The proof of the following statement is an immediate consequence from the
definition.

\

\begin{lemma} \label{lemma-fully-faithful} The functor $\mathfrak{N}$ is fully
  faithful and  its essential image consists of all \nc-Hodge
  structures that have regular singularities and monodromy $= \op{id}$ on
  $H^{0}$ and $= - \op{id}$ on $H^{1}$. 
\end{lemma} 

\

\begin{rem} \label{rem-fully-faithful} It is 
  straightforward to check that the functor $\mathfrak{N}$ can also be
  defined in families  and embeds the category
  of variations of Hodge structures (modulo the Tate twist) into the
  category of variations of \nc-Hodge structures.
\end{rem}

\

\bigskip

\subsubsection {\bfseries \nc-Hodge
  structures of exponential type} \label{sss:irregular} \ As we saw
in section \ref{sssec-usual-hodge} the usual Hodge structures give
rise to special \nc-Hodge structures with regular singularities. The
\nc-Hodge structures with regular singularities are also important
because they can serve as building blocks of general \nc-Hodge
structures. Let $(H,\mycal{E}_{B},\iso)$ be a \nc-Hodge structure, let
$(\mycal{H},\nabla) = \mathfrak{G}_{0}((H,\nabla))$ be the germ of
$(H,\nabla)$ at zero, and assume that $A_{-2} \neq 0$, i.e. $\nabla$
has a second order pole. According to Turrittin-Levelt formal
decomposition theorem (see e.g. \cite{malgrange-formal},
\cite{babbitt-varadarajan-memoir}, \cite[II.5.7 and
II.5.9]{sabbah-frobenius}) we can find a finite base change $p_{N} :
\mathbb{C} \to \mathbb{C}$, $p_{N}(t) := t^{N} = u$, so that
$p_{N}^{*}(\mycal{H},\nabla)[t^{-1}]$ is formally isomorphic to a
direct sum of regular singular connections on meromorphic bundles
multiplied by exponents of Laurent polynomials. More precisely we can
find polynomial tails $g_{i}(t) \in \mathbb{C}[t^{-1}]$,
$\mathbb{C}\{t\}[t^{-1}]$-vector spaces $\mycal{R}_{i}$ and
meromorphic connections
\[
(\nabla_{i})_{\frac{d}{dt}} :
\mycal{R}_{i} \to \mycal{R}_{i},
\] 
each with at most regular
singularity at $0$, so that we have an isomorphism of formal
meromorphic connections over
$\mathbb{C}((t))$:
\[
\Psi : p_{N}^{*}(\mycal{H},\nabla)\bigotimes_{\mathbb{C}\{t\}[t^{-1}]} 
\mathbb{C}((t)) \stackrel{\cong}{\longrightarrow} 
\left(\bigoplus_{i=1}^{m} \be^{g_{i}} \bigotimes_{\mathbb{C}\{t\}[t^{-1}]}
  (\mycal{R}_{i},\nabla_{i})\right) \bigotimes_{\mathbb{C}\{t\}[t^{-1}]} 
\mathbb{C}((t)).
\]
Here $\be^{f}$ denotes the rank one holomorphic bundle with
meromorphic connection \linebreak 
$\left( \mathbb{C}\{t\}, d - df\right)$, and
$(\mycal{R}_{i},\nabla_{i})$ denote meromorphic bundles with
connections having regular singularities.

\begin{rem} \label{rem-exp-basis} The bundle $\be^{f}$
has a non-vanishing horizontal section, namely $e^{f}$. In particular the
multivalued flat sections of $\be^{g_{i}}\otimes
(\mycal{R}_{i},\nabla_{i})$ are given by multiplying multivalued flat
sections of $(\mycal{R}_{i},\nabla_{i})$ by $e^{g_{i}}$.
\end{rem} 

\
\noindent
In the examples coming from Mirror Symmetry that we are interested in,
the base change $p_{N}$ is not needed for the decomposition to
work. In this case we can take $g_{i}(u) = \bc_{i}/u$ where 
$\bc_{1}, \ldots, \bc_{m} \in \mathbb{C}$ denote the  distinct
eigenvalues of $A_{-2}$. Because of this we introduce  the following
definition (see also \cite[Definition~8.1]{hertling-sevenchek}):

\begin{defi} \label{defi-no-base-change} We say that a \nc-Hodge
  structure $(H,\mycal{E}_{B},\iso)$ is of {\bfseries
   exponential type} if there exists a formal
  isomorphism
\[
\Psi : (\mycal{H}\otimes_{\mathbb{C}\{u\}} \mathbb{C}[[u]],\nabla)
\stackrel{\cong}{\to} 
\bigoplus_{i=1}^{m} \left(\be^{\bc_{i}/u}\otimes
(\mycal{R}_{i},\nabla_{i})\right)\otimes_{\mathbb{C}\{u\}}
\mathbb{C}[[u]] 
\]
where $(\mycal{R}_{i},\nabla_{i})$ are meromorphic bundles with
connections with regular singularities and $\bc_{1}, \ldots,
\bc_{m} \in \mathbb{C}$ denote the  distinct 
eigenvalues of $A_{-2}$.
\end{defi}

\

\begin{rem} \label{rem-decompose-with-lattices}
$\bullet$ \ There are various sufficient conditions that will
  guarantee that a given \nc-Hodge structure is decomposable without
  base change. For instance, this will be the case if $A_{-2}$ has
  distinct eigenvalues, or if $A_{-1} = 0$. More generally, it
  suffices to require that we can find holomorphic functions
  $\ell_{i}(u) \in \mathbb{C}\{u\}$ so that $\ell_{i}(0) =
  \bc_{i}$ for $i = 1, \ldots, m$ and the characteristic
  polynomial of $u^{2}A(u)$ is $\det\left(\bc\cdot \op{id} -
  u^{2}A(u)\right) = \prod_{i = 1}^{m}(\bc - 
  \ell_{i})^{\nu_{i}}$. 

\

\noindent
$\bullet$ \ Not every irregular connection with a pole of order two is
  of exponential type. Indeed the rank two connection
\[
\nabla = d - \begin{pmatrix} 0 & u^{-2} \\ u^{-1} &
  \frac{u^{-1}}{2}\end{pmatrix}
\]
has a horizontal section 
\[
\begin{pmatrix} e^{-2 u^{-\frac{1}{2}}}
  \\ u^{\frac{1}{2}}e^{-2 u^{-\frac{1}{2}}}\end{pmatrix}, 
\]
and so one needs a quadratic base change for the formal decomposition
to work for this connection.

\

\noindent
$\bullet$ \ If a \nc-Hodge structure $(H,\mycal{E}_{B},\iso)$ is of
exponential type, then one can check (see
\cite[Lemma~8.2]{hertling-sevenchek}) that for each $i = 1, \ldots, m$
we can find a unique holomorphic extension
$\mycal{H}_{\bc_{i}} \subset \mathcal{R}_{i}$ in which the connection
has a second order pole and so that $\Psi$ induces a
formal isomorphism of holomorphic bundles with meromorphic connections
\[
\Psi : (\mycal{H},\nabla)\otimes \mathbb{C}[[u]]
\stackrel{\cong}{\longrightarrow} 
\left(\bigoplus_{i=1}^{m} \be^{\bc_{i}/u}\bigotimes_{\mathbb{C}\{u\}}
(\mycal{H}_{\bc_{i}},\nabla_{i})\right)\otimes \mathbb{C}[[u]], 
\]
over $\mathbb{C}[[u]]$.
\end{rem}

\

\noindent
The \nc-Hodge structures with regular singularities or the \nc-Hodge
structures of exponential type comprise full subcategories 
\[
\QncHS^{\op{reg}} \subset \QncHS^{\exp} \subset \QncHS
\]
in $\QncHS$. In fact, in the exponential type case one can state the
\nc-Hodge structure axioms in an easier way. The simplification comes
from the fact that in this case the Deligne-Malgrange-Stokes
filtration is given by subsheaves $\stS_{\le \lambda}$ of $\stS$ that
are labeled by $\lambda\in \mathbb{R}$ and consisting of solutions
decaying faster than
$\bO\left(\exp\left(\frac{\lambda+o(1)}{r}\right)\right)$, $r=|u|$.
Indeed, tracing through the definition one sees that in the
exponential case for a ray defined by $\varphi$ the jumps of the steps
of the Deligne-Malgrange-Stokes filtration occur exactly at the
numbers $\op{Re}(\bc_{i}e^{-i\varphi})$.  Furthermore, the associated
graded pieces for the filtration are local systems on the circle and
in fact coincide with the regular pieces $(\mycal{R}_{i},\nabla_{i})$
that appear in the formal decomposition of the connection. Hence one
arrives at the following

\begin{defi} \label{defi:exp.type.axioms}
A  {\bfseries rational pure
    \nc-Hodge structure of exponential type} 
  consists of the data
  $(H,\mycal{E}_{B},\iso)$, where 
\begin{itemize}
\item $H$ is a  $\mathbb{Z}/2$-graded
algebraic vector bundle on $\mathbb{A}^{1}$. 
\item $\mycal{E}_{B}$ is a local system of finite dimensional
  $\mathbb{Z}/2$-graded $\mathbb{Q}$-vector spaces on \linebreak
  $\mathbb{A}^{1}-\{0\}$. 
\item $\iso$ is an {\em analytic}  isomorphism  of holomorphic vector
  bundles on  
$\mathbb{A}^{1} - \{ 0 \}$: 
\[
\iso : \mycal{E}_{B}\otimes \mathcal{O}_{\mathbb{A}^{1} - \{ 0 \}}
  \stackrel{\cong}{\longrightarrow}
  H_{|\mathbb{A}^{1} - \{ 0 \}}. 
\]
\end{itemize}

\

\noindent
These data have to satisfy the following axioms:
\begin{description}
\item[{\sf\bfseries (\nc-filtration axiom)$^{\bexp}$}] The connection
  $\nabla$ induced from $\iso$ is a meromorphic
  connection of exponential type on $H$ with a pole of order 
  $\leq 2$ at $u 
  = 0$ and a regular singularity at $\infty$. More precisely, there exist:
\begin{itemize}
\item  a  holomorphic frame of $\mycal{H}$ near $u = 0$ in which 
\[
\nabla = d + \left(\sum_{k \geq -2} A_{k} u^{k}\right)du
\]
with $A_{k} \in \op{Mat}_{r\times r}(\mathbb{C})$, $r =
\op{rank}_{\mathbb{C}\{u\}} \mycal{H}$.
\item a  holomorphic frame of $\mycal{H}$ near $u = \infty$ in which 
\[
\nabla = d + \left(\sum_{k \geq -1} B_{k} u^{-k}\right)d(u^{-1})
\]
and $B_{k} \in \op{Mat}_{r\times r}(\mathbb{C})$.
\item a formal isomorphism over $\mathbb{C}((u))$:
\[
(\mycal{H}[u^{-1}],\nabla) \stackrel{\cong}{\to}
\bigoplus_{i=1}^{m} \be^{\bc_{i}/u}\otimes
(\mycal{R}_{i},\nabla_{i})
\]
where $(\mycal{R}_{i},\nabla_{i})$ are meromorphic bundles with
connections with regular singularities and $\bc_{1}, \ldots,
\bc_{m} \in \mathbb{C}$ denote the  distinct 
eigenvalues of $A_{-2}$.
\end{itemize}
\item[{\sf\bfseries ($\mathbb{Q}$-structure axiom)$^{\bexp}$}] The
  $\mathbb{Q}$-structure $\mycal{E}_{B}$ on $(H,\nabla)$ is compatible
  with Stokes data in the following sense.  The
  filtration $\{ \stS_{\leq \lambda} \}_{\lambda \in \mathbb{R}}$ of
  $\stS$ by the subsheaves  $\stS_{\leq \lambda}$, whose stalk at
  $\varphi \in \bS^{1}$ is given by 
\[
\left(\stS_{\leq \lambda}\right)_{\varphi} := \left\{ s \in
\stS_{\varphi} = 
\Gamma\left(\mathbb{R}^{\times}_{+}e^{i\varphi},H\right) \;\;
\left| \;\; 
\text{\begin{minipage}[c]{2.5in}
$s$ is a $\nabla$-horizontal section of $H$ over the ray 
$\mathbb{R}^{\times}_{+}e^{i \varphi}$, for which 
\[
\left\|s\left(r e^{i \varphi}\right)\right\|
= \bO\left(\exp\left( \frac{\lambda + \bo(1)}{r}\right)\right).
\]
when $r \to 0$.
\end{minipage}
}\right. \right\}
\]
is defined over $\mathbb{Q}$, i.e.
\[
\left(\stS_{\leq \lambda} \cap \stS_{B}\right) \otimes_{\mathbb{Q}}
\mathbb{C} = \stS_{\leq \lambda}
\]
for all  $\lambda \in \mathbb{R}$.
\item[{\sf\bfseries (opposedness axiom)$^{\bexp}$ = (opposedness
  axiom)}] 
\end{description}
\end{defi}

\

\medskip

\begin{rem} \label{rem-switch.steps} It is instructive to understand
  more explicitly the behavior of the Deligne-Malgrange-Stokes
  filtration for \nc-Hodge structures (or more generally  irregular
  connections)  of exponential type. As before we denote by  $\stS$
  the  complex local system on the circle
  $\bS^{1}$ corresponding to a \nc-Hodge structure for which $A_{-2}$
  has distinct eigenvalues $\boldsymbol{c}_{1}, \ldots
  \boldsymbol{c}_{m}$. 

By definition, for every $\varphi$, the steps in the
  Deligne-Malgrange-Stokes filtration $(\stS_{\leq
  \lambda})_{\varphi}$ jump exactly when $\lambda$ crosses one of the
  numbers $\op{Re}(\bc_{k}e^{-i\varphi})$. More invariantly, the
  assignment $\varphi \in \bS^{1} \mapsto \{
  \op{Re}(\bc_{1}e^{-i\varphi}), \ldots, \op{Re}(\bc_{k}e^{-i\varphi}) \}
  \subset \mathbb{R}$ is a sheaf $\Lambda$ of finite sets of real
  numbers (possibly
  with repetitions) on $\bS^{1}$. For a general value of $\varphi$,
  the real numbers $\{ \op{Re}(\bc_{1}e^{-i\varphi}), \ldots,
  \op{Re}(\bc_{k}e^{-i\varphi}) \}$ are all distinct but for finitely
  many special values of $\varphi$ some of
  $\op{Re}(\bc_{1}e^{-i\varphi}), \ldots, \op{Re}(\bc_{k}e^{-i\varphi})$
  will coalesce. More precisely we have the  Stokes rays
  $\mathbb{R}_{>0} \cdot i(\bc_{b} - \bc_{a})$ and the associated set
  $\text{\sf\bfseries SD} \subset [0,2\pi)$ of  Stokes
  directions: i.e. $\varphi \in \text{\sf\bfseries SD} $, if and only
  if there is some pair $a \neq b$ s.t.  $\bc_{a} - \bc_{b} =
  re^{i\left(\frac{\pi}{2} + \varphi\right)}$ for some $r >
  0$. Clearly for every open arc $U \subset \bS^{1}$, which does not
  intersect $\text{\sf\bfseries SD}$ the restriction $\Lambda_{|U}$
  is a local system of finite sets of cardinality $m$. Moreover the
  values $\varphi \in \text{\sf\bfseries SD}$ are precisely the ones
  for which some of $\op{Re}(\bc_{1}e^{-i\varphi}), \ldots,
  \op{Re}(\bc_{k}e^{-i\varphi})$ become equal to each other. 

Now recall that for any given $\varphi \in \bS^{1}$, the subspaces
  $(\stS_{\leq \lambda})_{\varphi} \subset \stS_{\varphi}$ do not
  change if we move $\lambda \in \mathbb{R}$ continuously without
  passing through some element of $\Lambda_{\varphi}$. In other words,
  we can label the steps of the Deligne-Malgrange-Stokes filtration by
  local sections of $\Lambda$, and so that at each $\varphi \in
  \bS^{1}$ the steps are ordered according to the order on
  $\Lambda_{\varphi}$ induced from the embedding $\Lambda_{\varphi}
  \subset \mathbb{R}$. The finite set $\text{\sf\bfseries SD} \subset
  \bS^{1}$ of  Stokes directions breaks the circle into disjoint
  arcs. Over each such arc $U$ we have that $\Lambda_{|U}$ is a local
  system of finite sets of real numbers with $m$ linearly ordered flat
  sections and the steps Deligne-Malgrange-Stokes filtration of
  $\stS_{| U}$ are labeled naturally by these sections. If we move
  from $U$ to an adjacent arc $U'$ by passing across a Stokes
  direction $\phi \in \text{\sf\bfseries SD}$, then some of the
  elements in the labelling set get identified at $\phi$ and get
  reordered when we cross over to $U'$ (see
  Figure~\ref{fig:labelsLambda}).

\

\smallskip

\

\begin{figure}[!ht]
\begin{center}
\psfrag{l1}[c][c][1][0]{{$\lambda_{1}$}}
\psfrag{l2}[c][c][1][0]{{$\lambda_{2}$}}
\psfrag{l3}[c][c][1][0]{{$\lambda_{3}$}}
\psfrag{l4}[c][c][1][0]{{$\lambda_{4}$}}
\psfrag{l1'}[c][c][1][0]{{$\lambda_{1}'$}}
\psfrag{l2'}[c][c][1][0]{{$\lambda_{2}'$}}
\psfrag{l3'}[c][c][1][0]{{$\lambda_{3}'$}}
\psfrag{l4'}[c][c][1][0]{{$\lambda_{4}'$}}
\psfrag{U}[c][c][1][0]{{$U$}}
\psfrag{U'}[c][c][1][0]{{$U'$}}
\psfrag{S1}[c][c][1][0]{{$\bS^{1}$}}
\psfrag{R}[c][c][1][0]{{$\mathbb{R}$}}
\psfrag{aaaa}[c][c][1][0]{\begin{minipage}[c]{1in}
 Stokes \\
direction
\end{minipage}}
\psfrag{fi}[c][c][1][0]{{$\phi$}}
\epsfig{file=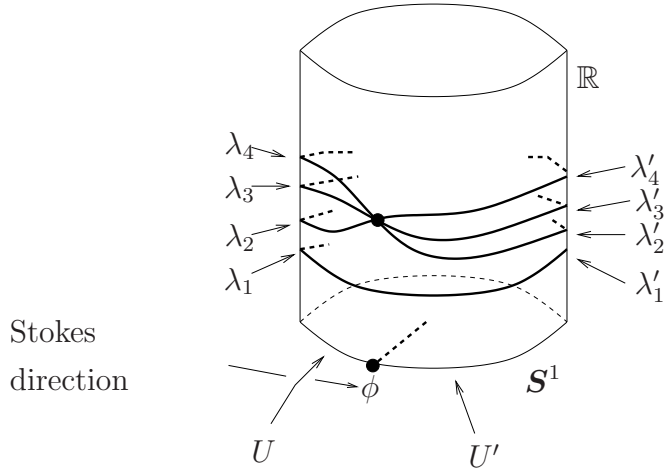,width=3in} 
\end{center}
\caption{The system of labels for the Deligne-Malgrange-Stokes
  filtration.}\label{fig:labelsLambda}  
\end{figure}

\

\smallskip

\noindent
In fact, if $\lambda_{1} < \ldots < \lambda_{m}$ are the ordered flat
sections of $\Lambda_{|U}$, and $\lambda_{1}' < \ldots < \lambda_{m}'$
are the ordered flat sections of $\Lambda_{|U'}$, then the transition
from the $\lambda$'s to the $\lambda'$'s is always such that certain
groups of consecutive $\lambda$'s are totally reordered into groups of
consecutive $\lambda'$'s. For instance in
Figure~\ref{fig:labelsLambda} the passage from $\{\lambda_{1},
\lambda_{2}, \lambda_{3}, \lambda_{4} \}$ to $\{\lambda_{1}',
\lambda_{2}', \lambda_{3}', \lambda_{4}' \}$ across the Stokes
point $\phi \in \text{\sf\bfseries SD}$ has the effect of
relabelling:
$\lambda_{1} \mapsto \lambda_{1}'$, $\lambda_{2} \mapsto
\lambda_{4}'$, $\lambda_{3} \mapsto \lambda_{3}'$, and $\lambda_{4}
\mapsto \lambda_{2}'$. 
\end{rem}

\

\noindent
This behavior of the labelling set and the behavior of the associated
filtration can be systematized in the following:

\begin{defi} \label{defi:abstractDMS} Let $\stS$ be a finite
  dimensional local system of $\mathbb{Z}/2$-graded complex vector
  spaces over $\bS^{1}$. Let $\bc_{1}, \ldots, \bc_{m}$ be
  distinct complex numbers, let $\Lambda$ be the sheaf of finite sets
  of real numbers on $\bS^{1}$ given by $\varphi \mapsto \{
  \op{Re}(\bc_{1}e^{-i\varphi}), \ldots,
  \op{Re}(\bc_{1}e^{-i\varphi})\}$, and let \linebreak $\text{\sf\bfseries SD}
  \subset \bS^{1}$ be the associated set of  Stokes
  directions. 

An {\bfseries abstract Deligne-Malgrange-Stokes filtration of $\stS$
  of exponential type and exponents $(\bc_{1}, \ldots, \bc_{m})$} is a
  filtration by subsheaves $\stS_{\leq\lambda}$ such that:
\begin{itemize}
\item  $\stS_{\leq\lambda}$ is labeled by local
  continuous sections $\lambda$ of $\Lambda$ and is locally constant
  on any arc which does not intersect $\text{\sf\bfseries SD}$.
\item  Suppose $\varphi
  \in \text{\sf\bfseries SD}$, and let $U, U' \subset \bS^{1}-
  \text{\sf\bfseries SD}$ be the two arcs adjacent at $\varphi$. Let 
$\lambda_{1} < \cdots < \lambda_{m}$ and $\lambda_{1}' < \cdots <
  \lambda_{m}'$ be the ordered flat sections of $\Lambda_{|U}$ and
  $\Lambda_{|U'}$ respectively. Trivialize $\stS$ on $U\cup U'\cup
  \{\varphi \}$ by
  identifying the flat sections with the elements of the fiber
  $\stS_{\varphi}$ and let 
\[
0 \subset F_{\leq \lambda_{1}} \subset \ldots \subset  F_{\leq
  \lambda_{m}} \subset \stS_{\varphi}, \qquad \text{and} \qquad  
0 \subset F'_{\leq \lambda_{1}'} \subset \ldots \subset  F'_{\leq
  \lambda_{m}'} \subset \stS_{\varphi}
\]
be the filtrations corresponding to this trivialization and the
 filtrations $\stS_{\leq \lambda}$ on $U$ and 
 $U'$ respectively, i.e. 
\[
F_{\leq \lambda_{i}}  := \lim_{\substack{\psi \in U \\ \psi \to
    \varphi}} \left( \stS_{\leq \lambda_{i}} \right)_{\psi} \qquad
    \text{and} \qquad 
F'_{\leq \lambda_{i}'}  := \lim_{\substack{\psi \in U' \\ \psi \to
    \varphi}} \left( \stS_{\leq \lambda_{i}'} \right)_{\psi}
\]
 Let 
$1 \leq i_{1} < j_{1} \leq i_{2} < j_{2} \leq \cdots \leq i_{k} <
 j_{k}  \leq m$ 
be the sequence of integers such that $\lambda_{a} = \lambda_{a}'$ for 
$a \not\in [i_{1},j_{1}]\cup [i_{2},j_{2}]\cup \cdots \cup
[i_{k},j_{k}]$, and for each interval $[i_{s}, j_{s}]$ we have
that $\lambda_{j_{s}}' = \lambda_{i_{s}}$, $\lambda_{j_{s} - 1} =
\lambda_{i_{s}+1}$, \ldots $\lambda_{i_{s}}' = \lambda_{j_{s}}$. Then
we require that:
\begin{itemize}
\item  for each $a \not\in [i_{1},j_{1}]\cup [i_{2},j_{2}]\cup \cdots \cup
[i_{k},j_{k}]$ we have $F_{\leq \lambda_{a}} = F'_{\lambda_{a}'}$;
\item for each $s = 1, \ldots, k$, $F_{\leq \lambda_{j_{s}}} = F_{\leq
  \lambda_{j_{s}}'}$ and the filtrations
\[
\xymatrix@C-1.5pc@R-1pc{
F_{\leq \lambda_{i_{s}}}/F_{\leq \lambda_{i_{s}-1}}  & \subset &  
F_{\leq \lambda_{i_{s}+1}}/F_{\leq \lambda_{i_{s}-1}} & \subset &
\cdots & \subset &
F_{\leq \lambda_{j_{s}}}/F_{\leq \lambda_{i_{s}-1}} \ar@{=}[d] \\
F'_{\leq \lambda_{i_{s}}'}/F'_{\leq \lambda_{i_{s}-1}'} & \subset &  
F'_{\leq \lambda_{i_{s}+1}'}/F'_{\leq \lambda_{i_{s}-1}'} & \subset &
\cdots & \subset &
F'_{\leq \lambda_{j_{s}}'}/F'_{\leq \lambda_{i_{s}-1}'}
}
\]
are $(j_{s} - i_{s})$-opposed. 
\end{itemize}
\end{itemize}
\end{defi}

\

\begin{rem} The above discussion generalizes immediately from 
  connections of exponential type to arbitrary meromorphic connections
 (see remark \ref{rem-Stokes-property}).  One gets a collection of
 curves drawn on the boundary of the cylinder which can be interpreted
 as a projection to $0$-jets of a Legendrian link in the contact
 manifold of $1$-jets of functions on $\bS^{1}$.
\end{rem}

\

\noindent
The categories of \nc-Hodge structures, of \nc-Hodge structures of
exponential type, or of \nc-Hodge structures with regular
singularities all behave similarly to ordinary Hodge structures. For
instance one can introduce the notion of polarization on \nc-Hodge
structures, which specializes to the usual notion in the case of
ordinary Hodge structures. (This will not be needed for our discussion
so we will not spell it out here. The interested reader may wish to
consult \cite{hertling-survey,hertling-sevenchek,kontsevich-lefschetz}
for the details of the definition.) In
fact we have the following

\begin{lemma} \label{lemma:abelian}
The categories 
$\QncHS^{\op{reg}} \subset \QncHS^{\exp} \subset \QncHS$
are $\mathbb{Q}$-linear abelian categories.
The respective categories of polarizable \nc-Hodge structures are
semi-simple. 
\end{lemma}
{\bfseries Proof:} The statement is a manifestation of Simpson's
{\bfseries Meta-Theorem} from \cite{simpson-mts}. The opposedness axiom
implies that the respective categories are abelian and the existence
of polarizations implies the semi-simplicity. The proofs follow
verbatim the argument in usual Hodge theory or the argument in
\cite{simpson-mts}. Alternatively one can use the comparison statement
\cite[Lemma~3.9]{hertling-sevenchek} identifying the \nc-Hodge
structures with pure twistor structures and then invoke
\cite[Lemma~1.3 and Lemma~3.1]{simpson-mts}. \ \hfill $\Box$

\

\noindent
The bundles with connections $(\mycal{H}_{\bc_{i}},\nabla_{i})$ can be
thought of as the regular singular constituents of the \nc-Hodge
structure $(H,\mycal{E}_{B},\iso)$. The
$(\mycal{H}_{\bc_{i}},\nabla_{i})$'s are invariants of the \nc-Hodge
structure but of course they do not give a complete set of invariants
(see the third point in \ref{rem-decompose-with-lattices}). As usual
we need additional Stokes data (see e.g. \cite{sabbah-frobenius}) in
order to reconstruct the pair $(\mycal{H},\nabla)$ from its regular
constituents.  To understand how the rest of \nc-Hodge structure
arises from the constituents we need to understand how the rational
structure $\mycal{E}_{B}$ interacts with the Stokes data. This process
is very similar to the interaction between Betti, de Rham and
Dolbeault cohomology in ordinary Hodge theory and we will describe it
in detail in section~\ref{sec-gluing}.

\

\noindent
The \nc-Hodge structures one finds in geometric examples are very
often regular (e.g. in the case of ordinary Hodge structures) or at
worst have exponential type. It is also expected that the \nc-Hodge
structures arising in mirror symmetry will always be of exponential
type but at the moment this is only supported by experimental
evidence. 

We will discuss in detail some of this  evidence in the
subsequent sections. Before we get to the examples  however, it will
be instructive to comment on the reason for introducing the 
\nc-Hodge structures at the first place. The geometric significance of
these structures stems from the fact that they appear naturally on the
cohomology of non-commutative spaces of categorical nature.

\subsection{\bfseries Hodge structures in \nc
  \ geometry} \label{subsec:ncgeometry} The version of non-commutative
geometry that is most relevant to \nc-Hodge structures is the one in
which a proxy for the notion of a non-commutative space (\nc-space) is
a category, thought of as the (unbounded) derived category of
quasi-coherent sheaves on that space.

\

\bigskip

\subsubsection {\bfseries Categorical \nc-geometry} \
\label{sssec:categorical} The basic notion here is:

\begin{defi} \label{defi:ncspace} A {\bfseries graded complex
    \nc-space} (respectively a {\bfseries complex \nc-space}) is a
    $\mathbb{C}$-linear differential graded (respectively
    $\mathbb{Z}/2$-graded) category $C$ which is homotopy complete and
    cocomplete.
\end{defi}

\

\noindent
{\em\bfseries Notation:} We will often write $C_{X}$ for the category to
signify that it describes the sheaf theory of some \nc-space $X$,
even when we do not have a geometric construction of $X$.

\

\noindent
The categorical point of view to non-commutative geometry goes back to
the works of Bondal \cite{bondal-mpi}, Bondal-Orlov
\cite{bondal.orlov-reconstruction,bondal.orlov-icm} with many
non-trivial examples computed in the later works of Orlov
\cite{orlov-LG,orlov-LG.equivalences,orlov-LG2}, Caldero-Keller
\cite{ck-cluster1,ck-cluster2}, Aroux, Orlov, and the first author,
\cite{auroux-katzarkov-orlov,auroux-katzarkov-orlov2}, Kuznetsov
\cite{kuznetsov-hpd,kuznetsov-quadrics,kuznetsov-grassmannian},
etc. More recently this approach to \nc-geometry became the central
part of a long term research program initiated by the second author
and was studied systematically in the works of the second author and
Soibelman \cite{ks-ncgeometry,kontsevich-lefschetz}, To\"{e}n
\cite{toen-morita}, and To\"{e}n-Vaquie \cite{toen-vaquie}.

\

\begin{rem} \label{rem:ncspace} {\bfseries (i)} \ We do not spell out here
  the notions of homotopy completeness and cocompleteness in dg
  categories since on one hand they are quite technical and on the
  other hand will not be used later in the paper. It is worth
  mentioning though that some effort is required to define these
  notions. In the original approach of the second author described in
  his 2005 IAS lectures and in his 2007 course at the University of
  Miami the homotopy completeness and cocompletness in $C$ was defined
  by a universal property for homotopy coherent diagrams of objects in
  the dg category labeled by simplicial sets.  Alternatively
  \cite{toen-letter} one may use the model category
  $(C^{\op{op}}-\text{{\sf mod}})$ of $C^{\op{op}}$-dg modules, whose
  equivalences are the quasi-isomorphisms, and whose fibrations are
  the epimorphisms. In these terms one says that $C$ is homotopy
  complete if the full subcategory of $(C^{\op{op}}-\text{{\sf mod}})$
  consisting of quasi-representable objects is preserved by all small
  homotopy limits (defined via the given model structure). Similarly
  we say that $C$ is homotopy cocomplete if $C^{\op{op}}$ is homotopy
  complete.

\

\noindent
{\bfseries (ii)} \ Note that in the above definition the category $C$ is
automatically triangulated as follows already from the
existence of {\em finite} homotopy limits, and Karoubi closed by the
standard mapping telescope construction \cite{boekstedt-neeman}.
\end{rem}

\

\begin{ex} \label{ex:ncspaces} The two main types  of  \nc-spaces are
  the following: 
\begin{description} 
\item[usual schemes:] Usual complex schemes can be viewed as (graded)
  \nc-spaces. Given a scheme $X$ over $\mathbb{C}$, the corresponding
  category $C_{X}$ is the derived category $D(\text{\sf Qcoh}(X)w)$ of
  quasi-coherent sheaves on $X$ taken with an appropriate dg
  enhancement (see \cite{bondal.kapranov-enhanced}). In particular,
  the closed point $\op{pt} = \op{Spec}(\mathbb{C})$ corresponds to
  the category $C_{\op{pt}}$ of complexes of $\mathbb{C}$-vector spaces.
\item[modules over an algebra:] If $A$ is a differential graded (or
  $\mathbb{Z}/2$-graded) unital associative algebra over $\mathbb{C}$,
  then we get a \nc-space $\ncSpec(A)$ such that $C_{\ncSpec(A)} =
  (A-\text{{\sf mod}})$ is the category of dg modules over $A$ which
  admit an exhaustive increasing filtration whose associated graded
are sums of shifts of $A$.
\end{description}
\end{ex}

\

\noindent
To illustrate how the above notion of a \nc-space fits with the \nc Hodge
structures we will concentrate on the case of \nc-affine spaces,
i.e. \nc-spaces equivalent to  $\ncSpec(A)$ for some differential 
  $\mathbb{Z}/2$-graded algebra $A$ over $\mathbb{C}$. Note that
because of derived Morita equivalences an affine \nc-space $X$ does not
determine an algebra $A$ uniquely, i.e. different algebras can give
rise to the same \nc-space.

\

\begin{rem} \label{rem:affinenc} The condition
  is not as restrictive as it appears at a first glance.  In fact
  almost all \nc-spaces that one encounters in practice are affine. For
  instance usual quasi-compact quasi-separated schemes of finite type
  over $\mathbb{C}$ are affine when viewed as \nc-spaces. This follows
  from a deep theorem of Bondal and van den Bergh
  \cite{bondal-van-den-bergh} which asserts that for such a scheme $X$
  the category $C_{X} = D(\text{\sf Qcoh}(X))$ has a compact generator
  $\mathcal{E}$. That is, we can find an object $\mathcal{E} \in
  C_{X}$ so that
\[
\op{Hom}(\mathcal{E},\bullet) : C_{X} \to C_{\op{pt}}
\]
commutes with homotopy colimits and has a zero kernel. In particular
the dg algebra
computing the category $C_{X}$ is  given in terms of the generator  
$\mathcal{E}$, i.e. 
\[
C_{X} \cong (Hom(\mathcal{E},\mathcal{E})^{\op{op}}-\text{{\sf mod}}).
\]
\end{rem}

\

\noindent
Suppose now that $X = \ncSpec(A)$. Recall that an object
$\mathcal{E} \in C_{X} = (A-\text{{\sf mod}})$ is {\em\bfseries perfect} if
$\op{Hom}(\mathcal{E},\bullet)$ preserves small homotopy colimits. We
will write $\Perf_{X}$ for the full subcategory of perfect objects
in $C_{X}$. We now have the following definition (see
e.g. \cite{ks-ncgeometry,kontsevich-lefschetz} or
\cite{toen-vaquie}): 

\

\begin{defi} \label{defi:smooth.compact} A complex differential
  $\mathbb{Z}/2$-graded algebra is called
\begin{description}
\item[{\em\bfseries smooth:}] if $A \in \Perf_{\ncSpec(A\otimes
  A^{\op{op}})}$;
\item[{\em\bfseries compact:}] if $\dim_{\mathbb{C}} H^{\bullet}(A,d_{A})
  < + \infty$ or equivalently if $A \in \Perf_{\op{pt}}$. 
\end{description}
\end{defi}

\

\noindent
{\bfseries Note:} \ One can check (see e.g. \cite{ks-ncgeometry} or
\cite{toen-vaquie}) that the properties of $X$ being smooth ad compact
do not depend on the choice of the algebra $A$ which computes
$C_{X}$. Also, for a usual scheme $X$ of finite type over $\mathbb{C}$,
smoothness and compactness in the scheme-theoretic sense are
equivalent to smoothness and compactness in the \nc-sense.

\

\bigskip

\subsubsection {\bfseries The main conjecture} \
The analogy with commutative geometry suggests that one should look
for pure \nc-Hodge structures on the cohomology of smooth an proper
\nc-spaces. More precisely we have the following basic conjecture

\

\begin{conn} \label{con:main.conjecture} Let $X$ be a smooth and
  compact \nc-space over $\mathbb{C}$. Then the periodic cyclic
  homology $HP_{\bullet}(C_{X})$ of $C_{X}$ carries a natural functorial pure
  $\mathbb{Q}$-\nc-Hodge structure with regular singularities. 

Furthermore if the $\mathbb{Z}/2$-grading on $X$ can be refined to a
$\mathbb{Z}$-grading, then the \nc-Hodge structure on
$HP_{\bullet}(C_{X})$ is an ordinary pure Hodge structure,
i.e. belongs to the essential image of the functor $\mathfrak{N}$.
\end{conn}

\

\bigskip

\subsubsection {\bfseries Cyclic homology} \ \label{sssec:cyclic}
There are some natural candidates for the various ingredients of the
conjectural \nc-Hodge structure on $HP_{\bullet}(C_{X})$. Assuming
that $X \cong \ncSpec(A)$ is \nc-affine, we can compute
$HP_{\bullet}(C_{X})$ in terms of $A$. Namely
\[
HP_{\bullet}(C_{X}) = HP_{\bullet}(A) =
HP_{\bullet}\left(C_{\bullet}^{\op{red}}(A,A)((u)), \partial + u\cdot B\right),
\]
where
\begin{itemize}
\item $u$ is an even formal variable (of degree $2$ in the
  $\mathbb{Z}$-graded case);
\item $C_{-k+1}^{\op{red}}(A,A)((u)) := A\otimes (A/\mathbb{C}\cdot
  1_{A})^{\otimes k}\otimes \mathbb{C}((u))$, for all $k \geq 0$;
\item $\partial = b + \delta$, where 
\[
\begin{split}
b(a_{0}\otimes \cdots \otimes a_{n}) & := \sum_{i=0}^{n-1}
(-1)^{\deg(a_{0}\otimes \cdots \otimes a_{i})} a_{0}\otimes \cdots
\otimes a_{i}a_{i+1} \otimes \cdots \otimes a_{n} \\
& \qquad +
(-1)^{\deg(a_{0}\otimes \cdots \otimes a_{n})(\deg(a_{n})+ 1) + 1}
a_{n}a_{0}\otimes \cdots \otimes a_{n-1},
\end{split}
\]
is the Hochschild differential, and 
\[
\delta(a_{0}\otimes \cdots \otimes a_{n}) :=
\sum_{i=0}^{n}(-1)^{\deg(a_{0}\otimes \cdots \otimes a_{i-1})}
a_{0} \otimes \cdots \otimes d_{A}a_{i} \otimes \cdots \otimes a_{n}
\]
is the differential induced from $d_{A}$ via the Leibniz rule;
\item 
\[
\begin{split}
B(a_{0}\otimes & \cdots \otimes a_{n}) := \\
& \sum_{i=0}^{n}
(-1)^{(\deg(a_{0}\otimes \cdots \otimes a_{i}) - 1)(\deg(a_{i+1}\otimes
  \cdots \otimes a_{n}) - 1)} 1_{A}\otimes a_{i+1}\otimes \cdots \otimes
  a_{n} \otimes a_{0} \otimes \cdots \otimes a_{i},
\end{split}
\]
is Connes' cyclic differential.
\end{itemize}

\

\bigskip

\subsubsection [{\bfseries The degeneration conjecture}] 
{\bfseries The degeneration conjecture and the vector bundle
  part of the \nc-Hodge structure} \ \label{sssec:degeneration}
Note that by construction $HP_{\bullet}(C_{X})$ is a module over
$\mathbb{C}((u))$. We can also look at the negative cyclic
homology $HC^{-}_{\bullet}(C_{X})$ of $C_{X}$.  By definition
$HC^{-}_{\bullet}(C_{X})$  is
the cohomology of the complex
\[
\left( C_{\bullet}^{\op{red}}(A,A)[[u]],\partial + u\cdot B\right),
\]
and so is a module over $\mathbb{C}[[u]]$. 
The
specialization 
\[
HC^{-}_{\bullet}(C_{X})/u HC^{-}_{\bullet}(C_{X})
\]
of this module at $u = 0$ maps to
the cohomology of the complex
\[
(C_{\bullet}^{\op{red}}(A,A), \partial)
\]
of reduced Hochschild chains for $A$ which by definition is the Hochschild
homology $HH_{\bullet}(A)$ of $A$.  The Hochschild-to-cyclic spectral
sequence implies that 
\begin{equation} \label{eq:HCinequality}
\dim_{\mathbb{C}((u))} HP_{\bullet}(A) \leq \dim_{\mathbb{C}}
HH_{\bullet}(A) 
\end{equation}
If $X$ is
a smooth and compact \nc-space, the Hochschild chain complex  of 
$C_{X}$ is the derived tensor product over $A\otimes A^{\op{op}}$ of a
perfect complex with finite dimensional cohomology with itself. In
particular  $HH_{\bullet}(C_{X}) := HH_{\bullet}(A)$ is a finite dimensional 
$\mathbb{C}$-vector space, and so by \eqref{eq:HCinequality} we have
that $HP_{\bullet}(C_{X})$ is finite dimensional over $\mathbb{C}((u))$. 
Thus the $\mathbb{C}[[u]]$-module $HC^{-}_{\bullet}(C_{X})$ is
finitely generated and so corresponds to the formal germ at $u = 0$ of an
algebraic $\mathbb{Z}/2$-graded coherent sheaf on
$\mathbb{A}^{1}_{\mathbb{C}}$. The fiber of this sheaf at $u = 0$  is
$HH_{\bullet}(C_{X})$ and the generic fiber is $HP_{\bullet}(C_{X})$. In
\cite{ks-ncgeometry,kontsevich-lefschetz}
the second author proposed the so called {\em\bfseries degeneration
conjecture} asserting that for a smooth and compact \nc-space $X =
\ncSpec(A)$ we have an equality of dimensions in
\eqref{eq:HCinequality}. In other words the degeneration conjecture
assert that for a smooth and compact \nc-space the
$\mathbb{C}[[u]]$-module $HC^{-}_{\bullet}(C_{X})$ is free of finite
rank and thus corresponds to an algebraic vector bundle on the one
dimensional  formal
disc $\mathbb{D} := \op{Spf}(\mathbb{C}[[u]])$.

\

\begin{rem} \label{rem:degeneration}
There is a lot of evidence supporting the validity of this conjecture.
The work of Weibel \cite{weibel-HC} shows that if $X$ is a usual
quasi-compact and quasi-separated complex scheme the Hochschild and
periodic cyclic homology of $X$ viewed as a \nc-space can be
identified with the algebraic de Rham and Dolbeault cohomology of $X$
respectively. Combined with the degeneration of the Hodge-to-de-Rham
spectral sequence in the smooth proper case this shows that the
degeneration conjecture holds true for usual schemes. Also recently in
a very exciting sequence of papers
\cite{kaledin-cartier,kaledin-degeneration} Kaledin proved the
degeneration conjecture for graded \nc-spaces $X = \ncSpec(A)$ for
which $A$ is concentrated in non-negative degrees. The case of graded
\nc-spaces $X = \ncSpec(A)$ for which $A$ is concentrated in
non-positive degrees was also settled by Shklyarov
\cite{shklyarov}. The general graded case and the
$\mathbb{Z}/2$-graded case are still wide open.
\end{rem}

\

\bigskip

\subsubsection {\bfseries The meromorphic connection in the $u$-direction} \ 
\label{sssec:mero.u} The next observation is that the
$\mathbb{C}\{u\}[u^{-1}]$-module $HP_{\bullet}(C_{X})$ comes equipped
with a natural meromorphic connection. Indeed, recall that by the work
of Getzler \cite{getzler-gm} there is a version of the Gauss-Manin
connection which exists on the periodic cyclic homology of any flat
family of differential graded algebras (see also
\cite{tsygan-gm,kaledin-cyclic}). An analogous statement holds in the
$\mathbb{Z}/2$-graded case as explained e.g. in
\cite[Section~11.5]{ks-ncgeometry}. The Gauss-Manin connection for any
family of dg algebras $\mathcal{A}_{x}$ over the formal disc
$\op{Spf}\mathbb{C}[[x]]$ with a formal parameter $x$, is an operator
\[ 
\nabla^{\op{GM}}_{u\frac{\partial}{\partial x}}:
 H^\bullet( C^{\op{red}}(\mathcal{A}_{x},\mathcal{A}_{x})[[u,x]],
\partial_{\mathcal{A}_x}+u\cdot B_{\mathcal{A}_x})\to H^\bullet(
C^{\op{red}}(\mathcal{A}_{x},\mathcal{A}_{x})[[u,x]],
\partial_{\mathcal{A}_x}+u\cdot B_{\mathcal{A}_x})
\] 
satisfying the
Leibniz rule with respect to the multiplications by $u$ and $x$
(compare this with the {\bfseries (Griffiths transversality axiom)} in
Definition~\ref{defi:variations} from Section~\ref{sssec:variations}).

Suppose now $A$
is a differential $\mathbb{Z}/2$-graded algebra with product $m_{A}$,
differential $d_{A}$, and a strict unit $1_{A}$. Then we can form a
flat family $\mathcal{A} \to \mathbb{A}^{1} - \{ 0 \}$ of differential
$\mathbb{Z}/2$-graded algebras parameterized by the punctured affine
line $\mathbb{A}^{1} - \{ 0 \}$.  The fiber $\mathcal{A}_{t}$ of
$\mathcal{A}$ over a point $t \in \mathbb{A}^{1}_{\mathbb{C}} -
\{ 0 \}$ is the d($\mathbb{Z}/2$)g algebra for which the underlying
$\mathbb{Z}/2$-graded vector space is $A$ and
\[
\begin{split}
m_{\mathcal{A}_{t}} & := t\cdot m_{A}, \\
d_{\mathcal{A}_{t}} & := t\cdot d_{A}, \\
1_{\mathcal{A}_{t}} & := t^{-1}\cdot 1_{A}.
\end{split}
\]
Looking at the scaling properties of $\partial$ and $B$ we see
that the identity morphism on the level of cochains induces  a natural
isomorphism  
\begin{equation} \label{eq:lambda-1}
\xymatrix@1{
H^{\bullet}\left(C^{\op{red}}_{\bullet}(\mathcal{A}_{t},
\mathcal{A}_{t})[[u]],
\partial_{\mathcal{A}_{t}}+u\cdot
B_{\mathcal{A}_{t}}\right) \ar[r]^-{\cong} &
H^{\bullet}(C^{\op{red}}_{\bullet}(A,A)[[u]],\partial + ut^{-2}\cdot B).
}
\end{equation}
This isomorphism does not come from a quasi-isomorphism of complexes,
 as the identity map is not a morphims of complexes: the differentials
 do not coincide but differ by the factor $t$. 
If $A$ is smooth and compact, then the negative cyclic homology of the
family of algebras $\mathcal{A}_{t}$ gives rise to an algebraic
vector bundle $\widetilde{HC}^{-}$ on the product $(\mathbb{A}^{1} -
\{ 0 \})\times \mathbb{D}$. Here $\mathbb{D} :=
\op{Spf}\mathbb{C}[[u]]$ denotes the one dimensional formal disc. We
will write $(t,u)$ for the coordinates on
$(\mathbb{A}^{1} - \{ 0 \})\times \mathbb{D}$. We will be interested
 in fact only in the formal neigborhood of point $t=1$ where we can
 choose as a local coordinate $x:=\op{log}(t)$. The 
Getzler-Gauss-Manin connection then can be viewed as a relative
holomorphic connection $\nabla^{\op{GM}}$ on $\widetilde{HC}^{-}$
which differentiates only along $\mathbb{A}^{1}_{\mathbb{C}} - \{ 0
\}$. On the other hand the formal completion of the group
$\mathbb{C}^{\times}$ at $1$ acts on $(\mathbb{A}^{1}_{\mathbb{C}} -
\{ 0 \})\times \mathbb{D}$ by $(t,u) \mapsto
(\mu t,\mu^{2}u)$ for $\mu \in \mathbb{C}^{\times}$. The
 isomorphism \eqref{eq:lambda-1} gives rise to a
$\mathbb{C}^{\times}$-equivariant structure on the vector bundle
$\widetilde{HC}^{-}$ and the infinitesimal action of $d/d\mu$
associated with this equivariant structure gives  a holomorphic
differential operator $\boldsymbol{\Lambda} \in \op{Diff}^{\leq
  1}(\widetilde{HC}^{-},\widetilde{HC}^{-})$ with symbol equal to
\[
\left(t\frac{\partial}{\partial{t}} +
2u\frac{\partial}{\partial u}\right)\cdot
\op{id}_{\widetilde{HC}^{-}}.
\]
Hence 
\[
\nabla_{u^2\frac{\partial}{\partial u}} :=
\frac{u}{2}\cdot \boldsymbol{\Lambda} -
\nabla^{\op{GM}}_{\frac{ut}{2}\frac{\partial}{\partial t}} 
\]
is a first order differential operator on $\widetilde{HC}^{-}$ with symbol 
\[ u^2
\frac{\partial}{\partial u}\cdot
\op{id}_{\widetilde{HC}^{-}}
\]
and so after restricting $\widetilde{HC}^{-}$ to $\{ 1 \}\times
\mathbb{D}$ this operator gives a meromorphic connection $\nabla$ on
the $\mathbb{C}[[u]]$-module $HC_{\bullet}^{-}(C_{X})$ with at most a
second order pole at $u=0$. Note also that if the algebra $A$ is
$\mathbb{Z}$-graded, then the family
$\mathcal{A}_{t}$  t is easily seen to be trivial and the connection
$\nabla$ has the first 
order pole at $u=0$ with monodromy equal to $(-1)^{\op{parity}}$.

\

\

\bigskip

\subsubsection {\bfseries The $\mathbb{Q}$-structure} \ The categorical
origin of the rational (or integral) structure of the conjectural
\nc-Hodge structure is more mysterious. Conceptually the correct
rational structure should come from the Betti cohomology or, say, the
topological K-theory of the \nc-space.  There are two natural approaches
to constructing the rational structure $\mycal{E}_{B} \subset
HP_{\bullet}(C_{X})$:

\

\noindent
{\bfseries (a)} \ {\bfseries The soft algebra approach
(\cite{kontsevich-lefschetz}).} Let again $X = \ncSpec(A)$ be an
affine \nc-space, and assume $X$ is compact. By analogy with the
classical geometric case one expects that there should exist a nuclear
Frech\'{e}t d$(\mathbb{Z}/2)$g algebra $A_{C^{\infty}}$ so that
\begin{itemize}
\item The K-theory of $A_{C^{\infty}}$ satisfies Bott periodicity,
  i.e. $K_{i}(A_{C^{\infty}}) = K_{i+2}(A_{C^{\infty}})$ for all $i
  \geq 0$. 
\item There is a homomorphism $\varphi : A \to A_{C^{\infty}}$
of  d$(\mathbb{Z}/2)$g algebras  for which 
$\varphi_{*} : HP_{\bullet}(A) \to HP_{\bullet}(A_{C^{\infty}})$
is an isomorphism, and the image of the Chern character map
\[
ch : K_{\bullet}(A_{C_{\infty}}) \to HP_{\bullet}(A_{C^{\infty}})
\]
is an integral lattice, and hence gives a rational structure 
$\mycal{E}_{B} \subset HP_{\bullet}(A)$.
\end{itemize}

\

\noindent
{\bfseries Note:} If $X$ is a smooth and compact complex variety and
if $\mathcal{E} \in \Perf(X)$ is a vector bundle generating $C_{X}$, then
one may take 
\[
\begin{split}
A & := A^{0,\bullet}(X,\mathcal{E}^\vee
\otimes\mathcal{E}),\bar{\partial}) \\
A_{C^{\infty}} & :=
A^{0,0}(X,\mathcal{E}^\vee \otimes\mathcal{E}).
\end{split}
\] 
Note that the algebra $A_{C^{\infty}}$ is Morita equivalent to $C^{\infty}(X)$.

\

\medskip

\noindent
{\bfseries (b)} {\bfseries The semi-topological K-theory approach
  (Bondal, To\"{e}n, \cite{toen-st}).} \ Assume again that $X =
  \ncSpec(A)$ is 
a smooth and compact graded \nc-affine \nc-space. Consider the moduli stack
$\mycal{M}_{X}$ of all objects in $\Perf_{X}$.  This is an
$\infty$-stack which by a theorem of To\"{e}n and Vaquie
\cite{toen-vaquie} is locally geometric and locally of finite
presentation. Moreover for any $a, b \in \mathbb{N}$ the substack
$\mycal{M}_{X}^{[a,b]} \subset \mycal{M}_{X}$ consisting of objects of
amplitude in the interval $[a,b]$ is a geometric $b-a+1$-stack. The
functor sending a complex scheme to the underlying topological space in
the analytic topology gives rise by a left Kan extension to a
topological realization functor
\[
|\bullet | : \op{Ho}\left(\text{\sf{Stacks}}/\mathbb{C}\right) \to
\op{Ho}(\text{\sf Top})
\]
from the homotopy category of stacks to the homotopy category of
complex spaces. Following Friedlander-Walker
\cite{friedlander.walker-st} we define the semi-topological K-group
of the \nc-space $X$ to be 
\[
K^{st}_{0}(X) := \pi_{0}(|\mycal{M}_{X}|).
\]
The group structure here is induced by the direct sum $\oplus$ of
$A$-modules: the monoid \linebreak
$\left(\pi_{0}(|\mycal{M}_{X}|),\oplus\right)$ is actually a group. To
see this note that for any $A$-module $E$ we have that $\left[E\oplus
  E[1]\right] = 0$ in $\pi_{0}(|\mycal{M}_{X}|)$. Indeed we have
distinguished triangles 
\[
\xymatrix@R-1pc{
E \ar[r] & 0 \ar[r] & E[1] \ar[r] & E[1] \\
E \ar[r] & E\oplus E[1]  \ar[r] & E[1] \ar[r] & E[1] 
}
\]
the first of which corresponds to $\op{id} \in \op{Ext}^{1}(E[1],E) =
\op{Hom}(E,E)$, and the second corresponds to $0 \in \op{Ext}^{1}(E[1],E) =
\op{Hom}(E,E)$. Since $\op{Ext}^{1}(E[1],E) =
\op{Hom}(E,E)$ is a vector space, it follows that $\op{id}$ deforms to
$0$ continuously and so the second terms in the above triangles
represent the same point in $\pi_{0}(|\mycal{M}_{X}|)$. 

More generally $\oplus$ makes $|\mycal{M}_{X}|$ into an $H$-space
$\mathbb{K}^{st}(X)$ which is the degree zero part of a natural
spectrum. Using this one can define $K^{st}_{i}(X)$ for all $i \geq
0$.

Next note that since $C_{X}$ is triangulated it is a
module over the category $\Perf_{\op{pt}}$ of complexes of
$\mathbb{C}$-vector spaces with finite dimensional total
cohomology. In particular $K^{st}_{\bullet}(X)$ is a graded module
over $K^{st}_{\bullet}(\op{pt})$. It can be checked that 
\[
\mathbb{K}^{st}(\op{pt}) = BU = \mathbb{K}^{\op{top}}(\op{pt}),
\]
and so $K^{st}_{\bullet}(X)$ is a graded $\mathbb{Z}[u]$-module
($\deg u = 2$). 

Now we can define 
\[
K^{\op{top}}_{\bullet}(X) := K^{st}_{\bullet}(X)[u^{-1}] =
K^{st}_{\bullet}(X) \otimes_{\mathbb{Z}[u]} \mathbb{Z}[u,u^{-1}].
\]
Again one expects that there is a Chern character map
\[
ch : K^{\op{top}}_{\bullet}(X)  \to HP_{\bullet}(C_{X})
\]
whose image gives a rational structure $\mycal{E}_{B}$ on $HP_{\bullet}(C_{X})$.

\

\noindent
{\bfseries Note:} If $X$ is a smooth and compact complex variety, then
the Friedlander-Walker comparison theorem
\cite{friedlander.walker-comparison} implies that 
$K^{\op{top}}(D(\text{\sf QCoh}(X))) \cong
K^{\op{top}}(X^{\op{top}})$, where $X^{\op{top}}$ is the topological
space underlying $X$.

\

\bigskip

\subsubsection {\bfseries Questions} \ Even though we have some good
candidates for the ingredients  $H$, $\nabla$,
$\mycal{E}_{B}$ of the conjectural \nc-Hodge structure associated with
a \nc-space, there are several important problems that need to be
addressed before one can prove Conjecture~\ref{con:main.conjecture}:
\begin{itemize}
\item show that the connection $\nabla$ has regular singularities
  (this is automatically true  in the $\mathbb{Z}$-graded case);
\item show that $\nabla$ preserves the rational structure;
\item show that the opposedness axiom hold. 
\end{itemize}

\

\noindent
One can show that for a smooth compact nc-space the coefficient
 $A_{-2}$ in the $u$-connection is a nilpotent operator, which gives
 an evidence in favor of the regular singularity question.

\

\noindent
In fact Conjecture~\ref{con:main.conjecture} and the above questions
are special cases of a general conjecture which predicts the existence
of a general \nc-Hodge structure on the periodic cyclic homology of a
curved d$(\mathbb{Z}/2)$g category which is formally smooth and
compact. We will not discuss the general conjecture or the relevant
constructions here but we will revisit these questions in some
interesting geometric examples in Section~\ref{sec:examples}.

\subsection{\bfseries Gluing data} \label{sec-gluing}

In this section we discuss how general \nc-Hodge structures of
exponential type can be glued together out of \nc-Hodge structures
with regular singularities and additional gluing data.

\

\bigskip

\subsubsection {\bf \nc-De Rham data}  \ \label{sssec-deRham-gluing} 
The de
Rham part of a \nc-Hodge structure of exponential type can be
prescribed in three equivalent ways:

\

\smallskip

\noindent
{\bf \nc dR(i)} A pair $(\mycal{M},\nabla)$, where $\mycal{M}$ is a
finite dimensional vector space over $\mathbb{C}\{u\}[u^{-1}]$ and
$\nabla$ is a meromorphic connection. These data should satisfy the
following

\begin{description}
\item[{\bf Property \nc dR(i):}] \ There exist:

\

\noindent
$\bullet$ 
 a frame $\underline{e} =
  (e_{1}, \ldots, e_{r})$ of $\mycal{M}$ over
  $\mathbb{C}\{u\}[u^{-1}]$ in which
\[
\nabla = d + \left(\sum_{k \geq -2} A_{k} u^{k}\right)du
\]
with $A_{k} \in \op{Mat}_{r\times r}(\mathbb{C})$, $r =
\op{rank}_{\mathbb{C}\{u\}[u^{-1}]} \mycal{M}$. In other
words, there is a holomorphic extension $\mycal{H} =
\mathbb{C}\{u\}e_{1}\oplus \ldots \oplus \mathbb{C}\{u\}e_{r}$ in
which $\nabla$ has at most a second order pole.

\

\noindent
$\bullet$  a formal isomorphism over $\mathbb{C}((u))$:
\[
(\mycal{M},\nabla)\otimes_{\mathbb{C}\{u\}[u^{-1}]}\mathbb{C}((u))
\stackrel{\cong}{\to} 
\bigoplus_{i=1}^{m} \be^{\bc_{i}/u}\otimes
(\mycal{R}_{i},\nabla_{i})
\]
where $(\mycal{R}_{i},\nabla_{i})$ are meromorphic bundles with
connections with regular singularities and $\bc_{1}, \ldots,
\bc_{m} \in \mathbb{C}$ denote the  distinct 
eigenvalues of $A_{-2}$.
\end{description}

\

\smallskip

\noindent
{\bf \nc dR(ii)} A pair $(M,\nabla)$, where $M$ is an algebraic vector
bundle on $\mathbb{A}^{1}-\{0\}$  and $\nabla$ is a  connection on
$M$. These data should satisfy the following

\begin{description}
\item[{\bf Property \nc dR(ii):}] $M$ can be extended to an algebraic
  vector bundle $\widetilde{M}$ on $\mathbb{P}^{1}$, and 

\

\noindent
$\bullet$ \ with respect
  to this extension and appropriate local trivializations at zero and
  infinity we must have
\[
\begin{split}
\nabla & =  d + \left(\sum_{k \geq -2} A_{k} u^{k}\right)du  \qquad
\text{near} \ 0, \\
\nabla & =  d + \left(\sum_{k \geq -1} B_{k} u^{-k}\right)d(u^{-1})  \qquad
\text{near} \ \infty.
\end{split}
\]
In other words $\nabla : \widetilde{M} \to
\widetilde{M}\otimes_{\mathcal{O}_{\mathbb{P}^{1}}}
  \Omega^{1}_{\mathbb{P}^{1}}(2\cdot \{0\} + \{\infty\})$. 

\

\noindent
$\bullet$ There is a formal isomorphism over $\mathbb{C}((u))$:
\[
(M,\nabla)\otimes_{\mathbb{C}[u,u^{-1}]}\mathbb{C}((u)) 
 \stackrel{\cong}{\to}
\bigoplus_{i=1}^{m} \be^{\bc_{i}/u}\otimes
(\mycal{R}_{i},\nabla_{i})
\]
where $(\mycal{R}_{i},\nabla_{i})$ are meromorphic bundles with
connections with regular singularities and $\bc_{1}, \ldots,
\bc_{m} \in \mathbb{C}$ denote the  distinct 
eigenvalues of $A_{-2}$.
\end{description}

\

\smallskip

\noindent
{\bf \nc dR(iii)} An algebraic holonomic  $\mathcal{D}$-module $\bM$ on
$\mathbb{A}^{1}$. The $\mathcal{D}$-module $\bM$ should also satisfy the
following
\begin{description}
\item[{\bf Property \nc dR(iii):}] $\bM$ has regular singularities and 
$H_{DR}^{\bullet}(\mathbb{A}^{1},\bM) = 0$.
\end{description}

\

\smallskip

\noindent
The \nc-de Rham data of types {\bf \nc dR(i)}, {\bf \nc dR(ii)}, and
  {\bf \nc dR(iii)} 
  form obvious full subcategories in the categories of meromorphic
  connections over $\mathbb{C}\{u\}[u^{-1}]$, algebraic vector bundles
  with connections on $\mathbb{A}^{1}-\{0\}$, and coherent algebraic
  $\mathcal{D}$-modules on $\mathbb{A}^{1}$ respectively. We have the
  following

\begin{lemma} \label{lem-dr-equivalence} The categories of \nc-de
  Rham data of types 
  {\bf \nc dR(i)}, {\bf \nc dR(ii)}, and \linebreak 
{\bf \nc dR(iii)} are all equivalent.
\end{lemma}
{\bf Proof.}  In essence we have already discussed the equivalence
{\bf \nc dR(i)} $\Leftrightarrow$ {\bf \nc dR(ii)} in
Remark~\ref{rem-algebraize-lattice}.  Explicitly   
$(\mycal{M},\nabla) =
\mathfrak{G}_{0}((M,\nabla)) =
(M\otimes_{\mathbb{C}[u,u^{-1}]}\mathbb{C}\{u\}[u^{-1}],
\nabla)$.  

We define a functor $\mathfrak{F} : (\text{data {\bf (iii)}}) \to
(\text{data {\bf (ii)}})$ as follows. Let $\bM$ be a regular holonomic
algebraic $\mathcal{D}$-module on $\mathbb{A}^{1}$  with
trivial de Rham cohomology. Denote the coordinate on $\mathbb{A}^{1}$
by $v$. The vanishing of de Rham cohomology means that  the
action $\frac{d}{dv} : \bM \to \bM$ is an invertible operator.  
Consider the algebraic Fourier transform $\boldsymbol{\Phi}\bM$ which
is the same 
vector space as $\bM$ endowed with action of the Weyl algebra defined by 
\[
\begin{split}
\tilde{v} & :=\frac{d}{dv} \\[0.5pc]
\frac{d}{d\tilde{v}} & :=-v
\end{split}
\] 
where $\tilde{v}$ is the coordinate on
the dual line.  By our assumptions $\boldsymbol{\Phi}{\bM}$ is a holonomic
$\mathcal{D}$-module on which $\tilde{v}$ acts invertibly.  Hence
$\boldsymbol{\Phi}{\bM}$ is the direct image of a holonomic $\mathcal{D}$-module
$\boldsymbol{\Phi}{\bM}'$ on $\mathbb{A}^{1} - \{0\}$ under the embedding
\[
\left(\mathbb{A}^{1} - \{0\}\right) \hookrightarrow
\mathbb{A}^{1}=\op{Spec}(\mathbb{C}[\tilde{v}])
\]
Finally, making the change of coordinates $u=1/\tilde{v}$ we obtain a
$\mathcal{D}$-module $M$ on $\mathbb{A}^{1} - \{0\}$ with
coordinate $u$.
     
We claim that $\mathfrak{F}(\bM) := M$ obtained in this way satisfies
the property {\bf \nc dR(ii)}, and that by this construction one
obtains all such modules.  It follows from the well-known properties
of the Fourier transform that $\boldsymbol{\Phi}{\bM}$ has no
singularities in $\mathbb{A}^{1} - \{0\}$ and the its singularity at
$\tilde{v}=0$ is regular. Hence $M$ is a vector bundle on
$\mathbb{A}^{1} - \{0\}$ endowed with connection with regular
singularity at $\infty$. It only remains to to use the well-known fact
(see e.g. \cite[Chapters~IX-XI]{malgrange-book} or
\cite[Theorem~2.10.16]{katz-exponential}) that the exponential type
property for $M$ is equivalent to the property of $\bM$ to have only
regular singularities. \ \hfill $\Box$

\

\medskip

\begin{rem} \label{rem:regular.sing} The characterization of the
  exponential type property in terms of the Fourier transform can be
  stated more precisely (see \cite[Chapters~IX-XI]{malgrange-book} or
\cite[Theorem~2.10.16]{katz-exponential}): For an algebraic
  holonomic $\mathcal{D}$-module  $M$ on the complex affine line,  the
  following two conditions are equivalent: 
\begin{itemize}
\item[1)] $M$ has regular singularities;
\item[2)] the Fourier transform $\boldsymbol{\Phi} M$ of $M$ has no
  singularities 
  outside $0$, its singularity  at $0$ is regular, and its
  singularity at infinity is of exponential type.
\end{itemize}

Explicitly $\boldsymbol{\Phi} M$ being of exponential type at infinity
 means that 
 if $x$ is a coordinate on $\mathbb{A}^{1}$ centered at
 $0$, then  after 
 passing to the formal completion 
$\left(\boldsymbol{\Phi} M\right) \otimes_{\mathbb{C}[x]}
  \mathbb{C}((x^{-1}))$  
the resulting module will be isomorphic to a finite sum 
\[
\bigoplus_{i=1}^{m} \be^{\bc_{i}x}\otimes
(\mycal{R}_{i},\nabla_{i})
\]
where $(\mycal{R}_{i},\nabla_{i})$ are  $\mathcal{D}$-modules with
a regular singularity at infinity. 
\end{rem}

\

\begin{rem} \label{rem-dR-other-fields} Note that the de Rham data {\bf
    \nc dR(i)} is analytic in nature, whereas  {\bf
    \nc dR(ii)} and {\bf 
    \nc dR(iii)} are algebraic. In fact from the proof it is clear
    that {\bf \nc dR(ii)} and {\bf 
    \nc dR(iii)} and their equivalence still make sense if we replace
  $\mathbb{C}$ with any field of characteristic zero.
\end{rem}

\

 \bigskip

\subsubsection {\bfseries \nc-Betti data} \ The (rational) Betti part of a
\nc-Hodge structure of exponential type can be prescribed in four ways:

\

\medskip

\noindent
{\bf \nc B(i)} A (middle perversity) perverse sheaf $\mycal{G}^{\bullet}$ of
$\mathbb{Q}$-vector spaces on the Riemann surface $\mathbb{C}$ (taken
with the analytic topology)  satisfying the
following 

\begin{description}
\item[{\bf Property \nc B(i):}] $R\Gamma(\mathbb{C},\mycal{G}^{\bullet}) =
  0$. 
\end{description}

\

\medskip

\noindent
{\bf \nc B(ii)} A constructible sheaf $\mycal{F}$ of
$\mathbb{Q}$-vector spaces on the Riemann surface $\mathbb{C}$ (taken
with the analytic topology)  satisfying the
following 

\begin{description}
\item[{\bf Property \nc B(ii):}] $R\Gamma(\mathbb{C},\mycal{F}) =
  0$. 
\end{description}

\

\medskip

\noindent
{\bf \nc B(iii)} A finite collection of distinct points $S =
\{\bc_{1},\ldots, \bc_{n} \} \subset \mathbb{C}$, and 
\begin{itemize}
\item a collection $U_{1}, U_{2}, \ldots, U_{n}$ of finite dimensional
non-zero $\mathbb{Q}$-vector spaces,
\item a collection of linear maps $T_{ij} : U_{j} \to U_{i}$, for all
  $i, j = 1, \ldots, n$,
\end{itemize}
satisfying the following 
\begin{description}
\item[{\bf Property \nc B(iii):}] $T_{ii} \in GL(U_{i})$. 
\end{description}

\

\medskip

\noindent
{\bf \nc B(iv)} A local system $\stS$ of $\mathbb{Q}$-vector spaces on
$\bS^{1}$ equipped with a filtration $\{\stS_{\leq \lambda}\}_{\lambda
  \in \mathbb{R}}$ by subsheaves of $\mathbb{Q}$-vector spaces,
satisfying the following  

\begin{description}
\item[{\bfseries Property \nc B(iv):}] The filtration $\{\stS_{\leq
  \lambda}\otimes \mathbb{C}\}_{\lambda \in \mathbb{R}}$ of
  $\stS\otimes \mathbb{C}$ is a
  Deligne-Malgrange-Stokes filtration of exponential type. In other
  words, there exist complex numbers $\bc_{1}, \ldots, \bc_{n} \in
  \mathbb{C}$ so that:
\begin{itemize}
\item For every $\varphi \in \bS^{1}$, the filtration $\{(\stS_{\leq
  \lambda}\otimes \mathbb{C})_{\varphi} \}_{\lambda \in \mathbb{R}}$
  of the stalk 
  $(\stS\otimes \mathbb{C})_{\varphi}$ jumps exactly at the real
  numbers $\{ \op{Re}\left(\bc_{k}e^{-i\varphi}\right) \}_{k = 1}^{n}$.
\item The associated graded sheaves of  $\stS\otimes \mathbb{C}$ with
  respect to $\{\stS_{\leq
  \lambda}\otimes \mathbb{C}\}_{\lambda \in \mathbb{R}}$ are local
  systems on $\bS^{1}$. 
\end{itemize}
\end{description}

\

\smallskip

\noindent
Again there are natural equivalences of the different types of Betti
data (for {\bf \nc B(iii)} the equivalence depends on
certain choices of paths as one can see from the proof of
Theorem~\ref{prop:B(i)=B(ii)} and the statement of
Lemma~\ref{lemma:B(ii)=B(iii)}.).   Consider the full subcategories ({\bf
  \nc B(i)}) and ({\bf \nc 
  B(ii)}) of  \nc-Betti data of types {\bf
  \nc B(i)} and {\bf \nc B(ii)} in the category of perverse sheaves of
$\mathbb{Q}$-vector spaces on $\mathbb{C}$ and in the category of
constructible sheaves of $\mathbb{Q}$-vector spaces on $\mathbb{C}$
respectively. We have the following

\begin{theo} \label{prop:B(i)=B(ii)}
The categories of \nc-Betti data of types {\em {\bf
  \nc B(i)}} and {\em {\bf \nc B(ii)}} are naturally equivalent.
More precisely, the natural functors
\[
\mathcal{H}^{-1} : D^{b}_{\op{constr}}(\mathbb{C},\mathbb{Q}) \to
\op{\sf{Constr}}(\mathbb{C},\mathbb{Q})  \qquad \text{and} \qquad  
[1] : \op{\sf{Constr}}(\mathbb{C},\mathbb{Q}) \to
D^{b}_{\op{constr}}(\mathbb{C},\mathbb{Q}) 
\]
induce mutually inverse equivalences of the full subcategories
$\text{{\em ({\bf \nc B(i)})}} \subset
D^{b}_{\op{constr}}(\mathbb{C},\mathbb{Q})$ and $\text{{\em ({\bf \nc
    B(ii)})}} \subset \op{\sf{Constr}}(\mathbb{C},\mathbb{Q})$.
\end{theo}  
{\bf Proof.} First we look at the data {\bf \nc B(i)} more closely.
Suppose $X$ is a complex analytic space underlying a complex
quasi-projective variety.  Recall (see e.g.  \cite{bbd,ks,dimca}) that
a bounded complex $\mycal{G}^{\bullet}$ of sheaves of
$\mathbb{C}$-vector spaces on $X$ is called a (middle perversity)
perverse sheaf if it has constructible cohomology sheaves
$\mathcal{H}^{k}(\mycal{G}^{\bullet})$ and if
\begin{itemize}
\item for all $k \in \mathbb{Z}$, we have $\dim_{\mathbb{R}} \{ x \in X | \;
  \mathcal{H}^{-k}(\bi_{x}^{*}\mycal{G}^{\bullet}) \neq 0\} \leq 2k$,
\item for all $k \in \mathbb{Z}$, we have $\dim_{\mathbb{R}} \{ x \in X | \;
  \mathcal{H}^{k}(\bi_{x}^{!}\mycal{G}^{\bullet}) \neq 0\} \leq 2k$.
\end{itemize}
Here $\bi_{x} : x \hookrightarrow X$ denotes the inclusion of the
point $x$ in $X$. 

For future reference we will write $D^{b}_{\op{constr}}(X,\mathbb{Q})$
for the derived category of complexes of $\mathbb{Q}$-vector spaces on
$X$ with constructible cohomology, $\perv(X,\mathbb{Q}) \subset
D^{b}_{\op{constr}}(X,\mathbb{Q})$ for the full subcategory of middle
perversity perverse sheaves, and $\op{\sf{Constr}}(X,\mathbb{Q}) \subset
D^{b}_{\op{constr}}(X,\mathbb{Q})$ for the full subcategory of
constructible sheaves.

From the definition it is clear that if $\mycal{G}^{\bullet}$ is a
perverse sheaf on $\mathbb{C}$, then $\mycal{G}^{\bullet}$ has at most
two non-trivial cohomology sheaves
$\mathcal{H}^{-1}(\mycal{G}^{\bullet})$ and
$\mathcal{H}^{0}(\mycal{G}^{\bullet})$. Moreover the support of
$\mathcal{H}^{0}(\mycal{G}^{\bullet})$ has dimension $\leq 0$. Now the
cohomology $R\Gamma^{\bullet}(\mycal{G}^{\bullet}) =
\mathbb{H}^{\bullet}(\mathbb{C},\mycal{G}^{\bullet})$ can be computed
via the hypercohomology spectral sequence
\[
E_{2}^{pq} = H^{p}(\mathbb{C},\mathcal{H}^{q}(\mycal{G}^{\bullet}))
\Rightarrow \mathbb{H}^{p+q}(\mathbb{C},\mycal{G}^{\bullet}).
\]
Since $\mycal{G}^{\bullet}$ has only two cohomology sheaves, the
$E_{2}$ sheet of this spectral sequence is

\

\def\csg#1{\save[].[drrr]!C*+<2pc>[F-,]\frm{}\restore}
\[
\xymatrix@=1.5pc{ 
0 & \csg1  
H^{0}({\mathbb C},\mathcal{H}^{0}(\mycal{G}^{\bullet})) \ar[drr]  &
H^{1}({\mathbb C},\mathcal{H}^{0}(\mycal{G}^{\bullet})) \ar[drr] &
H^{2}({\mathbb C},\mathcal{H}^{0}(\mycal{G}^{\bullet})) & \cdots \\  
-1 & H^{0}({\mathbb C},\mathcal{H}^{-1}(\mycal{G}^{\bullet})) 
& H^{1}({\mathbb C},\mathcal{H}^{-1}(\mycal{G}^{\bullet}))  &
H^{2}({\mathbb C},\mathcal{H}^{-1}(\mycal{G}^{\bullet})) & \cdots \\
  & 0 & 1 & 2 &  
}
\]
By Artin's vanishing theorem for constructible sheaves
\cite[Corollary~3.2]{artin-vanishing} he have $H^{p}({\mathbb
  C},\mathcal{H}^{q}(\mycal{G}^{\bullet})) = 0$ for all $q$ and all $p
> 1$. Furthermore since $\mathcal{H}^{0}(\mycal{G}^{\bullet})$ has at
most zero dimensional support we have  $H^{1}({\mathbb 
  C},\mathcal{H}^{0}(\mycal{G}^{\bullet})) = 0$. In particular the
spectral sequence degenerates at $E_{2}$ and the only potentially
non-trivial cohomology groups of $\mycal{G}^{\bullet}$ are
\[
\begin{split}
\mathbb{H}^{-1}(\mathbb{C},\mycal{G}^{\bullet}) & =
H^{0}(\mathbb{C},\mathcal{H}^{-1}(\mycal{G}^{\bullet})), \text{ and }
\\
\mathbb{H}^{0}(\mathbb{C},\mycal{G}^{\bullet}) & =
H^{1}(\mathbb{C},\mathcal{H}^{-1}(\mycal{G}^{\bullet})) \oplus
H^{0}(\mathbb{C},\mathcal{H}^{0}(\mycal{G}^{\bullet})). 
\end{split}
\]
Thus under the assumption that $R\Gamma(\mycal{G}^{\bullet}) = 0$ we
get that $H^{0}(\mathbb{C},\mathcal{H}^{0}(\mycal{G}^{\bullet})) = 0$,
i.e. that $\mathcal{H}^{0}(\mycal{G}^{\bullet}) = 0$. In other words 
$\mycal{G}^{\bullet} = \mycal{F}[1]$ for some constructible sheaf
$\mycal{F}$ with $R\Gamma(\mycal{F}) = 0$.

To finish the proof of the theorem we need to show that for every
constructible sheaf $\mycal{F}$ with $R\Gamma(\mycal{F}) = 0$, the
object $\mycal{F}[1]$ will be perverse (for the middle
perversity). For this we will have to look more closely at
constructible sheaves on the complex line.

Suppose $\mycal{F}$ is a constructible sheaf of $\mathbb{Q}$ vector
spaces on $\mathbb{C}$. Then there is a finite set $S =
\{\bc_{1},\ldots, \bc_{n}\}$ of points in $\mathbb{C}$ so that
$\mathbb{C}-S$ is the maximal open set on which $\mycal{F}$ restricts
to a local system. Let $\mathbb{F} := \mycal{F}_{|\mathbb{C} - S}$
denote this local system. Let $\mathbb{C}-S
\stackrel{\bj}{\hookrightarrow} \mathbb{C}
\stackrel{\bi}{\hookleftarrow} S$ be the natural inclusions and let 
$\varphi : \mycal{F} \to  j_{*}j^{*}\mycal{F} = j_{*}\mathbb{F}$ be
the adjunction homomorphism.

\begin{figure}[!ht]
\begin{center}
\psfrag{s}[c][c][1][0]{{$\bc_{i}$}}
\psfrag{o}[c][c][1][0]{{$\bc_{0}$}}
\psfrag{os}[c][c][1][0]{{$\bo_{i}$}}
\psfrag{Ds}[c][c][1][0]{{$\bD_{i}$}}
\psfrag{as}[c][c][1][0]{{$a_{i}$}}
\epsfig{file=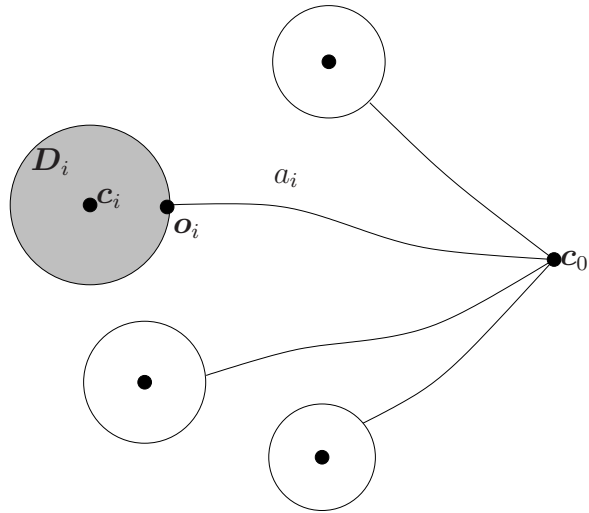,width=3in} 
\end{center}
\caption{A system of paths for $S\subset \mathbb{C}$.}\label{fig:pathstodiscs} 
\end{figure}

Before we can describe $\mathbb{F}$ and $\mycal{F}$ via the
quiver-like data of type {\bf \nc B(iii)} we will need to make some
rigidifying choices. First we fix a base point $\bc_{0} \in \mathbb{C} -
S$.  For $i = 1, \ldots, n$ we choose a collection of a small disjoint
discs $\bD_{i} \subset \mathbb{C}$, each $\bD_{i}$ centered at
$\bc_{i}$. For each disc we fix a point $\bo_{i} \in \partial \bD_{i}$
and denote by $\mathfrak{l}_{i}$  the loop starting and ending at $\bo_{i}$ and
tracing $\partial \bD_{i}$ once in the counterclockwise direction. We fix
an ordered system of non-intersecting paths $\{a_{i}\}_{i=1}^{n} \subset
\mathbb{C}-(\cup_{i=1}^{n} \bD_{i})$ which connect the base point
$\bc_{0}$ with the each of the $\bo_{i}$ as in
Figure~\ref{fig:pathstodiscs}.

Let $\op{mon}_{\mathfrak{l}_{i}} : \mathbb{F}_{\bo_{i}} \to \mathbb{F}_{\bo_{i}}$
be the monodromy operator associated with the local system
$\mathbb{F}$ and the loop $\mathfrak{l}_{i}$. The stalk
$(\bj_{*}\mathbb{F})_{\bc_{i}}$ of the constructible sheaf
$\bj_{*}\mathbb{F}$ at $\bc_{i}$ can be identified naturally with the
subspace $\mathbb{F}_{\bo_{i}}^{\op{mon}_{\mathfrak{l}_{i}}}$ of invariants for
the local monodromy.  Taking stalks at each $\bc_{i} \in S$ we get
$\mathbb{Q}$-vector spaces $\mycal{F}_{\bc_{i}}$ and the adjunction map
$\varphi : \mycal{F} \to \bj_{*}\mathbb{F}$ induces linear maps
\[
\varphi_{\bc_{i}} : \mycal{F}_{\bc_{i}} \to
\mathbb{F}_{\bo_{i}}^{\op{mon}_{\bc_{i}}} \subset \mathbb{F}_{\bo_{i}}.
\]
Note that, by descent, specifying the constructible sheaf $\mycal{F}$ is
equivalent to specifying the collection of points $S \subset
\mathbb{C}$, the local system $\mathbb{F}$ on $\mathbb{C} - S$, the
collection of vector spaces $\{ \mycal{F}_{\bc_{i}} \}_{i = 1}^{n}$ and
the collection of linear maps $\{ \varphi_{\bc_{i}} \}_{i = 1}^{n}$.  In
particular, the compactly supported pullback of $\mycal{F}[1]$ via the
inclusion $\bi_{\bc_{i}} : \{\bc_{i}\} \hookrightarrow \mathbb{C}$ can be
computed in terms of these linear algebraic data and is given
explicitly by the complex
\[
\bi_{\bc_{i}}^{!}(\mycal{F}[1]) = [\xymatrix@R-2pc@C+1pc{
    \mycal{F}_{\bc_{i}} 
  \ar[r]^-{\varphi_{\bc_{i}}} & \mathbb{F}_{\bo_{i}} \ar[r]^-{1 -
    \op{mon}_{\mathfrak{l}_{i}}} &  \mathbb{F}_{\bo_{i}}  \\ -1 & 0 & 1}].
\]
By definition $\mycal{F}[1]$ is a perverse sheaf iff for all $\bc_{i}
\in S$ the complex of vector spaces
$\bi_{\bc_{i}}^{!}(\mycal{F}[1])$ has no cohomology in strictly negative
degrees, i.e. iff $\varphi_{\bc_{i}}$ is injective for all $i = 1,
\ldots, n$. 

Next we rewrite the condition $R\Gamma(\mathbb{C},\mycal{F}) = 0$ in
terms of the descent data \linebreak $\left(\mathbb{F},\{\mycal{F}_{\bc_{i}}\},
\{\varphi_{\bc_{i}}\}\right)$. To simplify notation let $U :=
\mycal{F}_{\bc_{0}}$, $V_{i} = \mycal{F}_{\bc_{i}}$ for $i = 1, \ldots,
n$. Let $T_{i} : U \to U$ be the monodromy operator for the local
system $\mathbb{F}$ and the $\bc_{0}$-based loop $\gamma_{i}$ obtained
by first tracing the path $a_{i}$ from $\bc_{0}$ to $\bo_{i}$, then
tracing the loop $\mathfrak{l}_{i}$, and then tracing back $a_{i}$ in the
opposite direction.  Similarly we have linear maps $\psi_{i} : V_{i}
\to U^{T_{i}} \subset U$ obtained by conjugating $\varphi_{\bc_{i}} :
V_{i} \to \mathbb{F}_{\bo_{i}}$ with the operator of parallel
transport in $\mathbb{F}$ along the path $a_{i}$. 

The descent data for $\mycal{F}$ with respect to the open cover
$\mathbb{C} = (\mathbb{C}-S) \cup \left( \cup_{i = 1}^{n}
\bD_{i}\right)$ are now completely encoded in the linear algebraic data
$\left(U, \{V_{i}\}_{i=1}^{n},\{T_{i}\}_{i =
  1}^{n},\{\psi_{i}\}_{i=1}^{n}\right)$.  Cover $\mathbb{C}$ by the
two opens $\mathbb{C} - S$ and $\cup_{i = 1}^{n} \bD_{i}$. The
intersection of these two opens is the disjoint union of punctured
discs $\coprod_{i = 1}^{n} (\bD_{i} - \bc_{i})$. The Mayer-Vietoris
sequence for $\mycal{F}$ and this cover identifies
$R\Gamma(\mathbb{C},\mycal{F})$ with the complex:
\def\mvg#1{\save[].[dd]!C="mvg#1"*[F]\frm{}\restore}
\[
\xymatrix@C+4pc@R-2pc{
\mvg1  \bigoplus_{i=1}^{n} V_{i}
\ar[r]^-{\oplus_{i=1}^{n} \psi_{i}} & \mvg2
U^{\oplus n} \ar[r]^-{\oplus_{i=1}^{n}(1-T_{i})} & *+[F]{U^{\oplus n}} \\
\bigoplus & \bigoplus & \\
U \ar[uur]^-{\op{id}_{U}^{\oplus n}}
\ar[r]_-{\oplus_{i=1}^{n}(1-T_{i})} & U^{\oplus
  n}\ar[uur]^-{-\op{id}_{U}^{\oplus n}} & \\
& & \\
0 & 1 & 2
}
\]
In other words we have a quasi-isomorphism of complexes of
$\mathbb{Q}$-vector spaces:
\def\gag#1{\save[].[ddr]!C*+<2.4pc>[F-,]\frm{}\restore}
\[
R\Gamma(\mathbb{C},\mycal{F}) \cong \xymatrix@C+4pc@R-2pc{
\gag1 \mvg1  \bigoplus_{i=1}^{n} V_{i} \ar[dr]^-{\oplus_{i=1}^{n} \psi_{i}}
& \\
\bigoplus & *+[F]{U^{\oplus n}} \\
U \ar[ur]_-{\op{id}_{U}^{\oplus n}} & \\
& & \\
& & \\
0 & 1
}
\]
The acyclicity of this complex is equivalent to the conditions
\begin{itemize}
\item[(a)] the maps $\psi_{i} : V_{i} \to U$ are
injective for all $i = 1, \ldots, n$, and
\item[(b)] the map $U \to \oplus_{i=1}^{n} U/V_{i}$ is an isomorphism.
\end{itemize}

Thus the acyclicity of $R\Gamma(\mathbb{C},\mycal{F})$ implies
the perversity of $\mycal{F}[1]$. The theorem is proven. \newline
\mbox{\quad} \ \hfill
$\Box$

\

\bigskip

\noindent
The conditions (a) and (b) from the proof of
Theorem~\ref{prop:B(i)=B(ii)} suggest a better way of recording
the linear algebraic content of $\mycal{F}$. Namely, if we set $U_{i}
:= U/V_{i}$, then we can use (b) to identify $U$ with
$\oplus_{i=1}^{n} U_{i}$, $V_{i}$ with $\oplus_{j\neq i} U_{j}$ and the
map $\psi_{i} : V_{i} \hookrightarrow U$ with the natural inclusion
$\oplus_{j\neq i} U_{j} \subset \oplus_{i=1}^{n} U_{i}$. The only
thing left is the data of the monodromy operators $T_{i} \in GL(U)$,
$i = 1, \ldots, n$. However for each $i$ we have embedding
\[ 
\xymatrix@1{ 
V_i \hspace{0.5pc} \ar@{^{(}->}[r]^-{\psi_{i}} &  \op{Ker}\left[ U 
    \stackrel{(1-T_i)}{\longrightarrow} U\right]
}
\]
and so under the decomposition $U = \oplus_{i=1}^{n} U_{i}$ the
automorphism $T_{i}$ has a block form 
\[
T_{i} = \begin{pmatrix} 
1 & 0 & \cdots & 0 & T_{1i} & 0 & \cdots & 0  \\
0 & 1 & \cdots & 0 & T_{2i} & 0 & \cdots & 0  \\
\cdots & \cdots & \cdots & \cdots & \cdots  & \cdots & \cdots & \cdots
\\
0 & 0 & \cdots & 0 & T_{ii} & 0 & \cdots & 0  \\
0 & 0 & \cdots & 0 & T_{i+1,i} & 1 & \cdots & 0  \\
\cdots & \cdots & \cdots & \cdots & \cdots  & \cdots & \cdots & \cdots
\\
0 & 0 & \cdots & 0 & T_{ni} & 0 & \cdots & 1  \\
\end{pmatrix} 
\]
where $T_{i|U_{i}} = \sum_{j = 1}^{n} T_{ji}$, and $T_{ji} :U_{i} \to
U_{j}$. The linear maps $T_{ji}$ are unconstrained except for the
obvious condition that for all $i$ the map $T_{i}$ should be
invertible, which is equivalent to $T_{ii} : U_{i} \to U_{i}$ being
invertible for all $i = 1, \ldots, n$. Also since $S$ was chosen to be
such that $\mathbb{C} - S$ is the maximal open on which $\mycal{F}$ is
a local system, it follows that $U_{i} \neq \{ 0 \}$ for all $i = 1,
\ldots, n$.

In other words we have proven the following

\begin{lemma} \label{lemma:B(ii)=B(iii)}
Fix the set of points $S = \{ \bc_{1}, \ldots, \bc_{n}\}$ and  choose the
discs $\{ \bD_{i} \}_{i=n}$ and the system of paths $\{ a_{i}
\}_{i=1}^{n}$. The functor assigning to a constructible sheaf
$\mycal{F}$ with singularities at $S$ the data
$(\{U_{i}\}_{i=1}^{n},\{T_{ij}\})$ establishes an equivalence between
the groupoid of all data of type {\em {\bfseries \nc B(ii)}} with
singularities exactly at $S$ and all data  of type {\em {\bfseries \nc
  B(iii)}} with the given $S$.
\end{lemma}

\

\

\medskip

\noindent
The bridge between the \nc- de Rham and Betti data is provided as
usual by the Riemann-Hilbert correspondence. This is tautological but
we record it for future reference:

\begin{lemma} \label{lemma:RH} The de Rham functor:
\[
\bM \to \op{cone}\left( \xymatrix@1{\bM\otimes_{\mathbb{C}[u]}
\mathcal{O}^{\op{an}}_{\mathbb{A}^{1}} \ar[r]^-{\frac{\partial}{\partial u}}
 & \bM\otimes_{\mathbb{C}[u]} 
\mathcal{O}^{\op{an}}_{\mathbb{A}^{1}}} \right)
\] 
establishes an equivalence between the categories {\em ({\bfseries
    \nc dR(iii)})} 
and {\em ({\bfseries \nc B(i)})}$\otimes \mathbb{C}$.
\end{lemma}

\

\medskip

\begin{figure}[!ht]
\begin{center}
\psfrag{0}[c][c][1][0]{{$0$}}
\psfrag{el}[c][c][1][0]{{$\lambda e^{i\varphi}$}}
\psfrag{s}[c][c][1][0]{{$\lambda
    \geq 0$}}
\psfrag{ss}[c][c][1][0]{{$\lambda
    \leq 0$}}
\psfrag{Hlphi}[c][c][1][0]{{$\mathfrak{H}_{\varphi,\lambda}$}}
{\protect \epsfig{file=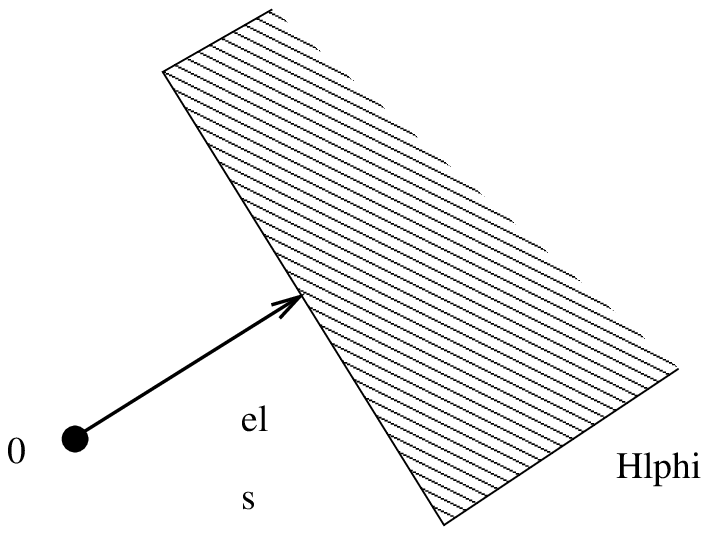,width=2in}} \qquad or \qquad 
{\protect \epsfig{file=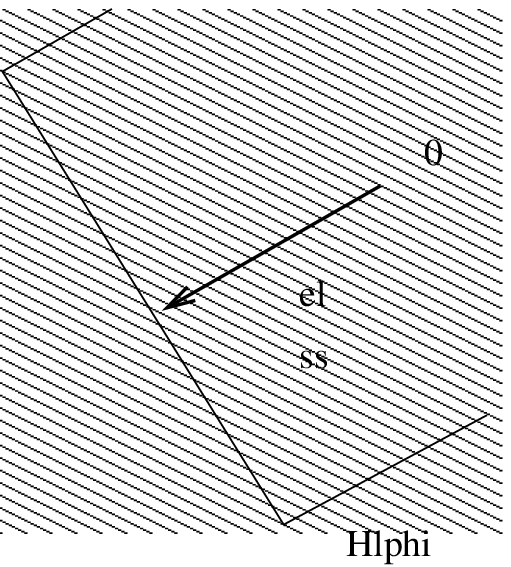,width=1.45in}}
\end{center}
\caption{The half-plane
  $\mathfrak{H}_{\varphi,\lambda}$.}\label{fig:halfplanes}  
\end{figure}

\noindent
Finally, note that Theorem~\ref{prop:B(i)=B(ii)}, together with
Lemma~\ref{lemma:RH}, and Deligne's classification
\cite[Theorem~4.7.3]{babbitt-varadarajan-memoir} of germs of irregular
connections give immediately:

\begin{lemma} The data data {\em ({\bfseries \nc
  B(ii)})} and {\em ({\bfseries \nc B(iv)})} are equivalent. 
\end{lemma}
{\bfseries Proof.} 
Let $\mycal{F}$ be  a constructible sheaf of
$\mathbb{Q}$-vector spaces on $\mathbb{C}$.  Define a local system
$\stS$ of $\mathbb{Q}$-vector spaces on $\bS^{1}$ as the restriction
of $\mycal{F}$ to the circle ``at infinity'', i.e. define the stalk of
$\stS$ at $\varphi \in \bS^{1}$ to be
\[
\stS_{\varphi} := \lim_{r \to +\infty} \mycal{F}_{re^{i\varphi}}.
\] 
Next, for any $\lambda \in \mathbb{R}$ and any $\varphi \in \bS^{1}$
consider the  half-plane 
\[
\mathfrak{H}_{\varphi,\lambda} := \left(\lambda + \left\{ u \in
\mathbb{C} | \ \op{Re}(u) \geq 0 \right\}\right)\cdot e^{i\varphi},
\]
as shown on  Figure~\ref{fig:halfplanes}. 

\

\noindent
Now suppose that $R\Gamma(\mathbb{C},\mycal{F}) = 0$. By the long
exact sequence for the cohomology of the pair
$\mathfrak{H}_{\varphi,\lambda} \subset \mathbb{C}$ we get that
$H^{i}(\mathbb{C},\mathfrak{H}_{\varphi,\lambda};\mycal{F}) = 0$
unless $i = 1$. The Deligne-Malgrange-Stokes filtration on $\stS$ is
then given explicitly by
\[
\stS_{\varphi,\leq \lambda} :=
H^{1}(\mathbb{C},\mathfrak{H}_{\varphi,\lambda};\mycal{F}) \subset
\stS_{\varphi}. 
\]
\  \hfill $\Box$

\

\medskip

\noindent
For the purposes of \nc-Hodge theory all these statements can be
summarized in the following

\begin{theo} \label{theo:summary} There is natural equivalence of
  categories 
\[
\left( 
\begin{minipage}[c]{2.1in} 
{\em
triples $(H,\mycal{E}_{B},\iso)$ satisfying \linebreak
the  
{\sf\bfseries
  (\nc-filtration axiom)$^{\bexp}$} 
and 
the 
{\sf\bfseries
  ($\mathbb{Q}$-structure axiom)$^{\bexp}$} 
}
\end{minipage}
\right)
\leftrightarrow
\left( 
\begin{minipage}[c]{3.3in} 
{\em quadruples $((H,\nabla),\mycal{F}_{B},\bbf)$, where
\begin{itemize}
\item $H$ is an algebraic $\mathbb{Z}/2$-graded vector bundle on
  $\mathbb{C}$ and $\nabla$ is a meromorphic connection on $H$
  satisfying the \linebreak {\sf\bfseries (\nc-filtration axiom)$^{\bexp}$};
\item $\mycal{F}_{B} \in \op{\sf{Constr}}(\mathbb{C},\mathbb{Q})$,
  satisfying $R\Gamma(\mathbb{C},\mycal{F}_{B}) = 0$;
\item $\bbf$ is an isomorphism  
\[
\bbf : \mycal{F}_{B}\otimes \mathbb{C}  \to DR\left(
\boldsymbol{\Phi}\left[\iota_{*} \left((H,\nabla)_{|\mathbb{A}^{1} - \{ 0
    \}}\right)\right]\right)  
\]
in $D^{b}_{\op{constr}}(\mathbb{C},\mathbb{C})$
\end{itemize}
}
\end{minipage}
\ 
\right)
\]
Here as before 
\begin{description}
\item[$DR$] is the de Rham complex functor from the derived category
of regular holonomic
$\mathbb{D}$-modules to the derived category of constructible
sheaves,
\item[$\iota$] is the inclusion map  $\iota : \mathbb{A}^{1} - \{ 0\}
  \hookrightarrow \mathbb{A}^{1}$ given by $\iota(v) = v^{-1}$, and 
\item[$\boldsymbol{\Phi}(\bullet)$] is the
  Fourier-Laplace transform for $\mathcal{D}$-modules on
  $\mathbb{A}^{1}$.
\end{description}
\end{theo}
{\bfseries Proof.} Follows immediately from previous
equivalences. \ \hfill $\Box$

\subsection{\bfseries Structure results} \label{subsec:structure}

In this section we collect a few results clarifying the structure
properties of the \nc-Hodge structures of exponential type.

\subsubsection {\bfseries A quiver description of \nc-Betti data} \ Since
the gluing data {\bfseries \nc B(iii)} are of essentially combinatorial
nature, it is natural to look for a quiver interpretation of this
data. To that end consider the algebra 
\begin{equation}
\boldsymbol{\mycal{A}}_{n} := 
\left\langle 
\text{
\begin{minipage}[c]{1.1in}
$\bp_{1}$,\ldots, $\bp_{n}$  \\
$T$,$T_{11}^{-1}$,\ldots,$T_{nn}^{-1}$
\end{minipage}}
\left| 
\text{
\begin{minipage}[c]{2.1in}
$\bp_{1}+\bp_{2}+\ldots+\bp_{n} = 1$\\
$\bp_{i}\bp_{j}=\bp_{j}\bp_{i}$ for $i\neq j$, $\bp_{i}^{2} = \bp_{i}$\\
$T_{ii}^{-1}\bp_{i}T\bp_{i} = \bp_{i}T\bp_{i}T_{ii}^{-1} = \bp_{i}$
\end{minipage}}
\right. \right\rangle 
\end{equation}
This is the path algebra of the complete quiver having $n$ ordered
vertices, $n^{2}-n$ arrows connecting all pairs of distinct vertices,
and $2n$-loops - two at each vertex, with the only relations being that
the two loops at every given vertex  are inverses to each other.

Note that our description of the gluing data {\bfseries \nc B(iii)}
now immediately gives the following
\begin{lemma} For a given set of points $S = \{ \bc_{1}, \ldots, \bc_{n}
  \} \subset \mathbb{C}$, the category of gluing data {\bfseries \nc
    B(iii)} with singularities at $S$  is equivalent to the category
  of finite dimensional representations of $\mycal{A}_{n}$. 
\end{lemma}

In particular since the braid group $\boldsymbol{B}_{n}$ on
$n$-strands acts naturally on the data {\bfseries \nc B(iii)} we get a
homomorphism $\boldsymbol{B}_{n} \to \op{Aut}(\mycal{A}_{n})$ from the
braid group to the group of algebra automorphisms of $\mycal{A}_{n}$.

\

\medskip

\subsubsection {\bfseries Gluing of \nc-Hodge
  structures} \label{sssec:gluenc} \ It is natural
to expect that the usual classification of connections with second order
poles in terms of formal regular type and Stokes multipliers can be
promoted to a similar classification of \nc-Hodge
structures. The search for such a classification leads naturally to
the following theorem:

\begin{theo} \label{theo:gluenc} Let  $\left\{\left(H,\mycal{E}_{B},
  \iso\right)\right\}$ be a \nc-Hodge structure 
 of exponential type.  Then specifying $\left\{\left(
  H, \mycal{E}_{B},\iso\right)\right\}$ is equivalent
  to specifying the following data:
\begin{description}
\item[{\bfseries (regular type)}]  A finite set $S =
  \{\bc_{1},\ldots,\bc_{n}\} \subset 
  \mathbb{C}$ and a collection $\left\{\left(
  (\mycal{R}_{i},\nabla_{i}),
  \mycal{E}_{B,i},\iso_{i}\right)\right\}_{i = 1}^{n}$ of \nc-Hodge
  structures with regular singularities.
\item[{\bfseries (gluing data)}] A base point $\bc_{0} \in
  \mathbb{C}-S$, a collection of discs $\{\bD_{i}\}_{i = 1}^{n}$ and
  paths $\{a_{i}\}_{i = 1}^{n}$, chosen as in the proof of
  Theorem~\ref{prop:B(i)=B(ii)}, and for every $i \neq j$, $i, j
  \in \{1, \ldots, n\}$ a map of rational vector spaces 
\[
T_{ij} : \left(\mycal{E}_{B,j}\right)_{\bc_{0}} \longrightarrow
\left(\mycal{E}_{B,i}\right)_{\bc_{0}} 
\]
\end{description}
\end{theo}
{\bf Proof.} It will be convenient to introduce formal counterparts
to the de Rham parts of the \nc-Hodge structures appearing in the
statement of the theorem. We consider the following:

 \

\smallskip

\noindent
{\bf formal(a)} A pair $(\mycal{M}^{\op{for}},\nabla^{\op{for}})$,
where $\mycal{M}^{\op{for}}$ is a finite dimensional vector space
over $\mathbb{C}((u))$ and
$\nabla^{\op{for}}$ is a meromorphic connection on
$\mycal{M}^{\op{for}}$ of exponential
  type. 

\

\medskip

\noindent
{\bf formal(b)} A finite set of points $S = \{\bc_{1},\ldots,\bc_{n}\}
\subset \mathbb{C}$ and a collection
$\{(\mycal{R}_{i}^{\op{for}},\nabla_{i}^{\op{for}})\}_{i=1}^{n}$ where
each $\mycal{R}_{i}^{\op{for}}$ is a non-zero finite dimensional
vector space over $\mathbb{C}((u))$ and each $\nabla_{i}^{\op{for}}$
is a meromorphic connection on $\mycal{R}_{i}^{\op{for}}$ with a
regular singularity.

\

\medskip

\noindent
{\bf formal(c)} A finite collection of points $S =
\{\bc_{1},\ldots, \bc_{n} \} \subset \mathbb{C}$, and 
\begin{itemize}
\item a collection $U_{1}, U_{2}, \ldots, U_{n}$ of finite dimensional
non-zero $\mathbb{Q}$-vector spaces,
\item a collection of linear maps $T_{ii} \in GL(U_{i})$, for all
  $i = 1, \ldots, n$,
\end{itemize} 

\

\bigskip

\noindent
By Remark~\ref{rem-decompose-with-lattices} the natural
functor from the category of data {\bf 
  formal(b)} to the category of data {\bf formal(a)}, which is given
by 
\[
\xymatrix@R-1pc{
(\text{{\bf  formal(b)}}) \ar[r] & (\text{{\bf  formal(a)}}) \\
\left(S;
\{(\mycal{R}_{i}^{\op{for}},\nabla_{i}^{\op{for}})\}_{i=1}^{n} \right)
\ar[r] & 
\bigoplus_{i=1}^{n}
\be^{\bc_{i}/u}\otimes (\mycal{R}_{i}^{\op{for}},\nabla_{i}^{\op{for}})
=: (\mycal{M}^{\op{for}},\nabla^{\op{for}}) 
}
\]
is an equivalence of categories.

Also we have the following 

\

\begin{lemma} \label{lemma:formal(b)=formal(c)} The categories of data 
{\bf  formal(b)} and {\bf  formal(c)} are naturally equivalent.
\end{lemma}
{\bf Proof.} Indeed, consider the category $\mathcal{C}$ of all data
consisting of
a finite set of points \linebreak $S = \{\bc_{1},\ldots,\bc_{n}\}
\subset \mathbb{C}$ and a collection
$\{(\mycal{R}_{i},\nabla_{i})\}_{i=1}^{n}$ where
each $\mycal{R}_{i}$ is a non-zero finite dimensional
vector space over $\mathbb{C}\{u\}[u^{-1}]$ and each $\nabla_{i}$
is a meromorphic connection on $\mycal{R}_{i}$ with a
regular singularity and non-trivial monodromy. Then we have natural
functors 
\[
\xymatrix{
& (\text{{\bf  formal(b)}}) \\
\mathcal{C} \ar[ru]^-{(\bullet)\otimes \mathbb{C}((u))}
\ar[rd]_-{\op{\text{{\bfseries\sf mon}}}} & \\
& (\text{{\bf  formal(c)}})
}
\]
where $(\bullet)\otimes \mathbb{C}((u))$ is the passage to a formal
completion and $\op{\text{{\bfseries\sf mon}}}$ is given by assigning to each
$(\mycal{R}_{i},\nabla_{i})$ the pair $(U_{i},T_{i})$, where $U_{i}$
is the fiber of the Birkhoff extension
$\mathfrak{B}_{0}(\mycal{R}_{i},\nabla_{i})$ of
$(\mycal{R}_{i},\nabla_{i})$ at $1 \in \mathbb{A}^{1}$, and $T_{i}$ is
the monodromy of $\mathfrak{B}_{0}(\mycal{R}_{i},\nabla_{i})$ around
the unit circle traced in the positive direction.

This proves the lemma since 
 $\op{\text{{\bfseries\sf mon}}}$ is an equivalence by the
Riemann-Hilbert correspondence and $(\bullet)\otimes \mathbb{C}((u))$
is an equivalence by the formal decomposition theorem
\cite[II.5.7]{sabbah-frobenius}). \\ \ \mbox{\ } \hfill $\Box$

\

Note that these equivalences are compatible with the corresponding
equivalence of analytic de Rham data and Betti data. More precisely we
have a commutative diagram of functors
\begin{equation} \label{eq:anfor}
\xymatrix{ 
(\text{{\bfseries \nc dR(i)}})
  \ar[r]
  \ar@{<.>}[d] & (\text{{\bfseries formal(a)}}) \ar@{<.>}[d] \\
(\text{{\bfseries \nc B(iii)}}) \ar[r] & (\text{{\bfseries formal(c)}})
}
\end{equation}
Here the right vertical equivalence is the composition of the
equivalences $(\text{{\bfseries formal(a)}}) \cong (\text{{\bfseries
    formal(b)}}) \cong (\text{{\bfseries formal(c)}})$ that we just
discussed. The left vertical equivalence is the composition of the
equivalence $(\text{\nc dR(i)}) \cong (\text{\nc dR(iii)})$ given in
Lemma~\ref{lem-dr-equivalence}, the equivalence $(\text{\nc dR(iii)})
\cong (\text{\nc B(i)})$ from Lemma~\ref{lemma:RH}, the equivalence
$(\text{\nc B(i)}) \cong (\text{\nc B(ii)})$ given in
Theorem~\ref{prop:B(i)=B(ii)}, and the equivalence 
 $(\text{\nc B(ii)}) \cong (\text{\nc B(iii)})$ from
Lemma~\ref{lemma:B(ii)=B(iii)}.

Horizontally we have the forgetful functors
\[
\xymatrix@R-2pc{ 
(\text{{\bfseries \nc dR(i)}})
  \ar[r]
& (\text{{\bfseries formal(a)}}) \\ 
(\mycal{M},\nabla)
  \ar@{|->}[r]
& (\mycal{M},\nabla)
\otimes_{\mathbb{C}\{u\}[u^{-1}]}\mathbb{C}((u)),}
\]
and
\[
\xymatrix@R-2pc{ 
(\text{{\bfseries \nc B(iii)}})
  \ar[r]
& (\text{{\bfseries formal(c)}}) \\
\left(S;\{U_{i}\}_{i=1}^{n},\{T_{ij}\}_{i,j=1}^{n}\right)
  \ar@{|->}[r]
& \left(S;\{U_{i}\}_{i=1}^{n},\{T_{ii}\}_{i=1}^{n}\right).}
\]
\

\noindent
Next we need the following 

\

\begin{lemma} Suppose that $(\mycal{M},\nabla)$ is some de Rham
  data of type {\em{\bfseries \nc dR(i)}} and let
\[
(\mycal{M}^{\op{for}},\nabla^{\op{for}}) =
  (\mycal{M},\nabla)
\otimes_{\mathbb{C}\{u\}[u^{-1}]}\mathbb{C}((u))
\]
  be the corresponding formal data. Then:
\begin{itemize}
\item[{\em\bfseries (a)}] the map
\[
\xymatrix@1@C+2pc{\left( \text{\begin{minipage}[c]{2.3in}
      $\mathbb{C}\{u\}$-submodules  
$\mycal{H} \subset \mycal{M}$, 
    on which $\nabla$ has a pole of order $\leq 2$ \end{minipage}}\right)
\ar[r]^-{(\bullet)\otimes \mathbb{C}[[u]]} & 
\left( \text{\begin{minipage}[c]{2.3in} $\mathbb{C}[[u]]$-submodules 
$\mycal{H}^{\op{for}} \subset \mycal{M}^{\op{for}}$, 
    on which $\nabla^{\op{for}}$ has a pole of order $\leq
    2$ \end{minipage}}\right),
}
\]
is bijective.
\item[{\em\bfseries (b)}] If $\Psi :
  (\mycal{M}^{\op{for}},\nabla^{\op{for}}) \to \bigoplus_{i=1}^{n}
\be^{c_{i}/u}\otimes (\mycal{R}_{i}^{\op{for}},\nabla_{i}^{\op{for}})$
is a formal isomorphism, then the map
\[
\xymatrix@1@C+2pc{ \left( \text{\begin{minipage}[c]{2.3in}
      $\mathbb{C}[[u]]$-submodules $\mycal{H}^{\op{for}}  \subset
      \mycal{M}^{\op{for}}$, on which $\nabla^{\op{for}}$ has a pole
      of order $\leq 2$ \end{minipage}}\right) & \ar[l]^-{\Psi}  
\left(
  \text{\begin{minipage}[c]{2.3in} $\mathbb{C}[[u]]$-submodules
      $\mycal{H}_{i}^{\op{for}}  \subset \mycal{R}^{\op{for}}_{i}$, for
      all $i = 1, \ldots, n$, on which $\nabla^{\op{for}}_{i}$ has a
      pole of order $\leq 2$ \end{minipage}}\right),}
\]
\end{itemize}
is bijective.
\end{lemma}
{\bf Proof.} {\bfseries (a)} \ 
Pick some frame $\underline{e}$ of $\mycal{M}$
over $\mathbb{C}\{u\}[u^{-1}]$ and let $\mycal{H}^{0} :=
\mathbb{C}\{u\}\cdot \underline{e} \subset \mycal{M}$ be the
submodule of all sections in $\mycal{M}$ that are
holomorphic in this frame. Now any $\mathbb{C}\{u\}$-submodule  
$\mycal{H} \subset \mycal{M}$, 
    on which $\nabla$ has a pole of order $\leq 2$ will be a 
    $\mathbb{C}\{u\}$-submodule of  $\mycal{M}$ which is
    commensurable with $\mycal{H}^{0}$, i.e. we will have 
$u^{N}\mycal{H}^{0} \subset \mycal{H} \subset u^{-N}\mycal{H}^{0}$ for
    $N \gg 1$. However the formal completion functor
    $(\bullet)\otimes_{\mathbb{C}\{u\}} \mathbb{C}[[u]]$  establishes an
    isomorphism between the Grassmanian
    $GL_{r}(\mathbb{C}\{u\}[u^{-1}])/GL_{r}(\mathbb{C}\{u\})$ and the
    affine Grassmanian
    $GL_{r}(\mathbb{C}((u))/GL_{r}(\mathbb{C}[[u]])$. But this
    map preserves the condition that a submodule $\mycal{H}$ is
    invariant under
    $\nabla_{u^{2}d/du}$ which proves {\bfseries (a)}.  

\

\noindent
{\bfseries (b)} As already mentioned in
Remark~\ref{rem-decompose-with-lattices} this is proven in 
\cite[Lemma~8.2]{hertling-sevenchek}. Alternatively we can reason as
in the proof of part {\bfseries (a)}.  Let 
 $\mycal{H} $ be a 
$\mathbb{C}[[u]]$-submodule in $\mycal{M}^{\op{for}}$ which is 
commensurable with $\mycal{H}^{0,\op{for}} $ and  preserved by
$\nabla_{u^{2}\frac{d}{du}}$.  The operator 
$\nabla_{u^{2}\frac{d}{du}}$ acts on the infinite-dimensional
 topological complex vector space $\mycal{M}^{\op{for}}$ with
 finitely many infinite Jordan blocks with eigenvalues 
$\{\bc_1,\dots,\bc_n\}$. The corresponding generalized eigenspaces are
 exactly modules $\be^{\bc_i/ u} \mycal{R}^{\op{for}}_{i}$. Hence 
\[
\mycal{H}^{\op{for}} =\oplus_i 
\left(\mycal{H}^{\op{for}} \cap
\be^{\bc_i/u}\mycal{R}^{\op{for}}_{i}\right)
\]
Therefore we obtain extensions $\mycal{R}^{\op{for}}_{i}$ with second
order poles and regular singularity. 

\

\noindent
Combining the previous lemma with the equivalences in diagram
\eqref{eq:anfor} and the description of \nc-Hodge structures from
Section~\ref{sss:irregular} gives the theorem. \ \hfill $\Box$

\subsection{\bfseries Deformations of \nc-spaces and gluing} 
\label{subsec:gluingdefos}

In this section we will briefly examine how the gluing construction
for \nc-Hodge structures varies with parameters. In particular, we
will look at deformations of \nc-spaces and the way the gluing data
for the \nc-Hodge structures on the cohomology of these spaces
interacts with the appearance of a curvature in the d$(\mathbb{Z}/2$g
algebra computing the sheaf theory of the space.

\subsubsection{\bfseries The cohomological Hochschild complex}
\label{sssec:HH.cochains} 

Suppose $X = \ncSpec A$ is a \nc-affine \nc-space. Recall that the
cohomological Hochschild complex is defined as 
\[C^{\bullet}(A,A):=  \prod_{n \geq 0}
\op{Hom}_{\mathbb{C}-\text{\sf{Vect}}}\left((\bPi A)^{\otimes n},A\right),
\]
Its shift $\bPi C^{\bullet}(A,A)$ is a Lie superalgebra with 
respect to the Gerstenhaber bracket \cite{murray}, and can be
interpreted as 
 the Lie algebra of continuous derivations of the free 
topological algebra
  $\prod_{n\ge 0} \left((\bPi A)^{\otimes n}\right)^\vee$.
The multiplication $m_{A}$ and differential $d_{A}$ of
$A$ combine into a cochain \linebreak $\gamma_A:= m_{A} + d_{A} \in
C^{\bullet}(A,A)$  satisfying $[\gamma_A,\gamma_A]=0$.

The formal
deformation theory of $X$ is controlled by a d$(\mathbb{Z}/2)$g Lie
algebra structure $\bPi C^{\bullet}(A,A)$ endowed with the differential
 $[\gamma_A,\bullet]$. It is convenient to consider also the reduced
Hochschild complex 
 \[C^{\bullet}_{\op{red}}(A,A):=  \prod_{n \geq 0}
\op{Hom}_{\mathbb{C}-\text{\sf{Vect}}}
\left(\left(\bPi \left(A/\mathbb{C}\cdot 
1_A\right)\right)^{\otimes n},A\right),
\]
which is naturally a subspace of $C^{\bullet}(A,A)$. The reduced
complex is (after the parity change) a dg Lie subalgebra in $\bPi
C^{\bullet}(A,A)$. Moreover it is quasi-isomorphic to $\bPi
C^{\bullet}(A,A)$. Hence, for deformation theory purposes one can
replace $\bPi C^{\bullet}(A,A)$ by $ \bPi
C^{\bullet}_{\op{red}}(A,A)$.

Let $\bbgamma = \sum_{i \geq 1} \gamma_{i} t^{i}  \in t
C^{\op{even}}_{\op{red}}(A,A)[[t]]$ be a formal path consisting of
 solutions of  the
Maurer-Cartan equation, i.e.
\[d\bbgamma + \frac{1}{2}
\left[\,\bbgamma, \bbgamma \,\right] = 0 \quad 
\left(\, \Leftrightarrow \;
\left[ \bbgamma + \gamma_{A}, \bbgamma + \gamma_{A} \right] = 0 \,
\right).
\]
Such a solution defines so called
 formal deformation of the d$(\mathbb{Z}/2)$g
algebra $A$ as a weak (or curved) $A_{\infty}$-algebra (see 
e.g. \cite{hasegawa} for the definition and \cite{schwartz-curved} for
 a more detailed analysis). We can use the cochain
 $\bbgamma + \gamma_{A} \in  C^{\op{even}}(A,A)[[t]]$   to twist the notion of
an $A$-module.  We will write $A_{\bbgamma}$ for the (weak)
$A_{\infty}$-algebra over $\mathbb{C}[[t]]$ corresponding to $A$ and
$\bbgamma + \gamma_{A}$ and $(A_{\bbgamma}-\text{\sf mod})$ for the
$\mathbb{C}[[t]]$-linear dg category of all modules over
$A_{\bbgamma}$. By definition  $(A_{\bbgamma}-\text{\sf mod})$ is the
category of dg modules over a bar-type resolution of $A_{\bbgamma}$
\cite{ks-ncgeometry}. As an algebra the relevant bar dg algebra is the
completed tensor product
\begin{equation} \label{eq:bar.algebra}
\prod_{n\geq 0} \left(
(\bPi A)^{\otimes n}\right)^{\vee} 
\widehat{\otimes} \mathbb{C}[[t]]
\end{equation}
where the algebra structure comes from the usual algebra structure on
$\mathbb{C}[[t]]$ and the tensor algebra structure on $\prod_{n\geq 0}
\left( (\bPi A)^{\otimes n}\right)^{\vee}$. Thus for every $\bbgamma \in t
C^{\op{even}}_{\op{red}}(A,A)[[t]]$ which solves the Mauer-Cartan
equation we get a differential $\bbgamma + \gamma_{A}$ on the graded
algebra \eqref{eq:bar.algebra}. The bar dg algebra of $A_{\bbgamma}$
is now defined as the dg algebra 
\[
B_{\bbgamma} := \left(\prod_{n\geq 0}
\left((\bPi A)^{\otimes n}\right)^{\vee}
\widehat{\otimes} \mathbb{C}[[t]],\bbgamma + \gamma_{A}\right).
\] 
The
dg category $\left(A_{\bbgamma}-\text{{\sf mod}}\right)$ is by
definition the category of dg modules over $B_{\bbgamma}$ which are
topologically free as modules of the underlying algebra, i.e. after
forgetting the differential, and also satisfying the condition of
unitality at $t=0$.

As before this category can be viewed as the category
$C_{\mathbb{X}_{\bbgamma}} := (A_{\bbgamma}-\text{\sf mod})$ of
quasi-coherent sheaves on a \nc-affine \nc-space
$\mathbb{X}_{\bbgamma} \to \mathbb{D}$ defined over the formal disc
$\mathbb{D} = \op{Spf}(\mathbb{C}[[t]])$.  More generally we will get
a \nc-space $\mathbb{X}$ over the formal scheme of solutions to the
Maurer-Cartan equation  and $\mathbb{X}_{\bbgamma} \to \mathbb{D}$ is
the restriction of $\mathbb{X}$ to the formal path $\bbgamma +
\gamma_{A}$ sitting inside that formal scheme.

Similarly we can use $\bbgamma$ to twist the notion of a Hochschild
cohomology class for $A$. Namely we can consider the Hochschild cohomology of
the $A_{\infty}$-algebra $A_{\bbgamma}$. It is given explicitly as the
cohomology 
\[
HH^{\bullet}_{\bbgamma}(A) :=
H^{\bullet}\left(C^{\bullet}(A,A)[[t]], \left[ \,
  \bbgamma + \gamma_{A}, \bullet \right]\right),
\] 
and is a commutative algebra with respect to the cup product.
Note also that the algebra  $HH^{\bullet}_{\bbgamma}(A)$ comes
equipped with a unit $\left[1_{A}\right]$ and a distinguished even element
$\left[\bbgamma+\gamma_{A} \right]$, i.e. a structure similar to the one
discussed in Section~\ref{sssec:mero.u}.

\begin{rem} \label{rem:nonzero.objects} $\bullet$ \ If $\bbgamma$ has
  no component of degree zero, i.e. if 
\[
\bbgamma \in t
  C^{\op{even}}_{\op{red},+}(A,A)[[t]],  \text{ where }
  C^{\bullet}_{\op{red},+}(A,A) = \prod_{n \geq 1}
\op{Hom}_{\mathbb{C}-\text{\sf{Vect}}}\left(\left(\bPi \left(A/\mathbb{C}\cdot
1_{A}\right)\right)^{\otimes n},A\right),
\]
then $A_{\bbgamma}$ is an honest (strong) $A_{\infty}$-algebra, and
the category $(A_{\bbgamma}-\text{\sf mod})$ will typically have many
interesting objects. Furthermore, in this case smoothness and
compactness are stable under deformations. That is, if $A$ is smooth
(respectively compact) over $\mathbb{C}$, then $A_{\bbgamma}$ is
smooth (respectively compact) over $\mathbb{C}[[t]]$.

\

\noindent
$\bullet$ \ If the $n = 0$ component of $\bbgamma$ is  non-trivial,
i.e. if the corresponding $A_{\infty}$ structure has a non-trivial
$m_{0}$, then  the category $(A_{\bbgamma}-\text{\sf mod})$
may have no non-zero objects. The basic example of this is when $A =
\mathbb{C}$ and $\bbgamma = t\cdot 1_{A}$.
\end{rem}

\smallskip

\

\noindent
If the original algebra $A$ has the degeneration property,
then it is easy to see that the Hodge-to-de Rham spectral sequence
will degenerate for the periodic cyclic homology of $A_{\bbgamma}$.
In other words the formal \nc-space $\mathbb{X}$ will give rise to a
variation of \nc-Hodge structures over the formal scheme of
solutions of the Maurer-Cartan equation for $A$.  When we have a
non-trivial $n = 0$ component in $\bbgamma$ this may lead to a
paradoxical situation in which we have a family of \nc-spaces over
$\mathbb{D}$ which has no sheaves over the generic point but has 
non-trivial de Rham cohomology (i.e. periodic cyclic homology)
generically. This suggests the following important

\begin{que} \label{que:dissapear} What is the geometrical meaning of
  $HH^{\bullet}_{\bbgamma}(A)$,  $HH_{\bullet}(A_{\bbgamma})$,
  $HH_{\bullet}^{-}(A_{\bbgamma})$, and $HP_{\bullet}(A_{\bbgamma})$, when
$\bbgamma$ has non-trivial $n = 0$ component and the objects of
  $(A_{\bbgamma}-\text{\sf mod})$ dissapear over $\mathbb{D}^{\times}$?
\end{que}

\

\begin{rem} \label{rem:add.one} Note that if $\bbgamma$ solves the
  Maurer-Cartan equation, then for any $c \in t\mathbb{C}[[t]]$, the
  cochain $\bbgamma + c\cdot 1_{A}$ will also solve the 
Maurer-Cartan equation\footnote{In fact this is the main reason for
  all the hassle with the unit and the reduced complex in this
  section.}. So we have a natural mechanism for modifying 
  formal paths of solutions of the Maurer-Cartan equation. We will
  exploit this mechanism in the next section.
\end{rem}

\

\subsubsection{{\bfseries Corrections by constants}}
\label{sssec:localizations} 

The unpleasant phenomenon of having \nc-spaces with no sheaves and
non-trivial cohomology at the generic point is related to the gluing
description for \nc-Hodge structures. The idea is that the $A$-modules
that dissapear at the generic point of $\mathbb{D}$ may reappear again
if we modify the weak $A_{\infty}$ algebra $A_{\bbgamma}$
appropriately. The periodic cyclic homologies of the different
admissible modifications of $A_{\bbgamma}$ then correspond to the
regular pieces in the gluing description of the \nc-de Rham data 
given by $HP_{\bullet}(A_{\bbgamma})$. More precisely we have the
following

\begin{conn} \label{conn:localizations} Suppose that $A$ is a smooth
  and compact d$(\mathbb{Z}/2)$g algebra. Let $\bbgamma  \in t
  C^{\op{even}}_{\op{red}}(A,A)[[t]]$  be a formal
  even path of solutions of the Maurer-Cartan equation for $A$. Then
  the periodic cyclic homology $HP_{\bullet}\left(A_{\bbgamma}\right)$
  carries a canonical functorial structure of a variation of
  $\mathbb{Q}$-\nc-Hodge structures of exponential type  over
  $\mathbb{D} = \op{Spf}(\mathbb{C}[[t]])$. 
  Furthermore 
  there exists a positive integer $N$ and a finite collection of
  pairwise distinct Puiseux
  series  
\[
\bc_{i} = \sum_{j \ge 1} c_{i,j} t^{\frac{j}{N}}, \quad c_{i,j} \in \mathbb{C}
\]
such that:
\begin{itemize}
\item The series $\bc_i$ are the 
distinct eigenvalues of the operator of mutiplication
   by the class $\left[\bbgamma + \gamma_A\right]$ in the
   supercommutative algebra 
  $HH^\bullet_{\bbgamma}(A)\widehat{\otimes}_{\mathbb{C}[[t]]}
   \mathbb{C}((t))$. 
\item For each $i$ the category $\left( A_{\bbgamma + \bc_{i}\cdot
  1_{A}}-\text{{\sf mod}}\right)$ is a non-trivial
  $\mathbb{C}\left[\left[t^{1/N}\right]\right]$-linear
  d$(\mathbb{Z}/2)$g category which 
  are smooth and compact over
  $\mathbb{C}\left[\left[t^{1/N}\right]\right]$ and is computed by a
  d$(\mathbb{Z}/2)$g algebra $B_{i}$ defined over
  $\mathbb{C}\left[\left[t^{1/N}\right]\right]$ and quasi-isomorphic
  to the (weak) $A_{\infty}$-algebra $A_{\bbgamma + \bc_{i}\cdot
  1_{A}}$.  
\item The Hochschild homologies $HH_{\bullet}(B_{i})$ are flat
  $\mathbb{C}\left[\left[t^{1/N}\right]\right]$-modules and we have 
\[
\sum_{i} \op{rk}_{\mathbb{C}\left[\left[t^{1/N}\right]\right]}
(HH_{\bullet}(B_{i})) =
\op{rk}_{\mathbb{C}\left[\left[t^{1/N}\right]\right]}
  HH_{\bullet}\left(A_{\bbgamma}\right) = \dim_{\mathbb{C}} HH_{\bullet}(A).
\]
\item The variation of \nc-Hodge structures $HP_{\bullet}\left(
  A_{\bbgamma} \right)$ viewed as a variation over
  $\mathbb{C}\left[\left[t^{1/N}\right]\right]$ has as regular
  constituents  the variations of \nc-Hodge 
  structures on $HP_{\bullet}(B_{i})$ whose existence is predicted by
  Conjecture~\ref{con:main.conjecture}.  
\end{itemize}
\end{conn}

\

\noindent
In particular Conjecture~\ref{conn:localizations} says that the
categorical and Hodge theoretic content of  the algebra
$A_{\bbgamma}$ consists of the following data:
\begin{description}
\item[{\bfseries (categories)}] A finite collection of smooth and
  compact $\mathbb{C}\left[\left[t^{1/N}\right]\right]$-linear 
 d$(\mathbb{Z}/2)$g categories $(B_{i}-\text{{\sf mod}})$.
\item[{\bfseries (gluing)}] A finite collection of distinct Piuseux
  series $\bc_{i} \in \mathbb{C}\left[\left[t^{1/N}\right]\right]$, and
  formal \nc-gluing data which glues the variations of regular \nc-Hodge
  structures on $HP_{\bullet}(B_{i})$ into a variation of \nc-Hodge
  structure of 
  exponential type over $\mathbb{C}\left[\left[t^{1/N}\right]\right]$.
\end{description}

\

\noindent
In the above discussion we have tacitly replaces the analytic setting
from Section~\ref{sec-gluing} 
by a formal setting. One can check that both the de Rham and Betti data
make sense here, e.g.  one can speak about homotopy classes of
non-intersecting paths to points $\bc_i$ thinking about $t$ as a small
real positive parameter.

\

\begin{rem} \label{rem:Milnor} This situation is analogous to a well
  known setup in singularity theory. Namely, if we have a germ of an
  isolated hypersurface singularity given by an equation $f = 0$, and
  if we have a deformation of $f$ which has several critical values,
  then the Milnor number of the original singularity is equal to the
  sum of the Milnor numbers of the simpler critical points of the
  deformed function. In fact, as we will see in
  section~\ref{subsec:B.model}  
  the singularity setup is a rigorous
  manifestation of the above conjectural picture.
\end{rem}

\subsubsection {{\bfseries Singular deformations}} 
\label{sssec:singular.defos} Suppose next that $A$ is compact but not
  smooth (or smooth but non-compact) d$(\mathbb{Z}/2)$g algebra and
  let again $\bbgamma \in t
  C^{\op{even}}_{\op{red}}(A,A)[[t]]$ be a formal path of solutions of
  the Maurer-Cartan equation. We expect that the usual definition of
  smoothness and compactness can be modified to give a notion of
  smoothness together with compactness of $A_{\bbgamma}$ at the
  generic point, i.e. over $\mathbb{C}((t))$, even when the objects in
  $\left( A_{\bbgamma}-\text{{\sf mod}}\right)$ dissapear over
  $\mathbb{C}((t))$.

In the case when $A_{\bbgamma}$ is smooth and compact over
$\mathbb{C}((t))$, i.e. when the deformation given by $\bbgamma$ is a
smoothing deformation, we also expect
Conjecture~\ref{conn:localizations} to hold at the generic point. More
precisely, we expect to have Puiseux series $\bc_{i}$ as above for which
the associated categories $\left( A_{\bbgamma +\bc_{i}1_{A}}-\text{{\sf
mod}}\right)$ are non-trivial and smooth and compact over
$\mathbb{C}\left(\left( t^{1/N} \right)\right)$. We also expect that
the periodic cyclic homology $HP_{\bullet}\left(A_{\bbgamma}\right)$
is equipped with a variation of \nc-Hodge structures of exponential
type over $\mathbb{C}(( t ))$ so that the periodic cyclic homologies
of the categories $\left( A_{\bbgamma + \bc_{i}1_{A}}-\text{{\sf
mod}}\right)$ are the regular pieces of this variation after we base
change to $\mathbb{C}\left(\left( t^{1/N} \right)\right)$. Finally,
the Puiseux series $\{ \bc_{i} \}$ should be the eigenvalues of the
operator of multiplication by $\left[\bbgamma + \gamma_{A} \right] \in
HH^{\bullet}\left(A_{\bbgamma}\right)\widehat{\otimes}_{\mathbb{C}[[t]]}
\mathbb{C}((t))$.

\section{Examples and relation to mirror symmetry} 
 \label{sec:examples}

In this section we discuss examples of \nc-Hodge structures arising
from smooth and compact Calabi-Yau geometries and we study how
these structures are affected by mirror symmetry. Specifically we
look at a generalization of Homological Mirror Symmetry which relates
categories of boundary topological field theories (or $D$-branes)
associated with the following two types of geometric backgrounds:
\begin{description}
\item[{\bfseries $A$-model backgrounds:}] Pairs $(X,\omega)$, where $X$
  is a compact $C^{\infty}$-manifold, and $\omega$ is a symplectic
  form on $X$ satisfying a convergence property (see below).
\item[{\bfseries $B$-model backgrounds:}] Pairs $\bw : Y \to
  \text{disc} \subset \mathbb{C}$, where $Y$ is a complex manifold
  with trivial canonical class, and $\bw$ is a proper holomorphic
  map. 
\end{description}
We will explain how each such background (both in the $A$ and the $B$
model) gives rise to the geometric and Hodge theoretic data described
in Section~\ref{sssec:localizations}. Namely we get:
\begin{itemize}
\item A finite collection $\{Z_{i}^{A/B}\}$ of smooth compact \nc-spaces. In
  fact $\{Z_{i}^{A/B}\}$ will  be (see
  Section~\ref{sssec:CYnc} for the definition) odd/even Calabi-Yau
  \nc-spaces of dimension $\left.\left(\frac{\dim_{\mathbb{R}} X}{2}
  \; \op{mod} \; 2\right)\right\slash \left( \dim_{\mathbb{C}} Y \;
  \op{mod} \; 2\right)$.
\item Complex numbers $\bc_{i}^{A/B}$ and Betti gluing data $\left\{
  T_{ij}^{A/B} \right\}$ for the regular
  \nc-Hodge structures on the periodic cyclic homology of $Z_{i}^{A/B}$. 
\end{itemize}
In particular the data
$\left(HC^{-}_{\bullet}\left(Z_{i}^{A}\right),
\left\{\bc_{i}^{A}\right\}, \left\{ T_{ij}^{A}\right\} \right)$ and 
$\left(HC^{-}_{\bullet}\left(Z_{i}^{B}\right),
\left\{\bc_{i}^{B}\right\}, \left\{ T_{ij}^{B}\right\} \right)$ each
glue into a \nc-Hodge structure of exponential type.  The {\bfseries
  generalized 
Homological Mirror Symmetry Conjecture} now asserts that if two
$A$/$B$-model 
backgrounds $(X,\omega)$/$(Y,\bw)$ are mirror to each other, then the
associated \nc-geometry and \nc-Hodge structure packages are
isomorphic:
\[
\left(Z_{i}^{A},
\left\{\bc_{i}^{A}\right\}, \left\{ T_{ij}^{A}\right\} \right) \cong
\left(Z_{i}^{B},
\left\{\bc_{i}^{B}\right\}, \left\{ T_{ij}^{B}\right\} \right).
\]

\subsection[{\bfseries $A$-model Hodge structures}]{\bfseries 
$A$-model Hodge structures: symplectic manifolds}
\label{subsec:symplectic} 

Suppose $(X,\omega)$ is a compact symplectic manifold of dimension
$\dim_{\mathbb{R}} X = 2d$. 
 In the case when $X$ is a Calabi-Yau variety (in particular $c_1(X)=0$) 
 one has a family of superconformal field theories attached to $X$ in the large volume limit
 (i.e. after the rescaling $\omega \to \omega/\hbar$ where $0<\hbar\ll 1$), and the 
 $A$-twist gives a topological quantum field theory (see \cite{mirrorbook}).
  In mathematical terms it means that we have Gromov-Witten invariants
  and a $\mathbb{Z}$-graded Fukaya category associated to $(X,\omega/\hbar)$.
   One the other side, Gromov-Witten invariants can be defined for an arbitrary
  compact symplectic manifold, not necessarily the one with $c_1(X)=0$.
   Our goal in this section is to describe what is an
 analog of the Fukaya category for general  $(X,\omega)$.
 
  Namely, it is expected that for $(X,\omega)$ of large
volume the Fukaya category of $(X,\omega)$ is a weak 
 $\mathbb{Z}/2$-graded
 $A_{\infty}$-category which will satisfy the generalized
smoothness and compactness properties conjectured in
Section~\ref{sssec:singular.defos}. Briefly this should work as
follows.  Following Fukaya-Oh-Ohta-Ono
\cite{fooo} consider a finite collection $\mathfrak{L} = \{ L_{i} \}$
of transversal oriented spin Lagrangian submanifolds in $X$ and form a
``degenerate'' version $\Fuk_{\mathfrak{L}}$ of Fukaya's category
which only involves the $L_{i}$. More precisely we take $\text{{\sf
Ob}}\left(\Fuk_{\mathfrak{L}}\right) = \{ L_{i} \}$, and define
\[
\op{Hom}_{\Fuk_{\mathfrak{L}}}\left(L_{i},L_{j}\right) = \begin{cases} 
\mathbb{C}^{L_{i}\cap L_{j}}, & i \neq j,  \\
A^{\bullet}(L_{i},\mathbb{C}), & i = j.
\end{cases}
\]
Here $\mathbb{C}^{L_{i}\cap L_{j}}$ is taken with the ordinary algebra
  structure but is put in degree equal to the Maslov grading
  $\op{mod} \; 2$, and $A^{\bullet}(L_{i},\mathbb{C})$ is the dg
  algebra of $C^{\infty}$ differential forms on $L_{i}$. 
  
  We consider a $1$-parameter family of symplectic manifolds
\begin{equation} \label{eq:lvfamily}
\left(X,\frac{\omega}{\hbar}\right), \quad \hbar \in \mathbb{R}_{>0},
\quad \hbar \to 0.
\end{equation} It will be convenient to introduce a new parameter $q :=
\exp(-1/\hbar)$ (note that $q \to 0$ when $\hbar \to 0$).  Denote by 
$\mathbb{C}_{q}$   the usual Novikov ring:
\[
\mathbb{C}_{q} := \left\{ \sum_{i=0}^{\infty} a_{i} q^{E_{i}}
\left| \begin{minipage}[c]{2.4in} formal series where $a_{i} \in
  \mathbb{C}$ and $E_{i} \in \mathbb{R}$ with $\lim_{i \to \infty}
  E_{i} = + \infty$.
\end{minipage}
\right. \right\}\]
In the case $[\omega] \in H^{2}(X,\mathbb{Z})$ one can replace
 the Novikov ring $\mathbb{C}_{q}$  by more familiar algebra $\mathbb{C}((q))$
  of Laurent series. The three-point genus zero Gromov-Witten
  invariants of the symplectic 
family \eqref{eq:lvfamily} give rise (see
e.g. \cite{kontsevich-manin},\cite{li-tian},\cite{siebert},
\cite{cox-katz},\cite{fukaya-ono}) to a $\mathbb{C}_{q}$-valued
(small) quantum deformation of the cup product on $H^{\bullet}(X,\mathbb{C})$:
\[
\bstar_{q} : H^{\bullet}(X,\mathbb{C})^{\otimes 2} \to
H^{\bullet}(X,\mathbb{C})\otimes \mathbb{C}_{q}
\] 
Conjecturally the series for the quantum product is absolutely convergent
 for sufficiently small $q$.

 What is constructed in \cite{fooo} is a solution $\bbgamma$ 
of the Maurer-Cartan
 equation in the cohomological Hochschild complex of  $\Fuk_{\mathfrak{L}}$ 
 with coefficients in the series in $\mathbb{C}_{q}$ with strictly
 positive exponents (equal to the 
 areas of non-trivial pseudo-holomorphic discs).
 The meaning of the quantum product is the cup-product in the Hochschild 
 cohomology of the deformed weak category.

The d$\mathbb{Z}/2$g category $\Fuk_{\mathfrak{L}}$ over  $\mathbb{C}_{q}$ 
is compact but not
smooth. 
If the collection $\mathfrak{L}$ is
chosen  to be big enough, i.e. if it generates
the full Fukaya category, then $\Fuk_{\mathfrak{L}}$ is the
large volume limit of $\Fuk(X,\omega)$, i.e. the limit  in which all
disc instantons for $\omega$ are supressed.

Now the formalism of Section~\ref{sssec:singular.defos} should
associate with $\Fuk_{\mathfrak{L}} = (A-\text{{\sf mod}})$ and
$\bbgamma$ a finite collection 
$\{ \bc_{i} \}$ of formal series in positive powers of $q$ and a collection
$\{\Fuk_{i}\}$ of non-trivial smooth and compact modifications of the
Fukaya category whose Hochschild homologies are the regular singularity
constitutents of the Hochschild homology of the $q$-family of Fukaya
categories near the large volume limit. In this geometric context, we
expect that the $\{ \bc_{i} \}$ are the eigenvalues of the quantum
multiplication operator $c_{1}(T_{X})*_{q}(\bullet)$ acting on
$H^{\bullet}(X,\mathbb{C})\otimes \mathbb{C}[[u]]$. Some evidence for
this comes from the observation that when $c_{1}(T_{X})$ vanishes in
$H^{2}(X,\mathbb{Z})$, then  the Fukaya category is
$\mathbb{Z}$-graded thus is  a fixed point of
the renormalization group.  There is also a more explicit direct
argument identifying the class $c_{1}(T_{X})$ with the infnitesimal generator
of the renormalization group, but we
will not discuss it here.

The formalism of Section~\ref{sssec:singular.defos} now predicts that
the periodic cyclic homology of the Fukaya category, which additively
should be the same as the de Rham cohomology of $X$, should carry a
natural \nc-Hodge structure satisfying the degeneration conjecture
from Section~\ref{sssec:degeneration}.  This expectation is supported
by ample evidence coming from mirror symmetry for Calabi-Yau complete
intersections. Here we present further evidence by describing a
natural \nc-Hodge structure on the de Rham cohomology of a symplectic
manifold and by showing that as $\omega$ approaches the large volume
limit this structure fits in a natural variation of \nc-Hodge
structures.

Using the quantum product $\bstar_{q}$ we will attach to $(X,\omega)$
a variation $((\mycal{H},\nabla),\mycal{E}_{B},\iso)$ of \nc-Hodge
structures over a small disc $\{ q \in \mathbb{C} | \, |q| < r\}$ in
the $q$-plane. First we describe the \nc-Hodge filtration
$(\mycal{H},\nabla)$ and its variation in the $q$-direction:
\begin{itemize}
\item $\mycal{H} := H^{\bullet}(X,\mathbb{C})\otimes \mathbb{C}\{
  u,q\}$
and 
\[
\begin{split}
\mycal{H}^{0} & := \left(\bigoplus_{k = d \mod 2}
H^{k}(X,\mathbb{C})\right)\otimes \mathbb{C}\{ u,q\}\\
\mycal{H}^{1} & := \left(\bigoplus_{k = d + 1 \mod 2}
H^{k}(X,\mathbb{C})\right)\otimes \mathbb{C}\{ u,q\}
\end{split}
\]
\item $\nabla$ is a meromorphic connection on $\mycal{H}$ with poles
  along the coordinate axes $u = 0$ and $q=0$, given by
\[
\begin{split}
\nabla_{\frac{\partial}{\partial u}} & := \frac{\partial}{\partial u} +
u^{-2}\left(\kappa_{X}\bstar_{q} \bullet\right) + u^{-1}\sGr, \\[0.5pc]
\nabla_{\frac{\partial}{\partial q}} & := \frac{\partial}{\partial q}
- q^{-1}u^{-1}\left([\omega]\bstar_{q} \bullet\right),
\end{split}
\] 
where:
\begin{description} 
\item[$\kappa_{X} \in H^{2}(X,\mathbb{Z})$] 
 denotes the first Chern class of the cotangent
  bundle of $X$ computed w.r.t. any $\omega$-compatible almost complex
  structure, and  
\item[$\sGr : \mycal{H} \to \mycal{H}$] is the grading operator defined to be 
$\sGr_{|H^{k}(X,\mathbb{C})} := \frac{k-d}{2}
\op{id}_{H^{k}(X,\mathbb{C})}$. 
\end{description}
\end{itemize} 

\

\noindent
The data $(\mycal{H},\nabla)$ define a $q$-variation of (the de Rham
part of) \nc-Hodge structures. Defining the $\mathbb{Q}$-structure is
much more delicate. To gain some insight into the shape of the
rational local system $\mycal{E}_{B}$ one can look at the monodromy in
the $q$ direction of the algebraic bundle with connection
\[
(H,\nabla)_{|(\mathbb{A}^{1}-\{ 0\})\times \{ q \in \mathbb{C} | \ |q|
  < R\}}, \qquad (H,\nabla) =  
\mathfrak{B}_{\text{along } u}((\mycal{H},\nabla)).
\] 
In some cases the fact that $\mycal{E}_{B}$ should be preserved by
$\nabla$ and the Stokes filtration is rational with respect to
$\mycal{E}_{B}$ is enough to determine $\mycal{E}_{B}$ completely:

\begin{prop} \label{prop:Pn} Let $X = \mathbb{C}\mathbb{P}^{n-1}$ and
  let $\omega$ be the Fubini-Studi form. Let $(H,\nabla)$ be the
  holomorphic bundle with meromorphic connection on
  $(\mathbb{A}^{1}-\{ 0\})\times \{ q \in \mathbb{C} | \ |q| < R\}$
  defined above. Let $\psi \in H$ be a holomorphic section which is
  covariantly constant with respect to $\nabla$. Then 
\begin{itemize}
\item[(a)] For every $u\neq 0$, $\psi\neq 0$ the limit (in a sector of
  the $q$-plane)
\[
\psi_{\op{cl}}(u) = \lim_{q\to 0} \left( \exp \left( -
\frac{\log(q)}{u}\left( [\omega]\wedge (\bullet)\right)\right) \right)\psi
\]
exists. Furthermore, $\psi_{\op{cl}}$ satisfies the
differential equation
\[
\left(\frac{d}{du} + u^{-2}\kappa_{X}\wedge  +
u^{-1}\sGr\right)\psi_{\op{cl}} = 0.
\]
\item[(b)] 
The vector
\[
\psi_{\op{const}}(u) :=
\exp(\log(u)\sGr)\exp\left(\frac{\log(u)}{u}
\kappa_{X}\wedge(\bullet)\right)\psi_{\op{cl}} 
\in H^{\bullet}(X,\mathbb{C})
\]
is independent of $u$. 
\item[(c)] Define the rational structure $\mycal{E}_{B} \subset
  H^{\nabla}$ 
as the subsheaf of all covariantly constant sections $\psi$ for which
  the vector $\psi_{\op{const}} \in
  H^{\bullet}(X,\mathbb{C})$ belongs to the image of the map 
\[
\xymatrix@1@C+1.5pc{
H^{\bullet}(X,\mathbb{Q}) \ar[r]^-{\mathfrak{d}} &
H^{\bullet}(X,\mathbb{C}) \ar[r]^-{\widehat{\Gamma}(X)\wedge
  (\bullet)} & 
H^{\bullet}(X,\mathbb{C}),
}
\]
where $\mathfrak{d} \in GL(H^{\bullet}(X,\mathbb{C}))$ is the operator
of multiplication by $(2\pi i)^{k/2}$ on $H^{k}(X,\mathbb{C})$, and 
$\widehat{\Gamma}(X)$ is a new characteristic class of $X$ defined as 
\[
\widehat{\Gamma}(X) := \exp\left( \boldsymbol{C}ch_{1}(T_{X}) +
\sum_{n\geq 2} \frac{\zeta(n)}{n}ch_{n}(T_{X})\right),
\]
where 
\[
\boldsymbol{C} = \lim_{n\to \infty}\left( 1 + \frac{1}{2} + \cdots +
\frac{1}{n} - \ln(n)\right)
\] 
is Euler's constant, and $\zeta(s)$ is Riemann's zeta function.
\end{itemize}
Then the inclusion $\mycal{E}_{B} \subset
  H^{\nabla}$ is compatible with Stokes data, i.e. the rational
  structure $\mycal{E}_{B}$ satisfies
  {\sf\bfseries ($\mathbb{Q}$-structure axiom)}$^{\bexp}$. 
\end{prop}

\

\noindent
The calculation presented below was known already to B.Dubrovin
\cite[Section 4.2.1]{dubrovin-berlin}, where he also obtained a
   Taylor expansions of a power of a Gamma function in quantum
   cohomology, although he did not
identify it with a characteristic class.

\

\noindent
{\bf Proof of Proposition~\ref{prop:Pn}.} In the standard basis $\{1, h,
h^{2}, \ldots, h^{n-1} \}$ of $H^{\bullet}(\mathbb{P}^{n-1},\mathbb{C})$ 
the connection $\nabla$ on $H$ is given by 

\[
\begin{split}
\nabla_{\!\!\frac{\partial}{\partial u}} & = \frac{\partial}{\partial
  u} + u^{-2}
\begin{pmatrix}
0 & & & nq \\
n & 0  & & \\
  & \ddots & \ddots & \\
  & & n & 0
\end{pmatrix} + u^{-1}
\begin{pmatrix}
\frac{1-n}{2} & & & 0 \\
& \ddots & & \\
& & \ddots & \\
0 & & & \frac{n-1}{2}
\end{pmatrix}  \\
& \\
\nabla_{\!\!\frac{\partial}{\partial q}} & = \frac{\partial}{\partial q} -
q^{-1}u^{-1} \begin{pmatrix}
0 & & & q \\
1 & 0  & & \\
  & \ddots & \ddots & \\
  & & 1 & 0
\end{pmatrix},
\end{split}
\]
If $\psi = \sum_{i = 1}^{n} \psi_{i} h^{i-1}$ is a local section of
$H$, a straightforward check  shows that the condition on $\psi$ to be
 $\nabla$-horizontal is solved by the following ansatz:
 \[
\begin{split}
\psi_{n} & = u^{\frac{1-n}{2}} \int_{\Gamma_{u,q}}\displaylimits
 \exp(\mathcal{F})
\prod_{i=1}^{n-1} \frac{dz_{i}}{z_{i}}  \\[0.5pc]
\psi_{n-1} & = \left( uq \frac{\partial}{\partial q} \right) \psi_{n}
\\[0.5pc] 
\psi_{n-2} & = \left( uq \frac{\partial}{\partial q} \right)^{2}
\psi_{n} \\[0.5pc]
& \cdots \\[0.5pc]
\psi_{1} & = \left( uq \frac{\partial}{\partial q} \right)^{n-1}
\psi_{n}.
\end{split}
\]
Here $\mathcal{F}$ is the function on
$(\mathbb{C}^{\times})^{n-1}$ with coordinates $z_{1}, \ldots,
z_{n-1}$  depending on parameters $u,q\ne 0$ and given by
\[
\mathcal{F}(z_{1},z_{2},\dots,z_{n-1};u,q) := u^{-1}\left(z_{1} + z_{2}
+ \ldots z_{n-1} + \frac{q}{z_{1}z_{2}\ldots z_{n-1}} \right).
\]
The integral is taken over some fixed $(n-1)$-dimensional
semi-algebraic non-compact cycle $\Gamma_{u,q}$ in
$(\mathbb{C}^{\times})^{n-1}$ (depending on the parameters $u$,$q$)
which is going to infinity in directions where
$\op{Re}({\mathcal{F}})\to -\infty$.

More generally, the domain of integration $\Gamma_{u,q}$ used for
defining $\psi_{n}$ can be taken to be a $(n-1)$-dimensional rapid
decay homology chain in $(\mathbb{C}^{\times})^{n-1}$. The rapid decay
homology cycles on smooth complex algebraic varieties are the
natural domains of integration for periods of cohomology classes of
irregular connections. The rapid decay homology was introduced and
studied by Hien \cite{hien1,hien2}, following previous works of
Sabbah \cite{sabbah-2d} and Bloch-Esnault \cite{bloch-esnault-rapid}. In
particular by a recent work of Mochizuki
\cite{mochizuki-2d,mochizuki-wild} and Hien \cite{hien2}it follows
that (after a birational base change) taking periods induces a perfect
pairing between the de Rham cohomology of an irregular connection and
the rapid decay homology. This powerful general theory is not really
needed in our case where the manifold is the affine algebraic torus
 $(\mathbb{C}^{\times})^{n-1}$ but it does provide a useful
perspective. 

Explicitly the non-compact cycles that we will use to generate
horizontal sections of $(H,\nabla)$ will be the $(n-1)$-dimensional
relative cycles for a pair $(\blX,Z)$ constructed as follows. Start with
a smooth projective compactification $\mathfrak{X}$ of
$(\mathbb{C}^{\times})^{n-1}$ with a normal crossing boundary divisor
$D$ which is adapted to $\mathcal{F}$ in the sense that if $u$ and $q$
are nonzero, the divisors of zeroes and poles of $\mathcal{F}$ in
$\mathfrak{X}$ do not intersect with each other, and locally at points
of $D$ the function $\mathcal{F}$ can be written as a product of an
invertible holomorphic function and a monomial in the local
coordinates. Let $\blX$ be the real oriented blow-up of $\mathfrak{X}$
along the divisor $D$. Now consider the real boundary $\partial \blX$ of
$\blX$, i.e. the union of all the boundary divisors of the real oriented
blow-up. The boundary $\partial \blX$ contains a natural open real
semi-algebraic subset $Z \subset \partial \blX$ consisting of all points
$b \in \partial\blX$, such that $|\mathcal{F}(z;u,q)| \to \infty$ when
$z \to b$, and for points $z \in t(\mathbb{C}^{\times})^{n-1}$ near
$b$ the argument of $\mathcal{F}(z;u,q)$ lies strictly in the left
half-plane of $\mathbb{C}$. Note that the real blow-up $X$ has the
same homotopy type as $\mathfrak{X} - D =
(\mathbb{C}^{\times})^{(n-1)}$ and so relative cycles on $(\blX,Z)$ can
be thought of as non-compact cycles on
$(\mathbb{C}^{\times})^{(n-1)}$. Moreover since $Z$ is defined by our
condition on the argument of $\mathcal{F}$, it follows that relative
cycles with boundaries in $Z$ give rise to well defined integrals of 
$\exp(\mathcal{F})\prod z_{i}^{-1}dz_{i}$.

Next observe the integrals over relative cycles with integral coefficients,
i.e. elements in $H_{n-1}(\blX,Z;\mathbb{Z})$  give rise to a
covariantly constant integral lattice in the bundle $(H,\nabla)$.
Furthermore the
Deligne-Malgrange-Stokes filtration is integral with respect to this
lattice. Indeed if we fix a real number
$\lambda$, then whenever $\op{Re}\left(\mathcal{F}\right) <
\lambda\cdot |u|^{-1}$, it follows that $|\exp(\mathcal{F})| <
\exp(\lambda\cdot |u|^{-1})$ when $u \to 0$. Hence the steps of the
Deligne-Malgrange-Stokes filtration  of $(H,\nabla)$ 
are easy to describe in this language:  they correspond to
periods of $\exp(\mathcal{F})\prod z_{i}^{-1}dz_{i}$ on relative
cycles on $(\blX,Z)$ whose boundary is contained in half-planes of the
form $\op{Re}(\mathcal{F}) < \op{const}$.  The periods over cycles
with integral coefficients and the same boundary property then give a
full integral lattice in each such step.

Now to finish the proof of the proposition we have just to calculate
the limiting lattice (which is independent of $u$ and $q$) consisting
of vectors $\psi_{\op{const}} \in H^{\bullet}(X,\mathbb{C})$ defined
in terms of $\psi$ by the formula in part (b) of the statement of the
proposition.
 
 For a general $\nabla$-horizontal local section $\psi= \sum_{i =
 1}^{n} \psi_{i} h^{i-1}$ in a sector at $0$ in the $q$-plane (for given
 $u\ne 0$) one has an asymptotic expansion of $\psi$ at $q\to 0$
 given by:
 \begin{equation} \label{eq:psinexpand}
\psi_{n} = \sum_{i =0}^{n-1} a_{i}(u)(\log q)^{i} + \bO(q (\log q)^n)
+ \ldots,
\end{equation}
Then we have that the ``classical limit'' (at $q\to 0$ where the quantum
multiplication becomes classical) is given by
\[
\psi_{\op{cl}} (u)= \begin{pmatrix}  (n-1)!u^{n-1}a_{n-1}(u) \\
(n-2)!u^{n-2} a_{n-2}(u) \\ \vdots \\ 0! u^{0} a_{0}(u) \end{pmatrix}.
\]
Now we restrict to the case where all variables are real,
 $u<0$, $q>0$ and the contour of integration being
the positive octant $\{(z_1,\dots,z_n)\in
 \mathbb{C}^n|\,z_i>0\,\,\forall i\}$. 

Function $\psi_n=\psi_n(u,q)$ decays exponentially fast at $q\to +\infty$ for a
given $u<0$, hence one can
extract its asymptotic expansion at $q\to 0$ through the Mellin transform:
\[
\int_{0}^{+\infty}\displaylimits
 \psi_{n} q^{s}\frac{dq}{q} = \sum_{i =0}^{\infty} a_{i}(u)
\frac{i!(-1)^{i}}{s^{i+1}} + \bO(1),\,\,\,s\to 0.
\]

This integral can be calculated explicitly
\[
\begin{split}
\int_{0}^{+\infty}\displaylimits \psi_{n} q^{s}\frac{dq}{q} & =
u^{\frac{1-n}{2}} 
\underbrace{\int_{0}^{+\infty}\displaylimits \cdots
  \int_{0}^{+\infty}\displaylimits}_{n \text{ 
    times }} \frac{dq}{q}\prod_{i=1}^{n-1}
\frac{dz_{i}}{z_{i}} \exp\left(u^{-1}\left(z_{1} + z_{2} + \ldots
z_{n-1} + \frac{q}{z_{1}z_{2}\ldots z_{n-1}}\right)\right) q^{s} \\[0.5pc]
& = u^{\frac{1-n}{2}}  \underbrace{\int_{0}^{+\infty}\displaylimits\cdots
  \int_{0}^{+\infty}\displaylimits}_{n-1  \text{
    times }}\prod_{i=1}^{n-1}
\frac{dz_{i}}{z_{i}} \exp\left( u^{-1}\sum_{i}^{n-1}
z_{i}\right) \\[0.5pc]
& \qquad\qquad \cdot \underbrace{\int_{0}^{+\infty} \exp\left(
\frac{q}{uz_{1}z_{2}\ldots z_{n-1}} \right)
q^{s}\frac{dq}{q}}_{\substack{ || \\ \Gamma(s)(-u z_{1}z_{2}\ldots
  z_{n-1})^{s}}} \\
& = u^{\frac{1-n}{2}} (-u)^{s} \Gamma(s)
\int_{0}^{+\infty} \cdots \int_{0}^{+\infty}
 \prod_{i=1}^{n-1} \left(\frac{dz_{i}}{z_{i}}z_{i}^{s}\exp
  \frac{z_{i}}{u} \right)  \\
& = u^{\frac{1-n}{2}} (-u)^{s} \Gamma(s)  \left( (-u)^{s} \Gamma(s)
  \right)^{n-1} \\
&  = u^{\frac{1-n}{2}} (-u)^{ns} \Gamma(s)^{n}.
\end{split}
\]

The conclusion is that the chosen branch $\psi_{\op{cl}}(u)$ is 
completely defined by the expansion
\[
u^{\frac{1-n}{2}} (-u)^{ns} \Gamma(s)^{n} =
\frac{\psi_{\op{cl},n}(u)}{(-u)^{0}s} +
\frac{\psi_{\op{cl},n-1}(u)}{(-u)^{1}s^{2}} + \cdots +
\frac{\psi_{\op{cl},1}(u)}{(-u)^{n-1}s^{n}} + \bO(1),\,\,\,s\to 0
\]
 Furthermore,
all the other branches can be obtained by acting on the branch we know
by the monodromy transformations (around $q=0$)
\[
\frac{(2\pi \sqrt{-1})^{i}}{u^{i}}
\begin{pmatrix} 0 & & & 0 \\
1 & 0  & & \\
  & \ddots & \ddots & \\
  & & 1 & 0
\end{pmatrix}^{i},
\]
for $i = 0, \ldots, n-1$.

\

\noindent

Section $\psi_{\op{cl}}$ satisfies the
differential equation 
\[
\left(\frac{d}{du} + u^{-2}\kappa_{X}\wedge  +
u^{-1}\sGr\right)\psi_{\op{cl}} = 0.
\]
which is the classical limit (at $q\to 0$) of the equation 
\[\nabla_{\!\!\frac{\partial}{\partial u}}(\psi)=0 \]
One can check that the operator  $\frac{d}{du} + u^{-2}\kappa_{X}\wedge  +
u^{-1}\sGr$ can be written as
\[
\exp\left(-\frac{\log(u)}{u}
\kappa_{X}\wedge(\bullet)\right)
\exp(-\log(u)\sGr)
\circ \frac{d}{du}\circ
\exp(\log(u)\sGr)\exp\left(\frac{\log(u)}{u}
\kappa_{X}\wedge(\bullet)\right)\]
This follows from the commutation relation 
\[ [\kappa_{X}\wedge(\bullet), \sGr] =-\kappa_{X}\wedge(\bullet)\]

Finally, in the above formulas one can replace $\log(u)$ by $\log(-u)$
  (and also $u^{\frac{1-n}{2}}$ by $(-u)^{\frac{1-n}{2}}$ with principal values
   at the domain $u<0$. Having this modification in mind, we conclude that
   the vector
\[
\psi_{\op{const}}=\psi_{\op{const}}(u) :=
\exp(\log(-u)\sGr)\exp\left(\frac{\log(-u)}{u}
\kappa_{X}\wedge(\bullet)\right)\psi_{\op{cl}} 
\in H^{\bullet}(X,\mathbb{C})
\]
is independent of $u$, and in particular it coincides with 
$\psi_{\op{cl}}(-1)$,
 as for $u=-1$ the correction matrices relating $\psi_{\op{const}}(u)$ and
 $\psi_{\op{cl}}(u)$ are identity matrices.  Therefore the vector
 $\psi_{\op{const}}$ is given by Taylor coefficients 
 \[ \psi_{\op{const},1}s^0+\dots +\psi_{\op{const},n} s^{n-1}=
  s^n\Gamma(s)^n +\bO(s^n)= \Gamma(1+s)^n+\bO(s^n)\] We see that
$\psi_{\op{const}}\in H^{\bullet}(X,\mathbb{C})$ (after rescaling by
operator $\mathfrak{d}$ from the Proposition) with the value of the
multiplicative characteristic class associated with the series
$\Gamma(1+s)=1+\bO(s)\in \mathbb{C}[[s]]$ and the tangent bundle
$T_X$, because $[T_X]=n[\mathcal{O}(1)]-[\mathcal{O}]$ for
$X=\mathbb{C}\mathbb{P}^{n}$, and by the classical expansion
\[
\log(\Gamma(1+s))=    \boldsymbol{C} s+\sum_{k\ge 2} \frac{\zeta(k)}{k}
   s^k
\]
The action of the monodromy corresponds (up to torsion)
to the multiplication
by $\kappa_{X} \in H^{\bullet}(X,\mathbb{Z})$. \ \hfill $\Box$

\

\medskip

\noindent
The previous proposition suggests the following general definition:

\begin{defi} \label{defi:Qstructure} The rational structure on
  $(H,\nabla)$ is the local subsystem $\mycal{E}_{B} \subset
  H_{|\mathbb{A}^{1}-\{0\}}$ of
  multivalued $\nabla$-horizontal sections whose values at $1$ belong
  to the image of 
\[
\xymatrix@1@C+1.5pc{
H^{\bullet}(X,\mathbb{Q}) \ar[r]^-{\mathfrak{d}} &
H^{\bullet}(X,\mathbb{C}) \ar[r]^-{\widehat{\Gamma}(TX)\wedge
  (\bullet)} & 
H^{\bullet}(X,\mathbb{C}),
}
\]
where $\mathfrak{d} \in GL(H^{\bullet}(X,\mathbb{C}))$ is the operator
of multiplication by $(2\pi i)^{k/2}$ on $H^{k}(X,\mathbb{C})$, and 
$\widehat{\Gamma}(TX)$ is a new characteristic class of $X$ defined as 
\[
\widehat{\Gamma}(TX) := \prod_{i=1}^{d}\Gamma(1 + \lambda_{i}),
\]
where $\Gamma(s)$ is the classical gamma function and $\lambda_{i}$
are the Chern roots of $T_{X}$ computed in any $\omega$-admissible almost
complex structure.
\end{defi}

\

\begin{rem} Apart from the calculation in Proposition~\ref{prop:Pn}
  there are a few other (loose) motivations for this definition:

\

\noindent
$\bullet$ \ The class $\widehat{\Gamma}$ appears in the context of
deformation quantization in the work of the second author
\cite[Section~4.6]{maxim-dq.motives}.

\noindent
$\bullet$ \ The number $\chi(X)\zeta(3)$ appears in the mirror formula
for the quintic threefold.

\noindent
$\bullet$ \ Golyshev's description \cite{golyshev-rr,golyshev-mirror}
of the \nc-motives associated with the Landau-Ginzburg mirror of a
toric Fano involves similar hypergeometric series.

\noindent
$\bullet$ \ The same class $\widehat{\Gamma}$ was derived and a
definition similar to Definition~\ref{defi:Qstructure} was proposed
in the recent work of Iritani \cite{iritani-integral} for the case of
toric orbifolds by tracing out the mirror image of rational structure
of the mirror Landau-Ginzburg model.
\end{rem}

\

\begin{conn} \label{con:symplectic} The triple
  $(H,\mycal{E}_{B},\iso)$ associated above with a symplectic manifold
  $(X,\omega)$ is a variation of \nc-Hodge structures of exponential
  type. 
\end{conn}

\

\begin{rem} \ {\bfseries (i)} \ In general it is not clear if the 
{\bfseries ($\mathbb{Q}$-structure axiom)$^{\bexp}$} holds in this
case. It does hold trivially in the graded case, i.e. when $X$ is a
Calabi-Yau. 

\

\noindent
{\bfseries (ii)} \ At the moment the ``exponential type'' part of the
conjecture is not supported by any evidence beyond the graded case in
which the \nc-Hodge structure is regular. It is possible that for
non-K\"{a}hler symplectic manifolds the \nc-Hodge structure on the de
Rham cohomology is not of exponential type.
\end{rem}

\

\subsection[{\bfseries $B$-model Hodge structures.}]{\bfseries 
$B$-model Hodge structures: holomorphic Landau-Ginzburg
  models} \label{subsec:B.model}

Suppose we have an
algebraic map  $\bw : Y \to \mathbb{C}$, where $Y$ is a smooth
quasi-projective manifold and $\bw$ has a compact critical locus
$\op{crit}(\bw) \subset Y$. Let $S = \{\bc_{1}, \ldots, \bc_{m} \} \subset
\mathbb{C}$ denote the critical values of $\bw$.

A pair $(Y,\bw)$ like that is called a {\em\bfseries holomorphic
  Landau-Ginzburg model} and often arises (see
  e.g. \cite{hori-vafa,mirrorbook}) as the mirror of a symplectic
  manifold underlying a hypersurface, or a complete intersection in a
  toric variety. Remarkably the pair $(Y,\bw)$ give rise to a natural
  \nc-space \nc$(Y,\bw)$. The category $C_{\text{\nc}(Y,\bw)}$ can be
  described in two equivalent ways (in fact these descriptions are
  valid even if the critical locus of $\bw$ is not compact).  First
  note that it is enough to define $\Perf_{C_{\text{\nc}(Y,\bw)}}$
  since that the category $C_{\text{\nc}(Y,\bw)}$ can be thought of as
  the homotopy colimit completion of
  $\Perf_{C_{\text{\nc}(Y,\bw)}}$. For the latter we have two models:
\begin{description}
\item[$\Perf_{C_{\text{\nc}(Y,\bw)}}$ as a category of matrix
  factorizations:]  This model was proposed originally by
the second author as a 
mathematical description of the category of $D$-branes and was
subsequently studied extensively in the physics and mathematics
literature, see \cite{kapustin.li-algebraic,kapustin.li-minimal} and
\cite{orlov-LG,orlov-LG.equivalences,orlov-LG2}.

A matrix factorization on $(Y,\bw)$ is a pair 
$\left(E = E^{0}\oplus E^{1}, d_{E} \in
  \op{End}(E)^{\op{opp}}\right)$, where 
\begin{itemize}
\item[$E$] is a $\mathbb{Z}/2$-graded
  algebraic vector bundle on $Y$, and 
\item[$d_{E}$] is an odd endomorphism
  satisfying $d_{E}^{2} = \bw\cdot \op{id}_{E}$.  
\end{itemize}
 In the case when $Y$ is affine the $\mathbb{Z}/2$-graded complex
  $\underline{\op{Hom}}((E,d_{E}),(F,d_{F})$ of homomorphisms between
  two matrix factorizations is defined as
  $\underline{\op{Hom}}((E,d_{E}),(F,d_{F}) := (\op{Hom}(E,F),d)$
  where for a $\varphi : E \to F$ we have $d\varphi := \varphi\circ
  d_{E} - d_{F}\circ \varphi$. For general $Y$ the same definition
  works if we replace $\op{Hom}(E,F)$ by some acyclic model, e.g. if
  we use the Dolbeault resolution. The resulting category $\MF(Y,\bw)$ of
  matrix factorizations is a $\mathbb{C}$-linear d($\mathbb{Z}/2$)g
  category. We define $\Perf_{C_{\text{\nc}(Y,\bw)}}$ to be the
  derived category $D^{b}(\MF(Y,\bw))$ of the category of matrix
  factorizations.

To construct $D^{b}(\MF(Y,\bw))$
  one notes that in addition to being a
  d($\mathbb{Z}/2$)g category $\MF(Y,\bw)$ can also be viewed as a
  curved d($\mathbb{Z}/2$)g category with central curvature $\bw$ (see
  e.g. \cite{pol.pol} for the definition) or as a
  $\mathbb{Z}/2$-graded weak $A_{\infty}$-category, i.e. an
  $A_{\infty}$ category with an $m_{0}$-operation given by $\bw$ (see
  e.g. \cite{schwartz-curved,hasegawa} for the definition).  In
  particular we can form 
  the associated homotopy category (in the $A_{\infty}$-sense) which
  by definition will be the derived category of matrix
  factorizations.

Alternatively, one can use the following two step construction
   proposed by Orlov. First we pass to the homotopy category of
   $\MF(Y,\bw)$, i.e. we consider the category whose objects are
   matrix factorizations and whose morphisms are given by the quotient
   of $\underline{\op{Hom}}((E,d_{E}),(F,d_{F}))$ by homotopy
   equivalences. Next (following the standard wisdom) we need to
   quotient $\op{Ho}(\MF(Y,\bw))$ by the subcategory of acyclic
   factorizations. Since the matrix factorizations are not complexes,
   they do not have cohomology and so we can not define acyclicity in
   the usual way. But there is another point of view on acyclicity. If
   we have a short exact sequence of usual complexes, then the total
   complex of this diagram will be an acyclic complex. So we define
   {\em\bfseries acyclic} matrix factorizations as the total matrix
   factorization of an exact sequence of factorizations. With
   this definition we get a thick subcategory in the homotopy category
   $\op{Ho}(\MF(Y,\bw))$ matrix factorizations and then we can pass to
   the Serre quotient of $\op{Ho}(\MF(Y,\bw))$ by this thick
   subcategory. We set $D^{b}(\MF(Y,\bw))$ to be this Serre quotient.

\item[$\Perf_{C_{\text{\nc}(Y,\bw)}}$ as a category of
  singularities:] This model was proposed originally by D. Orlov as
  an alternative to the matrix factorization description which is
  localized near the critical set of $\bw$. Orlov proved the
  equivalence of the two models, various versions of the localization
  theorem, and proved several duality statements relating derived
  categories of singularities to other familiar categories 
  \cite{orlov-LG,orlov-LG.equivalences,orlov-LG2}.

Suppose $Z$ is a quasi-projective complex scheme. The {\em\bfseries
  derived category $D^{b}_{\op{Sing}}(Z)$ of singularities} of $Z$ is
defined as the quotient
\[
D^{b}_{\op{Sing}}(Z) := D^{b}(\text{\sf Coh}(Z))/\Perf_{Z}
\]
of the (dg enhancement of the) bounded
derived category $D^{b}(\text{\sf Coh}(Z))$ of coherent sheaves on $Z$
by the thick subcategory of perfect complexes on $Z$. The syzygy
theorem implies  that $D^{b}_{\op{Sing}}(Z) = 0$ whenever $Z$ is
smooth and so  $D^{b}(\text{\sf Coh}(Z))$ can be thought of as an
invariant of the singularities of $Z$.
\end{description}

\

\medskip

\noindent
If now $\bw : Y \to \mathbb{C}$ is a holomorphic Landau-Ginzburg model
we  write $Y_{c}$ for the fiber $\bw^{-1}(c)$ and set
\[
\Perf_{C_{\text{\nc}(Y,\bw)}} := 
D^{b}_{\op{Sing}}(Y_{0}).
\]
Note that if $0 \in \mathbb{A}^{1}$ is not a critical value of $\bw$,  
 then with this 
 definition we will get $\Perf_{C_{\text{\nc}(Y,\bw)}} = 0$. In order
 to get non-trivial categories we will use the critical values $S =
 \{\bc_{1}, \ldots, \bc_{n}\}$ to shift the potentail $\xymatrix@1{\bw
 \ar@{~>}[r] & \bw-\bc_i}$ and associate with $\text{\nc}(Y,\bw)$
 honest categories $\Perf_{i} := \Perf_{C_{\text{\nc}(Y,\bw
 -\bc_{i})}} = D^{b}_{\op{Sing}}(Y_{\bc_{i}})$.  Conjecturally, these
 categories are smooth and compact.

\

\noindent
Mirror symmetry suggests that the \nc-space $\text{\nc}(Y,\bw)$ gives
rise to the $B$-model geometric and Hodge theoretic data described in
Section~\ref{sssec:localizations}, and in particular that the periodic
cyclic homology of $C_{\text{\nc}(Y,\bw)}$ carries a canonical
\nc-Hodge structure. In fact we have already described the geometric
part of the data, namely the numbers $\{\bc_{i}\}$ and the categories
$\{ \Perf_{i} \}$. These data of course fix the regular type (in the
sense of Theorem~\ref{theo:gluenc}) of the
\nc-Hodge structure but we are still missing the gluing data.  Here we
propose a construction of the Hodge structure on the periodic cyclic
homology of $C_{\text{\nc}(Y,\bw)}$ but similarly to the $A$-model we
have to rely on the actual geometry of $(Y,\bw)$ in order to produce
the gluing data. At present it is not clear if the gluing data can be
reconstructed from the category $C_{\text{\nc}(Y,\bw)}$ or more
generally from its one parameter deformation.

\

\medskip

\noindent
First we discuss the appropriate cohomologies of the Landau-Ginzburg
model. Let  
\[
\begin{split}
\mycal{H}^{\bullet}_{\op{for}} & := H_{DR}^{\bullet}((Y,\bw);
\mathbb{C}) \\
& = \mathbb{H}^{\bullet \text{ mod } 2}_{\text{Zar}}\left(Y,
\left(\Omega^{\bullet}_{Y}[[u]],
u d_{DR} + d\bw\wedge\right)\right)
\end{split}
\]
be the $\mathbb{Z}/2$-graded $\mathbb{C}[[u]]$-module
of algebraic de Rham cohomology of the potential $\bw$.  In the case
when $\text{crit}(\bw)$ is compact, the 
$\mathbb{C}[[u]]$-module $\mycal{H}^{\bullet}_{\op{for}}$ is known to
be free by 
the work of Barannikov and the second author (unpublished), Sabbah
\cite{sabbah-twisted}, or
Ogus-Vologodsky \cite{ogus.vologodsky}.  This implies the
following

\begin{lemma} \label{lemma:an.alg} Assume that $Y$ is
  quasi-projective and the critical locus of $\bw$ is compact. Then  we have: 
\begin{itemize}
\item[{\em\bfseries (i)}] The fiber of
 $\mycal{H}^{\bullet}_{\op{for}}$ at $u=0$ is 
 the algebraic Dolbeault cohomology 
\[
\mathbb{H}^{\bullet}_{\text{Zar}}\left(Y,
\left(\Omega^{\bullet}_{Y},
d\bw\wedge 
\right)\right)\cong
\mathbb{H}^{\bullet}_{\text{an}}\left(Y,
\left(\Omega^{\bullet}_{Y},
d\bw\wedge 
\right)\right)
\]
of the potential $\bw$. 
\item[{\em\bfseries (ii)}] There is a canonical isomorphism 
\[
\mathbb{H}^{\bullet}_{\text{Zar}}\left(Y,
\left(\Omega^{\bullet}_{Y}[[u]],
u d_{DR} + d\bw\wedge 
\right)\right) \cong \mathbb{H}^{\bullet}_{\text{an}}\left(Y,
\left(\Omega^{\bullet}_{Y}[[u]],
u d_{DR} + d\bw\wedge 
\right)\right)
\]
\item[{\em\bfseries (iii)}] If the map $\bw$ is proper then 
$\mycal{H}^{\bullet}_{\op{for}}$ is the 
formal germ at $u = 0$ of an algebraic vector bundle on the affine line
 \[\mycal{H}^{\bullet}_{\op{alg}}:= 
 \mathbb{H}^{\bullet \text{ mod } 2}_{\text{Zar}}\left(Y,
\left(\Omega^{\bullet}_{Y}[u],
u d_{DR} + d\bw\wedge\right)\right)\]
\end{itemize}
\end{lemma}
{\bfseries Proof.} The cohomology sheaves of the complex
$\left(\Omega^{\bullet}_{Y},d\bw\wedge \right)$ are supported
on the critical locus of $\bw$ and so, by our compactness assumption,
must be  coherent sheaves on $Y$ both in the
analytic and in the Zariski topology. The hypercohomology spectral
sequence then implies that the hypercohomology of the complex
$\left(\Omega^{\bullet}_{Y},d\bw\wedge \right)$ is finite
dimensional and the spectral sequence associated with the filtration
induced by multiplication by $u$ implies that
$\mathbb{H}^{\bullet}_{\text{Zar/an}}\left(Y,
\left(\Omega^{\bullet}_{Y}[[u]], u d_{DR} + d\bw\wedge 
\right)\right)$ is a finite rank
$\mathbb{C}[[u]]$-module. Furthermore, the same spectral sequence
implies that 
\[
\dim_{\mathbb{C}((u))}
\mathbb{H}^{\bullet}_{\text{Zar/an}}\left(Y,
\left(\Omega^{\bullet}_{Y}((u)), u d_{DR} + d\bw\wedge 
\right)\right) \leq 
\dim_{\mathbb{C}}\mathbb{H}^{\bullet}_{\text{Zar/an}}\left(Y,
\left(\Omega^{\bullet}_{Y}, d\bw\wedge  \right)\right).
\]  
The
freeness statement of Barannikov and the second author (see
e.g. \cite{sabbah-twisted})  now
gives that these two dimensions are equal and so 
$\mathbb{H}^{\bullet}_{\text{Zar}}\left(Y,
\left(\Omega^{\bullet}_{Y}[[u]], u d_{DR} + d\bw\wedge 
\right)\right)$ is a free finite rank module over
$\mathbb{C}[[u]]$. This proves part {\bfseries (i)} of the lemma.

For
part {\bfseries (ii)}  we only need to notice that the two spaces in
question are computed by spectral
sequences associated with the filtrations by the powers of $u$ and
that these spectral sequences have
$E_{2}$-sheets  whose entries are  finite sums of copies of
$\mathbb{H}^{\bullet}_{\text{Zar}}\left(Y, 
\left(\Omega^{\bullet}_{Y}, d\bw\wedge  \right)\right)$
and  $\mathbb{H}^{\bullet}_{\text{an}}\left(Y, 
\left(\Omega^{\bullet}_{Y}, d\bw\wedge \right)\right)$
respectively. Each of these can in turn be computed from the
hypercohomology spectral sequence for the complex
$\left(\Omega^{\bullet}_{Y}, d\bw\wedge  \right)$ of (Zariski or
analytic) coherent sheaves. But the cohomology sheaves of this complex
are supported on the zero locus of $d\bw$ which by assumption is
projective. Hence by GAGA the Zariski and analytic cohomologies of
this complex are naturally isomorphic. This gives isomorphisms of the
hypercohomology and filtration spectral sequences in the Zariski and
the analytic setup respectively and so the two types of
hypercohomologies are isomorphic. 

Finally, part {\bfseries (iii)} was also proven by 
 Barannikov and the second author, and by Sabbah
\cite{sabbah-twisted}. \ \hfill $\Box$ 

\

\medskip

\begin{rem} The isomorphism in part {\bfseries (ii)} of
  the previous lemma is not convergent for $u \to 0$ in
  general. Indeed if $u \neq 0$ is a complex number, then the complex
  vector space \linebreak $\mathbb{H}^{\bullet}_{\text{an}}\left(Y,
  \left(\Omega^{\bullet}_{Y}, u d_{DR} + d\bw\wedge \right)\right)$ is
  the same as the usual de Rham cohomology
  $H^{\bullet}_{DR}(Y,\mathbb{C})$ of $Y$. Indeed, for such a fixed $u
  \neq 0$, the complex $\left(\Omega^{\bullet}_{Y}, ud_{DR} +
  d\bw\wedge \right)) \cong \left(\Omega^{\bullet}_{Y}, d_{DR} +
  u^{-1}d\bw\wedge \right))$ is the holomorphic de Rham complex of the
  local system $(\mathcal{O}_{Y}, d_{DR} + u^{-1}d\bw)$. But the
  multiplication by $\exp(-u^{-1}\bw)$ is an analytic automorphism of
  the line bundle $\mathcal{O}_{Y}$ which gauge transforms the
  connection $d_{DR} + u^{-1}d\bw$ into the trivial connection
  $d_{DR}$. Hence $\exp(-u^{-1}\bw)$ identifies
  $\left(\Omega^{\bullet}_{Y}, u d_{DR} + d\bw\wedge \right)$ with the
  holomorphic de Rham complex $(\Omega^{\bullet}_{Y}, d_{DR})$ and
  $\mathbb{H}^{\bullet}_{\text{an}}\left(Y,
  \left(\Omega^{\bullet}_{Y}, u d_{DR} + d\bw\wedge \right)\right)$
  with $H^{\bullet}_{DR}(Y,\mathbb{C})$.  On the other hand, the space
  \linebreak $\mathbb{H}^{\bullet}_{\text{Zar}}\left(Y,
  \left(\Omega^{\bullet}_{Y}, u d_{DR} + d\bw\wedge \right)\right)$
  depends on the potential in an essential way. For instance, if $\bw
  : Y \to \mathbb{A}^{1}$ is a Lefschetz fibration, then the complex
  $(\Omega^{\bullet}_{Y},d\bw\wedge)$ is just the Koszul \linebreak
  complex associated with the regular section $d\bw \in
  \Omega^{1}_{Y}$. In particular the space \linebreak
  $\mathbb{H}^{\bullet}_{\text{Zar}}(Y,(\Omega^{\bullet}_{Y}, d_{DR} +
  d\bw\wedge)) \cong
  \mathbb{H}^{\bullet}_{\text{Zar}}(Y,(\Omega^{\bullet}_{Y},
  d\bw\wedge))$ has dimension equal to the number of critical points
  of $\bw$. More generally $
  \mathbb{H}^{\bullet}_{\text{Zar}}(Y,(\Omega^{\bullet}_{Y}, d_{DR} +
  d\bw\wedge))$ can be identified (see e.g. \cite{misha}) with the
  cohomology of the perverse sheaf of vanishing cycles of $\bw$.
\end{rem}

\

\medskip

\begin{rem} Under our assumptions, the algebraic de Rham and Dolbeault
  cohomologies $H^{\bullet}_{DR}((Y,\bw);\mathbb{C})$ and
  $H^{\bullet}_{Dol}((Y,\bw);\mathbb{C})$ of the potential $\bw$ can
  be identified respectively with  the periodic cyclic and Hochschild
  homologies $HP_{\bullet}(C_{\text{\nc}(Y,\bw)})$ and
  $HH_{\bullet}(C_{\text{\nc}(Y,\bw)})$ of the \nc-space
  $C_{\text{\nc}(Y,\bw)}$  (more precisely, of the collection
   of categories $\Perf_{i}$ labeled by numbers
  $\{\bc_{i}\}$). This can be done, e.g. by choosing  strong
  generators $\mathcal{E}_i$ of $\Perf_i$,
  and then identifying $HP_{\bullet}(C_{\text{\nc}(Y,\bw)})$ and
  $HH_{\bullet}(C_{\text{\nc}(Y,\bw)})$ with the periodic cyclic and
  Hochschild homologies of the curved d($\mathbb{Z}/2$)g algebra,
  which consists of the  d($\mathbb{Z}/2$)g algebra
  $R\op{Hom}(\mathcal{E},\mathcal{E})$ and a central curvature given
  by $\bw$. A detailed proof of the comparison theorem giving the
  identifications $H^{\bullet}_{DR}((Y,\bw);\mathbb{C})  \cong
HP_{\bullet}(C_{\text{\nc}(Y,\bw)})$ and
$H^{\bullet}_{Dol}((Y,\bw);\mathbb{C})  \cong  
HH_{\bullet}(C_{\text{\nc}(Y,\bw)})$
can be found in the recent work of Junwu Tu \cite{junwu}.
\end{rem}

\

\medskip

\

\noindent
We will construct a \nc-Hodge structure on
$H^{\bullet}_{DR}((Y,\bw);\mathbb{C})$ by using the dual description
of \nc-Hodge structures given in Theorem~\ref{theo:gluenc}.  
Here we will assume that we choose an open subset (in the analytic topology) 
$Y'\subset Y$ such that 
\begin{itemize}
\item $\op{crit}(\bw) \subset Y'$,
\item $\bw(Y')$ is an open disc in $\mathbb{C}$,
\item the closure $\overline{Y}'$ of $Y'$ is a manifold with corners,
\item the restriction of $\bw$ to the part of the boundary 
of $\overline{Y}'$ lying over $\bw(Y')$  is a smooth fibration.
\end{itemize}

In the case when $\bw$ is already proper one can choose $Y'$ to be the 
pre-image under $\bw$ of an open disc in $\mathbb{C}$ containing
 all the critical values $\bc_i$.

\begin{figure}[!ht]
\begin{center}
\psfrag{s}[c][c][1][0]{{$\bc_{i}$}}
\psfrag{o}[c][c][1][0]{{$\bc_{0}$}}
\psfrag{Ds}[c][c][1][0]{{$\Gamma_{i}$}}
\psfrag{as}[c][c][1][0]{{$\gamma_{i}$}}
\epsfig{file=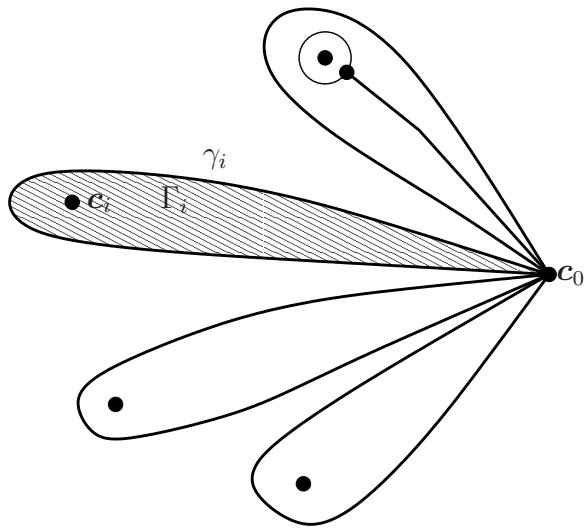,width=3in} 
\end{center}
\caption{A system of thickened loops for $S\subset \mathbb{C}$.}
\label{fig:thick.loops} 
\end{figure}

Label the
critical values of $\bw$: $S = \{\bc_{1}, \ldots, \bc_{n}\}$, and let
$\bc_{0} \in \bw(Y') - S$. Choose a system of paths
$\{a_{i}\}_{i=1}^{n}$ and discs $\bD_{i}$ as in the proof of
Theorem~\ref{theo:gluenc}. Choose $\bc_{0}$-based loops $\gamma_{1},
\ldots, \gamma_{n}$, so that $\gamma_{i}$ goes once around $c_{i}$ in
the counterclockwise direction, all $\gamma_{i}$ intersect only at
$\bc_{0}$, and each $\gamma_{i}$ encloses the path $a_{i}$ and the disc
$\bD_{i}$ (see Figure~\ref{fig:thick.loops}). Let $\Gamma_{i}$ denote
the closed region in $\mathbb{C}$ enclosed by $\gamma_{i}$. Adjusting
if necessary the choice of the $\gamma_{i}$ we can ensure also the
each $\Gamma_{i}$ is convex. From now on we will always assume that
this is the case.

\

\noindent
For $i = 1, \ldots, n$ set $Y_{i} := \bw^{-1}(\Gamma_{i})\cap Y'$ and
consider the $\mathbb{Q}$-vector spaces of
relative cohomology 
\[
U_{i} := H^{\bullet}(Y_{i},Y_{\bc_{0}};\mathbb{Q}),
\]
and
\[
\begin{split}
U & := \oplus_{i =1}^{n} U_{i} \\
& = H^{\bullet}(\bw^{-1}\left(\cup_{i =1}^{n}
\Gamma_{i}\right),Y_{\bc_{0}};\mathbb{Q}) \\
& =  H^{\bullet}(Y,Y_{\bc_{0}};\mathbb{Q}).
\end{split}
\]
Let $T_{i} : U \to U$  be the monodromy along $\gamma_{i}$. By
definition $T_{i}$ satisfies 
\[
\left(T_{i} - 1\right)_{|\oplus_{j \neq i} U_{j}} = 0
\]
and so we get operators $T_{ji} : U_{i} \to U_{j}$, such that
$T_{i|U_{i}} = \sum_{j=1}^{n} T_{ji}$. By construction the operator
$T_{ii}$ is the monodromy along $\gamma_{i}$ of the local system on
$\Gamma_{i}$ of local relative cohomology, i.e. the local system of
$\mathbb{Q}$-vector spaces whose fiber at $\bc \in \Gamma_{i}$ is
$H^{\bullet}(Y_{i},Y_{\bc};\mathbb{Q})$. Hence $T_{ii}$ is an
isomorphism, and so the data $\left( S, \left\{ U_{i}
\right\}_{i=1}^{n}, \left\{T_{ij}\right\}\right)$ are \nc-Betti data
of type {\bf \nc B(iii)}.

\begin{rem} \label{rem:relative}
(a) \ By Lemma~\ref{lemma:B(ii)=B(iii)} the data
$\left( S, \left\{ U_{i} \right\}_{i=1}^{n},
\left\{T_{ij}\right\}\right)$ are the same thing as a constructible
sheaf $\mycal{F}$ of $\mathbb{Q}$-vector spaces on $\mathbb{C}$,
satisfying $R\Gamma(\mathbb{C},\mycal{F}) = 0$. 
The sheaf $\mycal{F}$ can be described directly in terms of
the geometry of $(Y,\bw)$: for a $\bc \in \mathbb{C}$ the stalk
$\mycal{F}_{\bc}$ of $\mycal{F}$ at $\bc$ is the relative cohomology
$H^{\bullet}(Y,Y_{\bc};\mathbb{Q})$. 

\

\noindent
(b) \ The geometric construction of $\mycal{F}$ makes sense for
every cohomology theory $\text{\sf K}$. Indeed for every such
$\text{\sf K}$ we can form a constructible sheaf of abelian groups
$^{\text{\sf K}}\hspace{-0.3pc}\mycal{F}$ whose stalk at $\bc \in
\mathbb{C}$ is $\text{\sf 
  K}(Y,Y_{\bc})$ and which again satisfies
$R\Gamma(\mathbb{C},^{\text{\sf K}}\hspace{-0.3pc}\mycal{F}) 
= 0$. The vanishing of cohomology here is not obvious but can be
proven as follows. Given a disk $\bD \subset \bw(Y')\subset\mathbb{C}$
s.t. $\partial \bD \cap S = \varnothing$, and given any point $\bc \in
\partial \bD$ consider the abelian group $A(\bD,\bc) :=\text{\sf
  K}(\bw^{-1}(\bD),Y_{\bc})$. The collection of abelian groups
$A(\bD,\bc)$ satisfies:
\begin{itemize}
\item $A(\bD,\bc)$ are 
locally constant under small perturbations of $(\bD,\bc)$,
and 
\item for every decomposition $(\bD,\bc) = (\bD_{1},\bc)\cup
(\bD_{2},\bc)$ of $\bD$ obtained by cutting $D$ along a chord starting
at $c$, we have $A(\bD,\bc) = A(\bD_{1},\bc)\oplus A(\bD_{2},\bc)$.
\end{itemize}
This immediately gives
us an equivalent description of $^{\text{\sf
    K}}\hspace{-0.3pc}\mycal{F}$ via data of type {\bf \nc B(iii)},
which in turn yields the vanishing of cohomology of $^{\text{\sf
    K}}\hspace{-0.3pc}\mycal{F}$.
\end{rem}

\

\medskip

\noindent
Next, in order to complete the data {\bf \nc B(iii)} to a full-fledged
\nc-Hodge structure of exponential type, we need to construct:
\begin{itemize}
\item  a
collection $\{(\mycal{R}_{i},\nabla_{i})\}_{i=1}^{m}$ of holomorphic
bundles $\mycal{R}_{i}$ over $\mathbb{C}\{u\}$ equipped with meromorphic
connections $\nabla_{i}$ with at most second order pole and regular
singularities, and 
\item  for each
$i = 1, \ldots, m$, an isomorphism $\mathfrak{f}_{i}$
between the local system on $\bS^{1}$ induced from
$(\mycal{R}_{i},\nabla_{i})$ and the local system on $\bS^{1}$
corresponding to the vector space $U_{i}\otimes \mathbb{C}$ and the
monodromy operator $T_{ii}$.
\end{itemize}

\

\noindent
As explained above the local system on the circle corresponding to
the vector space $U_{i}\otimes \mathbb{C}$ and the monodromy operator
$T_{ii}$ can be described geometrically as the sheaf of complex vector
spaces on the loop $\gamma_{i}$, whose stalk at $\bc \in \gamma_{i}$ is 
$H^{\bullet}((Y_{i},Y_{\bc});\mathbb{C})$. We will exploit this
geometric picture to produce $(\mycal{R}_{i},\nabla_{i})$ and the
isomorphism $\mathfrak{f}_{i}$.  The most convenient way to define the
$\nabla_{i}$  is by using a Betti-to-de Rham cohomology
isomorphism given by oscillating integrals.

Fix $i \in \{1, \ldots, m\}$ and let 
$Z := Y_{i}$, $\bDelta :=
\Gamma_{i}-\bc_{i} \subset \mathbb{C}$, $\bbf := \bw - \bc_{i}$. By
construction  we have:
\begin{description}
\item[$Z$]
is a $C^{\infty}$-manifold with boundary which is the closure of an
open (in the classical topology) subset in the quasi projective complex
manifold $Y$.
\item[$\bDelta \subset \mathbb{C}$] is a closed disc containing zero.
\item[$\bbf : Z \to \bDelta$] is an analytic surjective map whose only
  critical value is zero and whose
  critical locus $\op{crit}(\bbf) \subset Z$ is compact.
\end{description}

\

\noindent
Consider now the $\mathbb{Z}/2$-graded $\mathbb{C}[[u]]$-module
$H^{\bullet}_{DR}((Z,\bbf);\mathbb{C})$ of de Rham cohomology of
$(Z,\bbf)$. By lemma \ref{lemma:an.alg}  we know that
$H^{\bullet}_{DR}((Z,\bbf);\mathbb{C})$ is a free
$\mathbb{C}[[u]]$-module which can be computed as the cohomology of
the complex $\left(\mathcal{A}^{\bullet}(Z)[[u]],d_{\op{tot}}\right)$,
where $\mathcal{A}^{\bullet}(Z)[[u]]$ are the global $C^{\infty}$
complex valued differential forms on $Z$, and $d_{\op{tot}} :=
\bar{\partial} + u\partial + d\bbf\wedge$. The
$\mathbb{C}[[u]]$-module $H^{\bullet}_{DR}((Z,\bbf);\mathbb{C})$
carries a natural meromorphic connection $\nabla$ differentiating in
the $u$-direction and having a second order pole at $u=0$. This
connection is induced from a connection $\boldsymbol{\nabla}$ on the
$\mathbb{C}[[u]]$-module $\mathcal{A}^{\bullet}(Z)[[u]]$ which also
has a second order pole and is defined by the formula
\[
\xymatrix@1{
\boldsymbol{\nabla}_{u^{2}\frac{d}{du}} := u^{2}\frac{d}{du} -
\bbf\cdot (\bullet) + u \sGr :  & \hspace{-2pc}
\mathcal{A}^{\bullet}(Z)[[u]] \ar[r] & \mathcal{A}^{\bullet}(Z)[[u]],
}
\]
where 
\[
\sGr_{|\mathcal{A}^{p,q}(Z)[[u]]} := \frac{q-p}{2}\cdot
\op{id}_{\mathcal{A}^{p,q}(Z)[[u]]}
\] 
is the grading operator coming from \nc-geometry
(compare with  \ref{sssec-usual-hodge}).

\

\noindent
With this definition we have

\begin{lemma} \label{lemma:LGconnection} The operator
  $\boldsymbol{\nabla}_{u^{2}\frac{d}{du}}$ satisfies:
\begin{itemize}
\item[{\em\bfseries (a)}] $\left[ \boldsymbol{\nabla}_{u^{2}\frac{d}{du}},
  d_{\op{tot}}\right] = \frac{u}{2}\cdot d_{\op{tot}} $. 
\item[{\em\bfseries (b)}] $\boldsymbol{\nabla}_{u^{2}\frac{d}{du}}$ preserves
  $\ker(d_{\op{tot}})$ and $\op{im}(d_{\op{tot}})$ and so induces a
  meromorphic connection $\nabla$ with a second order pole on the
  $\mathbb{C}[[u]]$-module $H^{\bullet}_{DR}((Z,\bbf);\mathbb{C})$. 
\end{itemize}
\end{lemma}
{\bfseries Proof.}  We compute 
\[
\begin{split}
\left[ \boldsymbol{\nabla}_{u^{2}\frac{d}{du}},
  d_{\op{tot}}\right] & 
= \left[u^{2}\frac{d}{du} - \bbf
  +u\sGr, \bar{\partial}  + u\partial + d\bbf\wedge\right]
\\[0.5pc]
& = \left[ u^{2}\frac{d}{du}, u\partial \right] - \left[\bbf,
    u\partial\right] + \left[u\sGr,  \bar{\partial}+u\partial +d\bbf\wedge\right]
  \\[0.5pc]
& = u^{2}\partial + ud\bbf\wedge + \frac{u\bar{\partial}}{2}-\frac{u d\bbf\wedge}{2}
-\frac{u^2\partial}{2} \\[0.5pc]
& = \frac{u}{2}\cdot d_{\op{tot}}.
\end{split}
\]

Part {\bfseries (b)} follows immediately from {\bfseries (a)}
 \ \hfill $\Box$

\

\medskip
\noindent
Suppose now that $\alpha = \alpha_{0} + \alpha_{1}u + \alpha_{2}u^{2}
+ \cdots \in \mathcal{A}^{\bullet}(Z)[[u]]$, $\alpha_{i} = \sum
\alpha_{i}^{p,q}$, $\alpha_{i}^{p,q} \in
\mathcal{A}^{p,q}(Z)$ is a $d_{\op{tot}}$-cocycle. Then the
differential $d + u^{-1}d\bbf\wedge = \bar{\partial} + \partial +
u^{-1}d\bbf\wedge =u^{-1/2} u^{\sGr} d_{\op{tot}} u^{-\sGr}$ 
will kill the element 
\[
u^{\sGr}\alpha := \sum_{\substack{i \geq 0 \\ 0 \leq p,q \leq
    \dim Z}} \alpha_{i}^{p,q} u^{i+\frac{q-p}{2}} 
    \in \mathcal{A}^{\bullet}(Z)((u^{1/2})).
\]
Therefore the expression
  $e^{\frac{\bbf}{u}}u^{\sGr}\alpha$ satisfies formally
\[
d\left(e^{\frac{\bbf}{u}}u^{\sGr}\alpha\right)=0,
\]
i.e. is $d$-closed.
Moreover, the action of the operator 
$\boldsymbol{\nabla}_{u^{2}\frac{d}{du}}$ on $\alpha$
translates to the action of $u^{2}\frac{d}{du}$ on the above 
expression modulo formally
exact forms.

\begin{figure}[!ht]
\begin{center}
\psfrag{s}[c][c][1][0]{{$0$}}
\psfrag{dd}[c][c][1][0]{{$\delta$}}
\psfrag{DD}[c][c][1][0]{{$\op{Sec}(\delta)$}}
\psfrag{Ds}[c][c][1][0]{{$\bDelta$}}
\psfrag{as}[c][c][1][0]{{$\partial \bDelta$}}
\epsfig{file=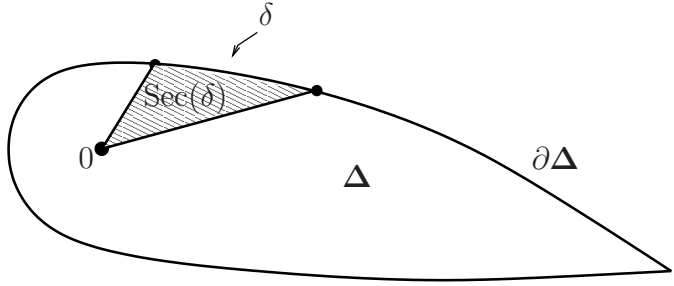,width=3.5in} 
\end{center}
\caption{A sector in $\bDelta$.}\label{fig:sector} 
\end{figure}

Consider now a closed connected arc $\delta \subset \partial \bDelta =
\gamma_{i}$ and let $\op{Sec}(\delta) \subset \bDelta$ be the
corresponding open sector (see Figure~\ref{fig:sector}) with vertex at $0
\in \bDelta$, and boundary made out of the arc $\delta$ and the
segments connecting $0$ with the end points of $\delta$. Note that the
convexity of $\bDelta$ assures that $\op{Sec}(\delta) \subset \bDelta$. Denote by 
$\op{Sec}(\delta)^\vee \subset \mathbb{C}$ the dual angle sector consisting
 of $u\in \mathbb{C}$ such that $\op{Re}(w/u)<0$ for all $w\in \op{Sec}(\delta)$. 

Clearly, for each class in the relative integral homology
$H_{\bullet}(Z,\bbf^{-1}(\delta);\mathbb{Z})$ we can choose a relative
chain $\mathfrak{c}$ representing it, so that $\mathfrak{c}$ satisfies:
\[
\tag{$\boldsymbol{\dagger}$}
\left\|
\text{
\begin{minipage}[c]{3.5in}
\begin{itemize}
\item $\mathfrak{c}$ is piece-wise real
analytic;
\item $\bbf(\op{supp}(\mathfrak{c})) \subset \op{Sec}(\delta)$;
\item $\bbf(\op{supp}(\partial \mathfrak{c}))\subset \delta$.
\end{itemize}
\end{minipage}
}
\right.
\]

\

\noindent
For every such relative chain $\mathfrak{c}$ we now have:

\begin{lemma} \label{lemma:asymptotic}
For every $d_{\op{tot}}$-closed formal power series of forms $\alpha
\in \mathcal{A}_{Z}^{\bullet}(Z)[[u]]$ and every relative chain
$\mathfrak{c} \in C_{\bullet}(Z,\bbf^{-1}(\delta);\mathbb{Z})$
satisfying $(\boldsymbol{\dagger})$ the oscillating integral 
\[
\int_{\mathfrak{c}}\displaylimits e^{\frac{\bbf}{u}}u^{\sGr}\alpha 
\]
is well defined as an asymptotic series in $u^{\mathbb{Q}}(\log
u)^{\mathbb{N}}$ in the sector $\op{Sec}(\delta)^\vee$. 
\end{lemma}
{\bfseries Proof.} Let $N \geq 0$ be a non-negative integer. Clearly the
expression
\[
e^{\bbf/u}u^{\sGr} \left(\sum_{\substack{0 \leq i \leq
    N}}\alpha_{i} u^{i}\right)
\]
is a well defined  analytic
function on $Z\times \op{Sec}(\delta)^\vee$. Using the fact that $(d +
u^{-1}d\bbf\wedge)u^{\sGr}\alpha = 0$ and the
Malgrange-Sibuya theory of 
asymptotic sectorial solutions to analytic differential equations, we
get that 
\begin{equation} \label{eq:truncate}
\int_{\mathfrak{c}}\displaylimits e^{\bbf/u}u^{\sGr}
\left(\sum_{\substack{0 \leq i \leq 
    N}}\alpha_{i} u^{i}\right) \simeq\sum_{j\in \mathbb{Q},k\in \mathbb{N}}
     c_{j,k} u^j (\log(u))^k
\end{equation}
is asymptotic to a series in $u^{\mathbb{Q}}(\log u)^{\mathbb{N}}$ in
which the logarithms enter with bounded powers. Thus the limit of
\eqref{eq:truncate} as $N\to \infty$ is asymptotic to a series in
$u^{\mathbb{Q}}(\log u)^{\mathbb{N}}$ on $\op{Sec}(\delta)^\vee$. \ \hfill
$\Box$

\

\medskip

\noindent
The previous lemma shows that the $\mathbb{C}[[u]]$-module with
connection $(H^{\bullet}_{DR}((Z,\bbf);\mathbb{C}),\nabla)$ is
formally isomorphic to a meromorphic local system of the form
$\be^{\bbf/u}\otimes (\mycal{R}_{i},\nabla_{i})$, where
$\mycal{R}_{i}$ is a free $\mathbb{C}[[u]]$-module, and $\nabla_{i}$
has regular singularities. Furthermore the lemma shows that the
oscillating integrals above identify the local system on $\gamma_{i}$
given by $(\bc \in \gamma_{i}) \mapsto
H^{\bullet}((Y_{i},Y_{\bc}),\mathbb{Q})$ with a rational structure on
$(\mycal{R}_{i}\otimes_{\mathbb{C}[[u]]}\mathbb{C}((u)),\nabla_{i})$. In
particular the data $\{(\mycal{R}_{i},\nabla_{i})\}_{i=1}^{m}$ and
$(S,\{U_{i}\},\{T_{ij}\})$ constitute the regular type and gluing data
(in the sense of Theorem~\ref{theo:gluenc}) 
of a \nc-Hodge structure of exponential type. 

\

\noindent
Usually if one tries to make a Landau-Ginzburg model with proper map
 $\bw$ from non-proper examples above, one gets new parasitic critical
points. Choosing an appropriate domain $Y'\subset Y$ one can define the
 gluing data for the relevant critical points.

\

\subsection{{\bfseries Mirror symmetry examples}}
\label{ssec:mirror.examples}

Finally, in order to give a general idea of the mirror correspondence, 
we briefly discuss three examples of Landau-Ginzburg models 
mirror dual to symplectic manifolds of positive, vanishing, and negative
anti-canonical class respectively.

\

\begin{itemize}
\item For $X=\mathbb{CP}^n$ one of the possible mirror dual
 Landau-Ginzburg models is given by 
 $Y=(\mathbb{C}^\times)^n$ endowed with potential
\[
\bw(z_1,\dots,z_n)=z_1+\dots+z_n+\frac{q}{z_1 \dots z_n}
\]
where $q\in \mathbb{C}^\times$ is a parameter.  In this model the map
$\bw$ is not proper. This can be repaired by compactifying the fibers
of $\bw$ to $(n-1)$-dimensional projective Calabi-Yau varieties. The
compactification is not unique, it depends on combinatorial data, but
the compactified space has the same critical points as $Y$. In
general, for symplectic manifolds $(X,\omega)$ with $\omega$
representing the anticanonical class, one can combine equations for
the connection in $q$ and $u$ directions and get a beautiful variation
of Hodge structures with strong arithmetic properties as predicted by our
considerations in section \ref{subsec:symplectic} (see also 
Golyshev's work \cite{golyshev-rr,golyshev-mirror}).
      
\item For a smooth projective Calabi-Yau variety $X$ one can take for
$Y$ the product $(X^\vee \times \mathbb{A}^{2N},\bw)$ where $X^\vee$
is a Calabi-Yau variety mirror dual to $X$, $N\ge 1$ is arbitrary
integer and $\bw$ is the pullback from $\mathbb{A}^{2N}$ of a
non-degenerate quadratic form. In general, the complex dimension of
the Landau-Ginzburg model is equal to half of the real dimension
of $X$ modulo $2$.

\item For $X$ being a complex curve of genus $g\ge 2$ (considered as a
symplectic manifold), the first author proposed several years ago a
mirror Landau-Ginzburg model $(Y,\bw)$ which is a complex algebraic
$3$-dimensional manifold with non-vanishing algebraic volume element,
such that locally (in the analytic topology) near each point the pair
$(Y,\bw)$ is isomorphic to
\[
\bw:\mathbb{C}^3\to \mathbb{C},\qquad \qquad 
	  (x,y,z)\mapsto xyz
\]
The set of critical point of $\bw$ is the union of $3g-3$ copies of
$\mathbb{CP}^1$ glued along points $0,\infty$ meeting $3$ curves at a
point. The graph obtained by contracting each copy of
$\mathbb{C}^\times$ to an edge is a connected $3$-valent graph with
$g$ loops, representing a maximal degeneration point in the
Deligne-Mumford moduli stack of stable genus $g$ curves.
\end{itemize}

\

\section{Generalized Tian-Todorov theorems and canonical
  coordinates} \label{sec:TT}

In this section we will examine more closely the other direction of
the mirror symmetry correspondence, i.e. the situation in which
symplectic Landau-Ginzburg models appear as mirrors of complex
manifolds with a fixed anti-canonical section. In order to understand
the Hodge theoretic implications of this process we first revisit a
classical concept in the subject: the notion of canonical coordinates.

\subsection[{\bfseries Canonical coordinates for Calabi-Yau
    variations}]
{\bfseries Canonical coordinates for Calabi-Yau variations of
  \nc-Hodge structures} \label{ssec:canonical}

\subsubsection {\bfseries Variations over supermanifolds} \
We begin with a reformulation of the definition of variations of
\nc-Hodge structures  (Definition~\ref{defi:variations}) to allow for
bases that are supermanifolds:

\begin{defi} \label{defi:super.variations} For a complex analytic
  supermanifold $S$, a {\bfseries variation of \nc-Hodge structures
    over $S$} (respectively a {\bfseries variation of \nc-Hodge
    structures over $S$ of exponential type}) is a triple
  $(H,\mycal{E}_{B},\iso)$, where
\begin{itemize}
\item  $H$ is a holomorphic  $\mathbb{Z}/2$-graded vector bundle
  on $\mathbb{A}^{1} 
  \times S$ which is algebraic in the $\mathbb{A}^{1}$-direction;
\item $\mycal{E}_{B}$ is a local system of $\mathbb{Z}/2$-graded
  $\mathbb{Q}$-vector spaces on $(\mathbb{A}^{1}-\{0\})\times S$;
\item $\iso$ is an analytic isomorphism of holomorphic vector bundles 
\[
\iso : \mycal{E}_{B}\otimes \mathcal{O}_{(\mathbb{A}^{1}-\{0\})\times
  S} \stackrel{\cong}{\to} H_{|(\mathbb{A}^{1}-\{0\})\times
  S};
\]
\end{itemize}
so that:
\begin{description}
\item[$\lozenge$]  the induced meromorphic connection $\nabla$ on
$H_{|(\mathbb{A}^{1}-\{0\})\times S}$ satisfies: locally on $S$, for
every section $\xi$ of $T_{S}$, the operators
$\nabla_{u^{2}\frac{\partial}{\partial{u}}}$, $\nabla_{u\xi}$ extend
to operators on $\mathbb{A}^{1}\times S$, and 
\item[$\lozenge$] the triple $(H,\mycal{E}_{B},\iso)$ satisfies 
the {\sf\bfseries ($\mathbb{Q}$-structure axiom)} and 
the {\sf\bfseries (Opposedness axiom)} (respectively $(H,\nabla)$ is
of exponential type and $(H,\mycal{E}_{B},\iso)$ satisfies 
the {\sf\bfseries ($\mathbb{Q}$-structure axiom)$^{\bexp}$} and 
the {\sf\bfseries (Opposedness axiom)$^{\bexp}$}).
\end{description}
\end{defi}

\

\begin{rem} From now on we will suppress the $\mathbb{Q}$-structure
  and the opposedness axioms  since they will not
  play any special role in our analysis. At any given stage of the
  discussion they can be added without any harm or alteration to the
  arguments. 
\end{rem}

\

\subsubsection {\bfseries Calabi-Yau variations} \
Suppose now $(H,\mycal{E}_{B},\iso)$ is a variation  of \nc-Hodge structures
over a supermanifold $S$. For any point $x \in S$ let $H_{0,x}$ denote
the fiber of $H$ at $(0,x) \in \mathbb{A}^{1}\times S$. We get a canonical
map
\[
\mu_{x} : T_{x}S \to \op{End}\left(H_{0,x}\right),
\]
defined as follows: Extend the tangent vector $v \in T_{x}S$ to some
analytic vector field $\xi$ defined in a neighborhood of $x$. Consider
the holomorphic first order differential operator $\nabla_{u\xi} : H
\to H$. By construction this operator has symbol $(u\xi) \otimes
\op{id}_{H}$. In particular, the restriction of $\nabla_{u\xi}$ to the slice
$\{ 0 \}\times S \subset \mathbb{A}^{1}\times S$ will have zero
symbol, and so will be an $\mathcal{O}$-linear endomorphism of $H_{|\{
  0 \}\times S}$. We define $\mu_{x}(v)$ to be the action of this
$\mathcal{O}$-linear map on the fiber $H_{(0,x)}$. It is
straightforward to check that this action is independent of the
extension $\xi$ and depends on $v$ only.

\begin{defi} \label{defi:CYvncHS} Let $S$ be a complex analytic
  supermanifold. We say that a variation $(H,\mycal{E}_{B},\iso)$ of
  \nc-Hodge structures on $S$ is of {\bfseries Calabi-Yau type at a
  point $x \in S$} if there exists an (even or odd) vector $h \in
  H_{(0,x)}$, so that the linear map 
\[
\xymatrix@R-1pc{
T_{x}S \ar[r] & H_{(0,x)} \\
v \ar@{|->}[r] & (\mu_{x}(v))(h)
}
\]
is an isomorphism. Such a vector $h$ will be called a {\bfseries
  generating vector for $H$ at $x$}.
\end{defi}

\

\noindent
It follows from the definition that if $S$ is the base of a variation
of \nc-Hodge structures which is of Calabi-Yau type at a point $x \in
S$, then the tangent space $T_x S$ is a unital commutative associative
algebra acting on $H_{0 ,x}$ via the map $\mu_x$ and such that $H_{0,x}$
is a free module of rank one.  The condition on a variation to have a
Calabi-Yau type (even or odd) is an open condition on $x\in S$.
Variations of \nc-Hodge structures of Calabi-Yau type should arise
naturally on the periodic cyclic homology of smooth and compact
d$(\mathbb{Z}/2)$g categories which are Calabi-Yau in the sense of
\cite{ks-ncgeometry}. The basic geometric example of a Calabi-Yau
variation is an extension of the setup we discussed in
section~\ref{subsec:symplectic}:

\begin{ex} \label{ex:extended.symplectic} Let $(X,\omega)$ be a
  compact symplectic manifold with $\dim_{\mathbb{R}} X =
  2d$. Conjecturally there exists a non-empty open subset $S
  \subset H^{\bullet}(X,\mathbb{C})$ so that the big quantum product
  $*_{x}$ is absolutely convergent for all $x \in S$ (the product is
  given by a formula similar to one on page \pageref{page:qproduct}).
  The manifold 
  $S$ has a natural structure of a supermanifold being an open subset
  in the affine super space $H^{\bullet}(X,\mathbb{C})$.  As in
  section~\ref{subsec:symplectic} we define a variation of \nc-Hodge
  structures $(H,\mycal{E}_{B},\iso)$ on $S$ by taking $H$ to be the
  trivial vector bundle on $\mathbb{A}^{1}\times S$ with fiber
  $H^{\bullet}(X,\mathbb{C})$, and defining the connection
  $\nabla$ on $H$ by the formulas:
\[
\begin{split}
\nabla_{\frac{\partial}{\partial u}} & := \frac{\partial}{\partial u} +
u^{-2}\left(\kappa_{X}\bstar_{x} \bullet\right) + u^{-1}\sGr, \\[0.5pc]
\nabla_{\frac{\partial}{\partial t^{i}}} & := \frac{\partial}{\partial t^{i}}
- q^{-1}u^{-1}\left(t_{i}\bstar_{x} \bullet\right),
\end{split}
\]
where the $(t_{i})$ form a basis on $H^{\bullet}(X,\mathbb{C})$, and
$(t^{i})$  are the dual linear coordinates. 

Clearly, if we restrict $(H,\nabla)$ to $S\cap H^{2}(X,\mathbb{C})$ we
will get back the bundle with connection we defined in
section~\ref{subsec:symplectic}.  We now define the integral lattice
$\mycal{E}_{B}$ and isomorphism $\iso$ on $S$ as the
$\nabla$-horizontal extensions of the integral lattice and isomorphism
we had defined on $S\cap H^{2}(X,\mathbb{C})$. Finally,   in order to
match the framework of \nc-geometry, we should 
 change the parity of the bundle $H$ in the case $d=1 \text{ mod } 2$. 
\end{ex}

\

\subsubsection {\bfseries Decorated Calabi-Yau variations} \
\label{sssec:decorate} The variations
of \nc-Hodge structures of Calabi-Yau type need to be decorated by a
few additional pieces of data before we can extract canonical
coordinates from them. To motivate our choice of such data we first
recall the Deligne-Malgrange classification of logarithmic holomorphic
extension of regular connections.

Let $S$ be a complex analytic supermanifold, let $\bD$ be a one
dimensional complex disc, and let $\mycal{E}$ be a complex local
system on $(\bD-\{\op{pt}\}) \times S$ and let $(\mycal{E},\nabla)$ be
the associated holomorphic bundle $E := \mycal{E}\otimes
\mathcal{O}_{(D-\{\op{pt}\}) \times S}$ on $(D-\{\op{pt}\}) \times S$
with the induced flat connection $\nabla$. Suppose $\widetilde{E}$ is
a holomorphic bundle on $\bD\times S$ which extends $E$ and on which
$\nabla$ has a logarithmic pole. The restriction
$\widetilde{E}_{|\{\op{pt}\}\times S}$ is a holomorphic bundle on $S$
and $\nabla$ induces: a holomorphic connection
${}^{\widetilde{E}}\nabla$ and an $\mathcal{O}_{S}$-linear residue
endomorphism $\op{Res}_{\widetilde{E}}(\nabla)$ on
$\widetilde{E}_{|\{\op{pt}\}\times S}$. Furthermore the integrability
of $\nabla$ on $(\bD-\{\op{pt}\}) \times S$ implies that
${}^{\widetilde{E}}\nabla$ is also integrable and that the
endomorphism $\op{Res}_{\widetilde{E}}(\nabla)$ is covariantly
constant with respect to ${}^{\widetilde{E}}\nabla$
\cite[Section~0.14b]{sabbah-frobenius}.

Recall next that by Deligne's
extension theorem (see e.g. \cite[Chapter~II.5]{deligne} or
\cite[Corollary~II.2.21]{sabbah-frobenius})  meromorphic bundles with
connections with regular singularities always admit functorial
holomorphic extensions across the pole divisor. Deligne's extension
procedure is not unique and depends on the choice of a set-theoretic
section of the quotient  map $\mathbb{C} \to \mathbb{C}/\mathbb{Z}$. 
We
fix $\mathcal{V}$ to be the unique Deligne extension of $E$ for which
$\nabla$ has a logarithmic pole at $\{\op{pt}\}\times S$ and a residue
with eigenvalues whose real parts are in the interval $(-1,0]$. Now
  the classification theorem of Deligne-Malgrange
  \cite[Theorem~III.1.1]{sabbah-frobenius} asserts that there is a
  natural equivalence of categories 
\[
\xymatrix@1{
\left(\text{\begin{minipage}[c]{2in}
Holomorphic extensions of $E$ to $\bD\times S$ for which  $\nabla$ has a
logarithmic singularity along $\{\op{pt}\}\times S$
\end{minipage}}\right)
\ar@{<->}[r] &  \left(\text{\begin{minipage}[c]{2in}
Decreasing filtrations of $\mycal{E}$ by $\mathbb{C}$-local subsystems on 
 $(\bD-\{\op{pt}\}) \times S$
\end{minipage}}\right).
}
\]
The equivalence depends on the chosen Deligne extension and is
explicitly given as follows. Let $t$ be a complex coordinate on $\bD$
which vanishes at $\op{pt} \in \bD$. Consider the restriction
$\mathcal{V}/t\mathcal{V}$ of $\mathcal{V}$ to $\{\op{pt}\}\times
S$. This is a holomorphic bundle on $S$ equipped as above with the
holomorphic connection ${}^{\mathcal{V}}\nabla$ and the covariantly
constant residue endomorphism
$\op{Res}_{\mathcal{V}}(\nabla)$. Suppose now that $\widetilde{E}$ is
another holomorphic bundle 
on $\bD\times S$ which extends $E$ and on which $\nabla$ has a
logarithmic pole. For any $k \in \mathbb{Z}$ we define a subbundle 
$(\mathcal{V}/t\mathcal{V})^{k} \subset \mathcal{V}/t\mathcal{V}$
by setting
\[
\left(\mathcal{V}/t\mathcal{V}\right)^{k} := \frac{\mathcal{V}\cap
t^{k} \widetilde{E}}{t\mathcal{V}\cap
t^{k} \widetilde{E}}
\]
where $\mathcal{V}$ and  $\widetilde{E}$ are viewed as subsheaves in
the meromorphic bundle $E$. 

By construction the sub-bundles
$\left(\mathcal{V}/t\mathcal{V}\right)^{k}$ are preserved both by 
${}^{\mathcal{V}}\nabla$ and by the residue endomorphism
$\op{Res}_{\mathcal{V}}(\nabla)$ and so give rise to
$\nabla$-covariantly constant meromorphic subbundles of $E$, or
equivalently to $\mathbb{C}$-local subsystems of $\mycal{E}$. 

Alternatively we can use a more intrinsic description of holomorphic
extensions of $(\mycal{E},\nabla)$ which is beter adapted to our
examples and in particular to Example~\ref{ex:can.basic}. Namely,
instead of  relying on the Deligne
extension and the induced filtration we can use 
 decreasing filtrations $\mycal{E}_{\leq
\lambda}$ of $\mycal{E}$ labeled by real numbers $\lambda\in
\mathbb{R}$ and such that on the associated graded pieces the
monodromy on $\bD-\{\op{pt}\}$ has eigenvalues in
$\mathbb{R}_{+}\times \exp(2 \pi i \lambda)$.

\

\noindent
We can now introduce the additional data that one
needs for the canonical coordinates 

\begin{defi} \label{defi:decorations}
Let $S$ be a complex supermanifold and let $(H,\mycal{E}_{B},\iso)$ be
a variation of \nc-Hodge structures of Calabi-Yau type on $S$. A
{\bfseries decoration} on $(H,\mycal{E}_{B},\iso)$ is a pair
$(\widetilde{H},\psi)$ where:
\begin{description}
\item[$\widetilde{H}$] is an extension of $H$ to
  $(\mathbb{Z}/2)$-graded vector bundle on $\mathbb{P}^{1}\times S$ so
  that $\nabla$ has a regular singularity at $\{\infty\}\times S$.
\item[$\psi$] is a ${}^{\widetilde{H}}\nabla$-covariantly constant
  section of $\widetilde{H}_{\{\infty\}\times S}$. 
\end{description}
A decoration is called {\bfseries rational} 
 iff the $\mathbb{R}$-filtration on the
local system
$\mycal{E}_{B}\otimes \mathbb{C}$ is compatible with the
 rational structure, and
  if the vector $\psi(x) \in \widetilde{H}_{\{\infty\}\times \{x\}}
=\op{gr} (\mycal{E}_{B}\otimes\mathbb{C})_{x}$ is
  rational, i.e. if $\psi(x) \in \op{gr}
 (\mycal{E}_{B})_{x}$.
\end{defi}

\

\noindent
The previous discussion applied to the local system
$\mycal{E}_{B}\otimes\mathbb{C}$, the disc $\bD = \{ |u| > 1 \} \cup
\{\infty\}$ and the point $\op{pt} = \infty$ shows that the data of a
decoration are equivalent to the data
$\left(\mycal{E}_{B}\otimes\mathbb{C})_{\leq \bullet},\psi\right)$,
where $\left(\mycal{E}_{B}\otimes \mathbb{C}\right)_{\leq \bullet}$ is
a decreasing filtration of $\mycal{E}_{B}\otimes \mathbb{C}$ (labeled
by real numbers) and $\psi$ is a covariantly constant section (along
$S$) of the corresponding logarithmic holomorphic extension of $H$.
We will freely go back and forth between these two points of view.

\

\noindent
Any decorated variation $(H,\mycal{E},\iso;\widetilde{H},\psi)$ of
\nc-Hodge structures of Calabi-Yau type gives rise to a natural open
domain $U \subset S$ defined by
\[
U := \left\{ x \in S \;  \left| \;
\text{\begin{minipage}[c]{3.5in} $\widetilde{H}_{\mathbb{P}^{1}\times
      \{x \}}$ is holomorphically trivial and if \linebreak 
$s \in \Gamma\left(
    \mathbb{P}^{1}\times 
      \{x \} \right)$ is such that $s_{x}(\infty) = \psi(\infty,x)$,
      then  $s_{x}(0)$ is a generating vector for $(H,\mycal{E},\iso)$.
\end{minipage}}\right.\right\}
\]
Furthermore for every $x\in U$ we get a natural map $\op{{\sf can}}_{x} :
T_{x}S \to \widetilde{H}_{\infty,x}$ defined as the composition
\[
\xymatrix@1@C+2.5pc{ T_{x}S \ar[r]^-{\mu_{x}(\bullet)(s_{x}(0))}
\ar@/_1.5pc/[rrr]_-{\op{{\sf can}}_{x}} & H_{0,x}
\ar[r]^-{\op{ev}_{(0,x)}^{-1}} & \Gamma\left(
\widetilde{H}_{|\mathbb{P}^{1}\times \{x\}}\right)
\ar[r]^-{\op{ev}_{(\infty,x)}} & \widetilde{H}_{\infty,x}.  }
\]
where $\op{ev}_{(t,x)} :  \Gamma\left(\mathbb{P}^{1},
  \widetilde{H}_{|\mathbb{P}^{1}\times \{x\}}\right) \to
  \widetilde{H}_{t,x}$ denotes the natural evaluation of sections,
  which is invertible by the triviality assumption on
  $\widetilde{H}_{|\mathbb{P}^{1}\times \{x\}}$. 

\

\noindent
The pullback of the flat connection ${}^{\widetilde{H}}\nabla$ by the
map $\op{{\sf can}}$ induces a flat connection on $TS_{|U}$. The
canonical coordinates on $S$ come from the following easy claim whose proof
we omit

\begin{claim} \label{claim:torsion.free} The flat connection $\op{{\sf
      can}}^{*}\left({}^{\widetilde{H}}\nabla\right)$ on
  $TS_{|U}$ is torsion free and so gives rise to a natural affine
  structure and affine coordinates on $U$. If the decoration is rational then
  the tangent bundle $TS_{|U}$ carries a natural rational structure.
\end{claim}

\

\begin{rem} \ {\bfseries (i)} The canonical coordinates on $U$
  corresponding to a decorated \nc-variation of Hodge structures are
  only affine coordinates and are defined only up to a translation.

\

\noindent
{\bfseries (ii)} For any $u \in \mathbb{A}^{1} - \{0\}$ we can
introduce another affine structure which is a {\bfseries vector
  structure}. In fact, we get an analytic isomorphism between $U$ and a
domain in $H_{u,\bullet} = (\mycal{E}_{B})_{u,\bullet}\otimes
\mathbb{C}$:
\[
x \in U \mapsto \op{ev}_{(u,x)}\op{ev}_{(\infty,x)}^{-1} (\psi(x)) \in
H_{(u,x)}.
\]
One can use this to show that the local Torelli theorem holds for
decorated Calabi-Yau variations of \nc-Hodge structures.
\end{rem}

\

\begin{ex} \label{ex:can.basic}
The setup of Example~\ref{ex:extended.symplectic} gives not only a
variation of \nc-Hodge structures but in fact gives a rationally 
decorated
\nc-Hodge structure of Calabi-Yau type. Indeed by definition the
fibers of $H$ are identified with
$\bPi^{d}H^{\bullet}(X,\mathbb{C})$. The monodromy of the connection
around $\infty \in \mathbb{P}^{1}$ is the operator acting by
$(-1)^{i+d}\exp(\kappa_{X}\wedge (\bullet))$ on
$H^{i}(X,\mathbb{C})$. Consider the monodromy invariant filtration on
$H^{\bullet}(X,\mathbb{C})$ whose step in degree $\frac{d-i}{2}$ is
$H^{\geq i}(X,\mathbb{C})$. Let $\widetilde{H}$ be the corresponding
logarithmic extension of $H$ and let $\psi$ be the section of
$\widetilde{H}$ corresponding to the image of $1 \in
H^{0}(X,\mathbb{C}) \subset H^{\bullet}(X,\mathbb{C})$. The bundle
$\widetilde{H}_{|\{\infty \} \times S}$ is trivialized and
$\nabla_{\frac{\partial}{\partial t^{i}}} = \frac{\partial}{\partial
  t^{i}}$ in this trivialization. This gives the desired decoration
$(\widetilde{H},\psi)$ and the associated canonical coordinates 
are the standard canonical coordinates in Gromov-Witten theory. 
\end{ex}

\

\medskip

\subsubsection {\bfseries Generalized decorations} \ The notion of a
decorated Calabi-Yau variation of \nc-Hodge structures can be
generalized in various ways. For instance, instead of specifying a
covariantly constant filtration on $H$ giving the extension
$\widetilde{H}$ we can start with any holomorphic bundle $H'$ defined
on  $\{u \in  \mathbb{P}^{1}| \, |u| \geq R \}$,
and an identification of $C^{\infty}$-bundles 
\[
p_{1}^{*}\left(H'_{|\{|u| = R\}}\right) \cong
(\mycal{E}_{B}\otimes \mathbb{C})_{\{|u| = R\}\times S},
\]
where $p_{1} : \{|u| = R\}\times S \to \{|u| = R\}$ is the projection
on the first factor. 

Furthermore (locally in $S$) the holomorphic bundle
$p_{1}^{*}H'$  on $\{u \in  \mathbb{P}^{1}| \, |u| \geq R \}\times S$
carries a flat connection defined along $S$ only. We can use the above
identification to glue this
together with $H$ along $\{|u| = R\}\times S$ to get a bundle
$\widetilde{H}$ on $\mathbb{P}^{1}\times S$ equipped with a flat
connection $\nabla_{/S}$ along $S$. This generalizes the first part of
the decoration. For the second part we will take a
$\nabla_{/S}$-covariantly constant section $\psi$ of
$\widetilde{H}_{|\{\infty \}\times S}$. Now the same definition of the set
$U$ and the canonical map $\op{{\sf can}}$ make sense in this
context. The resulting connection on $TS_{|U}$ is again torsion free.

\

\medskip

\subsubsection {\bfseries Formal variations of Calabi-Yau type} \ The notion
of a Calabi-Yau variation extends readily to the formal
context. Suppose $S = \op{Spf}
\mathbb{C}[[x_{1},\ldots,x_{N},\xi_{1},\ldots,\xi_{M}]]$ be a formal
algebraic supermanifold, where $x_{i}$ are even and $\xi_{j}$ are odd 
formal variables. The de Rham part of formal variation of
\nc-Hodge structures on $S$ is a pair $(H,\nabla)$ where $H$ is a
$(\mathbb{Z}/2)$-graded algebraic vector bundle over $\mathbb{D}\times
S$, where $\mathbb{D}$ is the one dimensional formal disc $\mathbb{D}
:= \op{Spf}(\mathbb{C}[[u]])$. Here $\nabla$ is a meromorphic
connection on $H$ such that
$\nabla_{u^{2}\frac{\partial}{\partial{u}}}$,
$\nabla_{u\frac{\partial}{\partial{x_{i}}}}$,
$\nabla_{u\frac{\partial}{\partial{\xi_{j}}}}$ are regular
differential operators on $H$.

We say that such a pair $(H,\nabla)$ has the Calabi-Yau property if we
can find a vector $h \in H_{0,0}$, so that the natural linear map 
$T_{0}S \to H_{0,0}$, $v \mapsto \mu_{0}(v)(h)$ is an isomorphism. 

Finally a decoration of a formal Calabi-Yau de Rham data  $(H,\nabla)$
is a pair $(\boldsymbol{e},h)$, where $\boldsymbol{e}$ is a
trivialization $\boldsymbol{e} : H_{|\mathbb{D}\times \{0\}} \to
H_{0,0}\otimes \mathcal{O}_{\mathbb{D}\times \{0\}}$, and $h \in
H_{0,0}$ is a generating vector for the Calabi-Yau property. 

Again a decorated de Rham data of Calabi-Yau type gives an affine
structure and canonical formal coordinates on $S$.

\subsection{\bfseries Algebraic framework: dg Batalin-Vilkovisky
  algebras} \label{ssec:dgBV} 

In this section we discuss the aspects of algebraic deformation theory
relevant to the study of \nc-Hodge structures. We will work over
$\mathbb{C}$ but all algebraic considerations in this section make
sense over any field of characteristic zero. 

\subsubsection {\bfseries Preliminaries on $L_{\infty}$ algebras} \ 
Our main objects of interest here will be differential $\mathbb{Z}/2$-graded
algebras over $\mathbb{C}$ or more generally $\mathbb{Z}/2$-graded
$L_{\infty}$-algebras over $\mathbb{C}$. We begin with a definition:

\begin{defi} \label{defi:habelian}
A complex differential $\mathbb{Z}/2$-graded Lie algebra $\mathfrak{g}$ 
(or a $\mathbb{Z}/2$-graded $L_{\infty}$-algebra) is called
{\bfseries homotopy abelian}  if it is $L_{\infty}$
quasi-isomorphic to an abelian d$(\mathbb{Z}/2)$g Lie algebra.
\end{defi}

\

\begin{rem} \label{rem:equiv.ha} Homotopy abelian differential
  $\mathbb{Z}/2$-graded Lie algebras can be characterized in a variety
  of ways. In particular we have the following statements that follow
  readily from the definition:
\begin{itemize} 
\item A differential $\mathbb{Z}/2$-graded Lie algebra $\mathfrak{g}$
  is homotopy abelian if and only if all the higher operations
  $m_{n}$ vanish on its $L_{\infty}$ minimal model
  $\mathfrak{g}^{\min} = H^{\bullet}(\mathfrak{g},d_{\mathfrak{g}})$,
  i.e. $m_{n} = 0$ for $n \geq 1$.
\item A differential $\mathbb{Z}/2$-graded Lie algebra $\mathfrak{g}$
  is homotopy abelian if and only if there exist
  d$(\mathbb{Z}/2)$g Lie algebras $\mathfrak{g}_{1}$ and
  $\mathfrak{g}_{2}$, and morphisms of  d$(\mathbb{Z}/2)$g
  Lie algebras:
\[
\xymatrix@R-1pc@C-1pc{
& \mathfrak{g_{1}} \ar[dl]_-{\cong} \ar[dr]^-{\cong} & \\
\mathfrak{g} & & \mathfrak{g}_{2}
}
\]
so that $\mathfrak{g}_{2}$ is an abelian d$(\mathbb{Z}/2)$g Lie algebra,
and the morphisms $\mathfrak{g}_{1} \to \mathfrak{g}$ and
$\mathfrak{g}_{1} \to \mathfrak{g}_{2}$ are quasi-isomorphisms. 
\item A differential $\mathbb{Z}/2$-graded Lie algebra $\mathfrak{g}$
  is homotopy abelian if and only if the Lie algebra cohomology
  algebra $H^{\bullet}(\mathfrak{g},\mathbb{C})$ is free, i.e. is
  isomorphic to the algebra of formal power series on some (possibly
  infinitely many) supervariables. Here the Lie algebra cohomology is
  defined as 
\[
H^{\bullet}(\mathfrak{g},\mathbb{C}) := H^{\bullet}\left(\prod_{n\geq
0}\op{Hom}_{(\mathbb{C}-\text{\sf Vect})}
(\op{Sym}^{n}\bPi\mathfrak{g},\mathbb{C})^{\bullet},d\right)
\]
where $d$ is the cochain Cartan-Eilenberg
differential. 
\end{itemize}
\end{rem}

\

\noindent
After the prioneering work of Deligne and Drinfeld in the 80's, it
is by now a common wisdom (see
e.g. \cite[Chapter~III.9]{manin-frobenius}) that dg Lie algebras give
rise to solutions of moduli problems. In particular a homotopy abelian
d$(\mathbb{Z}/2)$g Lie algebra $\mathfrak{g}$ gives rise to a moduli
space - the formal supermanifold
$\MM_{\left(\mathfrak{g},d_{\mathfrak{g}}\right)} :=
\op{Spf}(H^{\bullet}(\mathfrak{g},\mathbb{C}))$.

\

\medskip

\noindent
The property of being homotopy abelian is preserved by suitably
non-degenerate deformations and various other natural operations:

\begin{prop} \label{prop:ha.properties} \
{\em\bfseries (i)} \ Let $\mathfrak{g}$ be a flat family of
  d$(\mathbb{Z}/2)$g Lie algebras (or $(\mathbb{Z}/2)$-graded
  $L_{\infty}$ algebras) over $\mathbb{C}[[u]]$. That is
  $\mathfrak{g}$ is a flat $(\mathbb{Z}/2)$-graded
  $\mathbb{C}[[u]]$-module, and the Lie bracket and differential on
  $\mathfrak{g}$ are $\mathbb{C}[[u]]$-linear. Assume further that
\begin{itemize}
\item[(A)] $\mathfrak{g}_{\op{gen}} :=
  \mathfrak{g}\otimes_{\mathbb{C}[[u]]} \mathbb{C}((u))$ is
  homotopy abelian over $\mathbb{C}((u))$, and
\item[(B)] $H^{\bullet}\left(\mathfrak{g},d_{\mathfrak{g}}\right)$ is a flat 
$\mathbb{C}[[u]]$-module.
\end{itemize}
Then the special fiber $\mathfrak{g}_{0} :=
\mathfrak{g}\widehat{\otimes}_{\mathbb{C}[[u]]} \mathbb{C}$ is also a
homotopy abelian d$(\mathbb{Z}/2)$g Lie algebra over $\mathbb{C}$.

\

\noindent
{\em\bfseries (ii)} \ If $\mathfrak{g}$ is a homotopy abelian
d$(\mathbb{Z}/2)$g Lie algebra over $\mathbb{C}$, and
$\mathfrak{g}_{1} \to \mathfrak{g}$ is a morphism of
$L_{\infty}$-algebras inducing a monomorphism
$H^{\bullet}\left(\mathfrak{g}_{1},d_{\mathfrak{g}_{1}}\right)
\hookrightarrow H^{\bullet}\left(\mathfrak{g},d_{\mathfrak{g}}\right)$, then 
$\mathfrak{g}_{1}$ is homotopy abelian as well.

\

\noindent
{\em\bfseries (iii)} \ If $\mathfrak{g}$ is a homotopy abelian
d$(\mathbb{Z}/2)$g Lie algebra over $\mathbb{C}$, and $\mathfrak{g}
\to \mathfrak{g}_{2}$ is a morphism of $L_{\infty}$-algebras inducing
an epimorphism $H^{\bullet}\left(\mathfrak{g},d_{\mathfrak{g}}\right)
\twoheadrightarrow
H^{\bullet}\left(\mathfrak{g}_{2},d_{\mathfrak{g}_{2}}\right)$, then
$\mathfrak{g}_{2}$ is homotopy abelian as well.
\end{prop}
{\bf Proof.} The proof is standard so we only mention some of the
highlights of the argument. First note that parts {\bfseries (ii)} and
{\bfseries (iii)} follow immediately by passing to minimal models. For
part {\bfseries (i)} we note first that the assumption {\em (B)}
implies (and is in fact equivalent to) the existence of
$\mathbb{C}[[u]]$-linear quasi-isomorphisms $p_{1}$, $p_{2}$ of
complexes:
\[
\xymatrix@1{
\left(H^{\bullet}\left(\mathfrak{g}_{0},d_{\mathfrak{g}_{0}}\right)[[u]],
0\right) & \hspace{-2pc}
\cong  & \hspace{-2pc}
\left(H^{\bullet}\left(\mathfrak{g},d_{\mathfrak{g}}\right),
0\right) \ar@<1ex>[r]^-{p_{1}}  &
\left( \mathfrak{g}, d_{\mathfrak{g}} \right), \ar@<1ex>[l]^-{p_{2}}
} 
\]
and a $\mathbb{C}[[u]]$-linear homotopy $h$ so that 
\[
\begin{split}
p_{2}\circ p_{1} & = \op{id} \\[-0.5pc]
p_{1}\circ p_{2} & = \op{id} + \left[d_{\mathfrak{g}}, h\right].
\end{split}
\]
Next note that the homological perturbation theory of
\cite{kontsevich.soibelman-tori} carries over verbatim to the
$L_{\infty}$-context and gives explicit expressions
for the higher products $m_{n}$ on
$\left(H^{\bullet}\left(\mathfrak{g}_{0},d_{\mathfrak{g}_{0}}\right)[[u]],
0\right)$ as a polynomial expression in $p_{1}$, $p_{2}$ and $h$. In
particular the operations $m_{n}$ are all $\mathbb{C}[[u]]$-linear and
are given by universal expressions.  But by assumption {\em
(A)} we know that the higher operations are zero after tensoring with
$\otimes_{\mathbb{C}[[u]]} \mathbb{C}((u))$ and so $m_{n} = 0$ as
formal power series in $u$ for all $n \geq 1$. This implies that
$m_{n|u=0} = 0$ for all $n \geq 1$ and so the proposition is proven. \
\hfill $\Box$

\

\subsubsection {\bfseries DG Batalin-Vilkovisky algebras} \ Recall
\cite[Chapter~III.10]{manin-frobenius} the notion of a dg BV algebra:

\begin{defi} \label{defi:dgBV} A differential $\mathbb{Z}/2$-graded
  Batalin-Vilkovisky algebra over $\mathbb{C}$ is the data
  $(A,d,\Delta)$, where $A$ is a $\mathbb{Z}/2$-graded
  suppercommutative associative unital algebra, and $d : A \to A$,
  $\Delta : A \to A$ are odd $\mathbb{C}$-linear maps satisfying:
\begin{itemize}
\item $d(1) = \Delta(1) = 0$,
\item $d$ is a differential operator of order $\leq 1$ on $A$,
\item $\Delta$ is a differential operator of order $\leq 2$ on $A$,
\item $d^{2} = \Delta^{2} = d\Delta + \Delta d = 0$.
\end{itemize}
\end{defi}

\

\noindent
Note that the first two properties in the definition imply that $d$ is
a derivation of $A$. Also   $\mathfrak{g} := \bPi A$ together with
$[a,b] := \Delta(ab) - \Delta(a)b - (-1)^{\deg(a)} a\Delta(b)$ is
a Lie superalgebra with two anti-commuting differentials $d$ and
$\Delta$.

\

\begin{defi} \label{defi:degeneration} We will say that a
  d$(\mathbb{Z}/2)$g Batalin-Vilkovisky algebra $(A,d,\Delta)$ has the
  {\bfseries degeneration property}  if for every $N \geq 1$ we have
  that 
$H^{\bullet}(A[u]/(u^{N}),d+u\Delta)$ is a free
  $\mathbb{C}[u]/(u^{N})$-module. 
\end{defi}

\

\noindent
Equivalently $(A,d,\Delta)$ has the degeneration property iff
$H^{\bullet}(A[[u]],d+u\Delta)$ is a topologically free (flat)
$\mathbb{C}[[u]]$-module. This in turn is equivalent to the existence
of a (non-unique) isomorphism of topological $\mathbb{C}[[u]]$-modules:
\begin{equation} \label{eq:deg.iso}
\xymatrix@1{T : & \hspace{-1.5pc} H^{\bullet}(A[[u]],d+u\Delta)
  \ar[r]^-{\cong} 
  & H^{\bullet}(A,d)[[u]].
}
\end{equation}
In this situation we will always normalize $T$ so that $T_{|u = 0} =
\op{id}_{H^{\bullet}(A,d)}$. 

\
 
\noindent
The degeneration property for dg Batalin-Vilkovisky algebras defined
  above is weaker than the $\partial\bar{\partial}$-lemma used
  Barannikov and the second author in \cite{bk98} and by Manin in
  \cite{manin-frobenius,manin-3}. In particular it has potentially a
  wider scope of applications - a feature that we will exploit
  next. We begin with a general smoothness result which was also proven by
  J.Terilla \cite{terilla}.

\

\begin{theo} \label{theo:dgBV} 
Suppose $(A,d,\Delta)$ is a d$(\mathbb{Z}/2)$g Batalin-Vilkovisky
algebra  which has the degeneration property. Let $\mathfrak{g} :=
\bPi A$ be the associated super Lie algebra with anti-commuting
differentials $d$ and $\Delta$. Then:
\begin{itemize}
\item[{\em\bfseries (1)}] The d$(\mathbb{Z}/2)$g Lie algebra
  $(\mathfrak{g},d)$ is homotopy abelian, i.e. is quasi-isomorphic to
  $H^{\bullet}(\mathfrak{g},d)$ endowed with the trivial bracket and
  the zero differential. In particular the associated moduli space
  $\MM_{(\mathfrak{g},d)}$ is (non-canonically) isomorphic to a formal
  neighborhood of $0$ in the super affine space $\bPi
  H^{\bullet}(\mathfrak{g},d)$.
\item[{\em\bfseries (2)}] Every choice of a normalized
  degeneration isomorphism $T$ as in equation \eqref{eq:deg.iso} gives an
  identification of formal manifolds
\[
\xymatrix@1{
\Phi_{T} : & \hspace{-1.5pc}  \MM_{(\mathfrak{g},d)} \ar[r]^-{\cong} &
\left(\text{\begin{minipage}[c]{1.5in} formal neighborhood of $0$ in $\bPi
  H^{\bullet}(\mathfrak{g},d)$
\end{minipage}}\right)
}
\]
\end{itemize}
\end{theo}
{\bf Proof.} Part {\bfseries (1)} of the theorem follows immediately from

\begin{lemma} \label{lemma:dpha} The d$(\mathbb{Z}/2)$g Lie algebra
$\left(\mathfrak{g}((u)),d+u\Delta\right)$ is homotopy abelian over 
$\mathbb{C}((u))$. 
\end{lemma}
{\bf Proof.} Consider the formal completion at zero $\widehat{A}$ of
the vector superspace underlying $A =
\bPi\mathfrak{g}$ as an algebraic supermanifold, and let as before
$\mathbb{D} = \op{Spf}(\mathbb{C}[[u]])$ be the formal one dimensional
disc. The d$(\mathbb{Z}/2)$g Lie algebra structure on
$\mathfrak{g}[[u]]$ is  encoded in an odd vector field $\boldsymbol{\xi}
\in \Gamma(\widehat{A}\times \mathbb{D},T)$ on the supermanifold
$\widehat{A}\times 
\mathbb{D}$, defined by
\[
\dot{a} := \boldsymbol{\xi}(a) = da + u\Delta a + \frac{1}{2}[a,a].
\]
There is a natural automorphism (i.e. a formal change of coordinates)
$F : \widehat{A}\times \mathbb{D}^{\times} \to \widehat{A}\times
\mathbb{D}^{\times}$ on 
the formal supermanifold $\widehat{A}\times \mathbb{D}^{\times}$ given by
\[
F(a) := u\left(\exp\left(\frac{a}{u}\right) - 1\right) = a +
\frac{1}{u}\frac{1}{2!} a^{2} + \frac{1}{u^{2}}\frac{1}{3!} a^{3} +
\cdots,
\]
and in the new coordinates $b = F(a)$ the vector field $\boldsymbol{\xi}$
is linear:
\[
\begin{split}
\dot{b} & = \dot{a}\cdot
\exp\left(\frac{a}{u}\right) = \left( da
+ u\Delta a + \frac{1}{2}[a,a] \right)\cdot \exp\left(\frac{a}{u}\right)
\\[+0.5pc]
& = u\cdot \left(\frac{da}{u} + u \Delta\left(\frac{a}{u}\right)
+u\frac{1}{2}\left[\frac{a}{u}, \frac{a}{u}\right] \right)\cdot
\exp\left(\frac{a}{u}\right) \\[+0.5pc]
& = u\cdot (d + u\Delta) \exp\left(\frac{a}{u}\right) = (d + u\Delta)b.
\end{split}
\]
So in the $b$-coordinates, the vector field $\boldsymbol{\xi}$ depends
only on the differential $d + u\Delta$ and does not depend on any
higher operations.  Passing to the minimal model we see that
$(\mathfrak{g}((u)),d+u\Delta)$ is homotopy abelian, which proves the
lemma. \ \hfill $\Box$

\

\noindent
The lemma implies that the hypothesis {\em (A)} of
Proposition~\ref{prop:ha.properties} {\bfseries (i)} holds. On the
other hand the hypothesis {\em (B)} holds
by the degeneration assumption. Therefore by
Proposition~\ref{prop:ha.properties} {\bfseries (i)} we conclude that
$(\mathfrak{g},d)$ is homotopy abelian. This proves part {\bfseries
  (1)} of the theorem. 

\medskip

\

\noindent
Next we construct the identification $\Phi_{T}$. Given a formal path
in $\MM_{(\mathfrak{g},d)}$, i.e. a family of
solutions (up to guage equivalence)
\[
\begin{split}
a(\varepsilon) & = a_{1}\varepsilon + a_{2}\varepsilon^{2} +
a_{3}\varepsilon^{3} + \cdots \in \varepsilon A[[\varepsilon]] \\
d(a(\varepsilon & )) + \frac{1}{2}[a(\varepsilon),a(\varepsilon)] = 0
\end{split}
\]
of the Maurer-Cartan equation in $(\mathfrak{g},d)$,  we have to
construct the corresponding formal path through the origin in 
$H^{\bullet}(\mathfrak{g},d)$. 

As a first step choose a lift of  the formal arc $a(\varepsilon)$ to a
formal series in two variables $\tilde{a}(\varepsilon,u) \in \varepsilon
A[[\varepsilon,u]]$ such that
\[
\begin{split}
(d + & u\Delta)\tilde{a} + \frac{1}{2}[\tilde{a},\tilde{a}] = 0, \\
a(\varepsilon &,0) =  a(\varepsilon).
\end{split}
\]
Consider the reparameterization 
\[
\tilde{b} = F(\tilde{a}) = u\left(\exp\left(\frac{\tilde{a}}{u}\right)
- 1\right) 
\in \varepsilon A((u))[[\varepsilon]].
\]
Arguing as before we see that $\tilde{b}$ satisfies 
$(d +  u\Delta)\tilde{b} = 0$. So if we expand
\[
\tilde{b} = \tilde{b}_{1}\varepsilon + \tilde{b}_{2}\varepsilon^{2} +
\cdots, \qquad\qquad \text{ where } \tilde{b}_{n} \in A((u)), \text{
  satisfy }  (d +
u\Delta)\tilde{b}_{n} = 0,
\]
we can define cohomology classes $\left[\tilde{b}_{n}\right] \in
H^{\bullet}(A((u)), d + u\Delta)$. We can now apply the isomorphism 
$T\otimes_{\mathbb{C}[[u]]} \mathbb{C}((u))$ to the series 
\[
\sum_{n \geq 1}\displaylimits \left[\tilde{b}_{n}\right]
\varepsilon^{n} \in \varepsilon
H^{\bullet}(A((u)),d+u\Delta)[[\varepsilon]], 
\]
to obtain an element
\[
T\left( \sum_{n \geq 1}\displaylimits \left[\tilde{b}_{n}\right]
\varepsilon^{n}\right) \in \varepsilon
H^{\bullet}(A,d)((u))[[\varepsilon]].
\]
In fact one has the following lemma whose proof we will skip since it
is a somewhat tedious application of
homological perturbation theory: 

\begin{lemma} \label{lemma:choice.of.a}
There exists a lift $\tilde{a}(\varepsilon,u)$ of $a(\varepsilon)$
such that the associated class $T\left( \sum_{n \geq 1}\displaylimits
\left[\tilde{b}_{n}\right] \right)$ belongs to $\varepsilon
H^{\bullet}(A,d)[[\varepsilon]]\subset \varepsilon
H^{\bullet}(A,d)((u))[[\varepsilon]]$. Any such lift $\widetilde{a}$
produces the same class $T\left( \sum_{n \geq 1}\displaylimits
\left[\tilde{b}_{n}\right] \right)$ and this class depends only on the
gauge equivalence class of the original arc $a$, i.e. on the image
$\underline{a}(\varepsilon)$  of
$a(\varepsilon)$ in $\MM_{(\mathfrak{g},d)}$. 
\end{lemma}

\

\noindent
Now by definition the map $\Phi_{T}$ assigns the class 
$T\left( \sum_{n \geq 1}\displaylimits
\left[\tilde{b}_{n}\right] \right) \subset \varepsilon
H^{\bullet}(A,d)[[\varepsilon]]$ 
to the formal arc $\underline{a}(\varepsilon)$. \ \hfill $\Box$

\

\medskip

\subsubsection {\bfseries Geometric interpretation}
\begin{figure}[!ht]
\begin{center}
\psfrag{MM}[c][c][1][0]{{$\MM$}}
\psfrag{MMg}[c][c][1][0]{{$\MM_{(\mathfrak{g},d)}$}}
\psfrag{u}[c][c][1][0]{{$u$}}
\psfrag{Hgd}[c][c][1][0]{{$H^{\bullet}(\mathfrak{g},d)$}}
\psfrag{0}[c][c][1][0]{{$0$}}
\psfrag{inf}[c][c][1][0]{{$\infty$}}
\epsfig{file=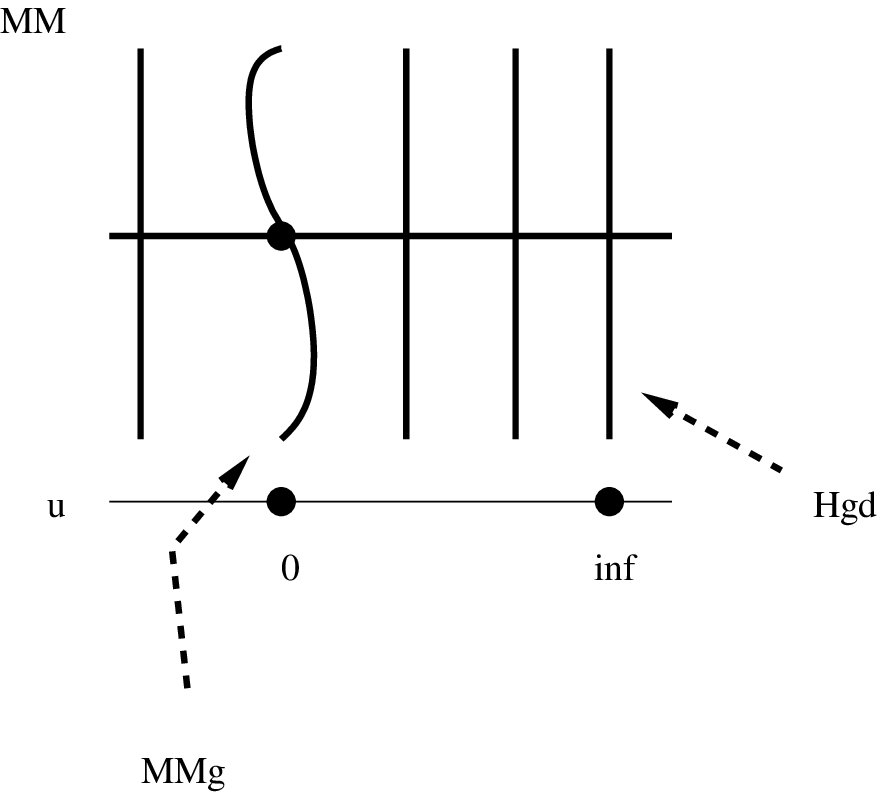,width=3in} 
\end{center}
\caption{The relative moduli $\MM \to \mathbb{P}^{1}$.}\label{fig:affine}  
\end{figure}

\ 
 \label{sss:geometry} The previous discussion can be
  repackaged geometrically as follows. A $(\mathbb{Z}/2)$-graded
  Batalin-Vilkovisky algebra $(A,d,\Delta)$, gives rise to a family
  $\mycal{M} \to \mathbb{D} = \op{Spf}(\mathbb{C}[[u]])$ of formal
  manifolds over the one dimensional formal disc. The family
  $\mycal{M}$ is the total space of the relative moduli space
  $\MM_{(\mathfrak{g}, d +u\Delta)}$ over $\mathbb{C}[[u]]$. If
  $(A,d,\Delta)$ has the degeneration property, then by
  Lemma~\ref{lemma:dpha} we have an affine structure on the generic
  fiber $\mycal{M}^{\op{gen}} := \mycal{M}\otimes_{\mathbb{C}[[u]]}
  \mathbb{C}((u))$ of the family (see Figure~\ref{fig:affine}) given
  by the map $F$.

\

\smallskip

\

\

\smallskip

\noindent
  Furthermore the map $T$ can be viewed as an
  extension of the affine bundle $\mycal{M}^{\op{gen}} \to
  \mathbb{D}^{\times}$ to a trivial bundle on $\mathbb{P}^{1}-\{ 0 \}$
  of formal super affine spaces, where the fiber at $\infty$ is the
  super affine space $H^{\bullet}(\mathfrak{g},d)$. This results into
  a family $\MM \to \mathbb{P}^{1}$ of formal super manifolds, which
  is a trivial vector bundle outside of zero but has a non-linear
  fiber at $0 \in \mathbb{P}^{1}$. Moreover by picking the closed
  point in each fiber we get a section of $\MM \to \mathbb{P}^{1}$,
  which is just the zero section of the vector bundle
  $\MM_{|\mathbb{P}^{1}-\{ 0 \}} \to \mathbb{P}^{1}-\{ 0 \}$. The
  normal bundle to this section in $\MM$ is trivial (hence $\MM$ is
  trival as a non-linaer bundle 
    over $\mathbb{P}^{1}$), and the map
  $\Phi_{T}$ gives a (non-linear) trivialization of $\MM$ over
  $\mathbb{P}^{1}$. This type of geometry was already discussed in
  \cite{chen-kontsevich-schwarz}. 

\

\medskip

\subsubsection {\bfseries Relation to Calabi-Yau variations of \nc-Hodge
  structures} \  Suppose $(A,d,\Delta)$ is a d$\mathbb{Z}/2$g
Batalin-Vilkovisky algebra which has the degeneration property. In
this generality one does not expect to find a natural connection on
$H^{\bullet}(A,d + u\Delta)$ along $u$, i.e. one does not
expect to have a general formal analogue of a \nc-Hodge structure.

However, a natural connection along the $u$-line may exist if we
specify some additional data on $(A,d,\Delta)$. Following the analogy
with the \nc-Hodge structure associated with a symplectic manifold and
the Gromov-Witten invariants,  it is sufficient to specify:
\begin{itemize}
\item an even element $\kappa \in A$, with $d\kappa = 0$, and
\item a grading operator $\sGr : A \to A$,
\end{itemize}
so that if we consider $\Gamma_{-1} := \sGr : A \to A$,
and $\Gamma_{-2} : A \to A$ - the operator of multiplication by
$\kappa$, then we have the commutation relations:
\[
\begin{split}
\left[\Gamma_{-1},\Delta\right] & =  -\frac{1}{2} \Delta\\
\left[\Gamma_{-2},d\right] & = 0 \\
d & = \left[\Gamma_{-1},d\right] + \left[\Gamma_{-2},\Delta\right].
\end{split}
\]
These commutation relations imply the identity
\[
\left[ u\frac{\partial}{\partial u}  + u^{-1}\Gamma_{-2} +
  \Gamma_{-1}, d + u\Delta \right] = \frac{1}{2}(d + u\Delta),
\]
which is consistent with the general formulas from
Section~\ref{sssec:mero.u}.  
In particular, we can define a connection on
$H^{\bullet}(A,d + u\Delta)$ along the $u$-line by setting
\[
\nabla_{\frac{\partial}{\partial u}} :=  \frac{\partial}{\partial u}
+ u^{-2}\Gamma_{-2} + u^{-1}\Gamma_{-1}.
\]
\

\begin{ex} \label{ex:CYdgBV}  
Let $Y$ be a (possibly non-compact) $d$-dimensional 
Calabi-Yau manifold
  with a fixed holomorphic volume form $\bOmega_{Y}$. Let
  $\bw : Y \to \mathbb{C}$ be a proper holomorphic function. This
  geometry gives rise to a natural dg Batalin-Vilkovisky algebra:
\[
\begin{split}
A & := \Gamma_{C^{\infty}}\left(Y,\wedge^{\bullet} T^{1,0}_{Y}\otimes
\wedge^{\bullet}A^{0,1}_{Y}\right), \\
d & := \bar{\partial} + \iota_{d\bw}, \\
\Delta & := \ddiv_{\bOmega_{Y}} = \iota_{\bOmega_{Y}}^{-1}\circ
\partial \circ \iota_{\bOmega_{Y}},
\end{split}
\]
where $\iota_{\bOmega_{Y}} : \wedge^{\bullet} T^{1,0}_{Y}  \to
\wedge^{d-\bullet} \Omega^{1,0}_{Y}$ denotes the contraction with
$\bOmega_{Y}$. 

As discussed in section \ref{subsec:B.model} in this situation we get
a connection along $u$ which conjecturally defines a \nc-Hodge
structure. The connection is defined the above formula with
$\Gamma_{-2} = $ the operator of multiplication by $-\bw$, and
$\Gamma_{-1} = \sGr : A \to A$, the grading operator which is equal to
$\frac{q+p-d}{2}\cdot \op{id}$ on $\Gamma_{C^{\infty}}\left(Y,\wedge^{p}
T^{1,0}_{Y}\otimes \wedge^{q}A^{0,1}_{Y}\right)$.

We will elaborate on this geometric picture in the next section.
\end{ex}

\subsection{\bfseries $B$-model framework: manifolds with anticanonical
  sections} \label{ssec:anticanonical}

\subsubsection {\bfseries The classical Tian-Todorov theorem.} \
\label{sssec:Tian.Todorov} Let $X$ be a
compact K\"{a}hler manifold. By Kodaira-Spencer
theory we know that the deformations of $X$ are controlled by the dg
Lie algebra
\[
\left( \mathfrak{g}^{(1)}, d_{\mathfrak{g}^{(1)}} \right) :=
\left(\Gamma_{C^{\infty}}\left(X, T^{1,0}_{X}\otimes_{\mathcal{C}^{\infty}_{X}}
  A^{0,\bullet}_{X}\right), \bar{\partial} \right).
\]
The classical Tian-Todorov theorem \cite{tian,andrey} can be
formulated as follows:

\

\begin{theo} \label{theo:Tian.Todorov}
 If $X$ is a compact  K\"{a}hler manifold with
  $c_{1}(X) = 0 \in \op{Pic}(X)$, then $\left( \mathfrak{g}^{(1)},
  d_{\mathfrak{g}^{(1)}} \right)$ is homotopy abelian. In
  particular the formal moduli space of $X$ is smooth.
 \end{theo}
{\bfseries Proof.}  Since $c_{1}(X) = 0 \in \op{Pic}(X)$ we can find a
unique up to scale holomorphic volume form $\bOmega_{X}$ on $X$. As in
 example \ref{ex:CYdgBV} the pair $(X,\bOmega_{X})$ gives rise to a 
dg Batalin-Vilkovisky algebra $(A,d,\Delta)$: 
\[
\begin{split}
A & := \Gamma_{C^{\infty}}\left(X,\wedge^{\bullet} T^{1,0}_{X}\otimes
\wedge^{\bullet}A^{0,1}_{X}\right) \\
d & := \bar{\partial} \\
\Delta & := \ddiv_{\bOmega_{X}} = \iota_{\bOmega_{X}}^{-1}\circ
\partial \circ \iota_{\bOmega_{X}}.
\end{split}
\]
Consider the associated dg Lie algebra $\left(\mathfrak{g},
d_{\mathfrak{g}} \right) := \left( \bPi A, d \right)$. We have a
natural inclusion of dg Lie algebras
\[
\xymatrix@R-1pc{
\left( \mathfrak{g}^{(1)}, d_{\mathfrak{g}^{(1)}} \right)
\ar@^{{(}->}[r] \ar@{=}[d] & \left(\mathfrak{g},
d_{\mathfrak{g}} \right) \ar@{=}[d]  \\
\left(\Gamma_{C^{\infty}}\left(X,
T^{1,0}_{X}\otimes_{\mathcal{C}^{\infty}_{X}} 
  A^{0,\bullet}_{X}\right), \bar{\partial} \right) \ar@^{{(}->}[r]
    &  \Gamma_{C^{\infty}}\left(X,\wedge^{\bullet} T^{1,0}_{X}\otimes
\wedge^{\bullet}A^{0,1}_{X}\right)
}
\]
which embeds $\left( \mathfrak{g}^{(1)}, d_{\mathfrak{g}^{(1)}}
\right)$ as a direct summand in $\left(\mathfrak{g}, d_{\mathfrak{g}}
\right)$, and so induces and embedding $H^{\bullet}\left(
\mathfrak{g}^{(1)}, d_{\mathfrak{g}^{(1)}} \right) \subset
H^{\bullet}\left(\mathfrak{g}, d_{\mathfrak{g}} \right)$ in
cohomology. So by Proposition~\ref{prop:ha.properties} it suffices to
check that 
$\left(\mathfrak{g}, d_{\mathfrak{g}} \right)$ is homotopy abelian. 

On the other hand the contraction map $\iota_{\bOmega_{X}}$ gives an
isomorphism of bicomplexes between the dg Batalin-Vilkovisky algebra
$(A,d,\Delta)$ and
the Dolbeault bicomplex
$(A^{\bullet}(X),\bar{\partial},\partial)$. Since $X$ is assumed
compact and K\"{a}hler, the Hodge-to-de Rham spectral sequence
degenerates on $X$ which is equivalent to the equality 
$\dim H^{k}_{dR}(X,\mathbb{C}) = \dim (\oplus_{p+q = k}
H^{p}(X,\Omega_{X}^{q}))$ which implies that the Dolbeault bicomplex
$(A^{\bullet}(X),\bar{\partial},\partial)$ has the degeneration
property. Thus by Theorem~\ref{theo:dgBV} {\bfseries (1)} it follows that 
$\left(\mathfrak{g}, d_{\mathfrak{g}}
\right)$ is homotopy  abelian. The theorem is proven. \ \hfill $\Box$

\

\subsubsection [{\bfseries Moduli of Calabi-Yau
  manifolds.}]{\bfseries Canonical coordinates on the moduli of Calabi-Yau
  manifolds.} \ \label{sssec:CYcan} Let $X$ be a Calabi-Yau manifold, i.e. a
  $d$-dimensional compact K\"{a}hler manifold with $c_{1}(X) =  0$ in
  $\op{Pic}(X)$. Let $(A,d,\Delta)$ be the dg Batalin-Vilkovisky
  algebra defined in section~\ref{sssec:Tian.Todorov}. The contraction
  map $\iota_{\bOmega_{X}}$ identifies the 
  $\mathbb{C}[[u]]$-module $H^{\bullet}\left(\mathfrak{g}[[u]], d +
  u\Delta\right)$ with the Rees module of the \nc-Hodge filtration on
  $H^{\bullet}_{dR}(X,\mathbb{C})$ for which $H^{p,q}(X) \subset
  F^{\frac{p-q}{2}}$. Now choose  one of the following
  equivalent pieces of data:
\begin{itemize}
\item a filtration $G^{\bullet}$ on
  $H^{\bullet}_{dR}(X,\mathbb{C})$ which is opposed to the \nc-Hodge
  filtration,
\item a splitting of the \nc-hodge filtration,
\item an extension of the associated \nc-Hodge structure to a trivial
  bundle on $\mathbb{P}^{1}$ such that the connection has at most a
  first order pole at infinity. 
\end{itemize}
Each such choice gives rise to an affine structure on
$\MM_{\left(\mathfrak{g}, d_{\mathfrak{g}}\right)}$.  This affine
structure is the same as the one described in
section~\ref{sssec:decorate} corresponding to the \nc-Hodge structure
above and the decoration $\psi$ given by the class
$\left[\bOmega_{X}\right]$ in the associated graded
$\op{{\sf gr}}_{G^{\bullet}} H_{dR}^{\bullet}(X,\mathbb{C})$.

In mirror symmetry considerations a choice of this type arises
naturally when $X$ is a Calabi-Yau manifold near a large complex
structure limit point. Concretely, suppose $X = X_{\mathfrak{z}}$ is
member in a
holomorphic family $\{X_{z}\}$ of compact $d$-dimensional Calabi-Yau
manifolds parameterized by $z$ in a polydisc $\prod_{i = 1}^{M} \{
z_{i} \in \mathbb{C} \, | \, 0 < |z_{i}| \ll 1 \}$, and such that:
\begin{itemize}
\item $M = \dim_{\mathbb{C}} H^{1}\left(X_{z},T_{X_{z}}\right)$;
\item for each $i = 1, \dots, M$ the monodromy operator $t_{i} \in
  GL\left(H^{1}\left(X_{\mathfrak{z}},T_{X_{\mathfrak{z}}}\right)\right)$
  assigned to the 
  circle (traced counterclockwise) 
\[
\gamma_{i} = \left\{ z \, \left| \, \text{\begin{minipage}[c]{1in}
$z_{j} = \mathfrak{z}_{j}$, $j \neq i$,
$|z_{i}| = |\mathfrak{z}_{i}|$
\end{minipage}
}
\right.
\right\}
\]
is unipotent of order $d$.
\end{itemize}
In this setup, the filtration $G^{\bullet}$ of
$H^{\bullet}\left(X_{\mathfrak{z}},\mathbb{C}\right)$ invariant under
all unipotent  operators
$\prod_{i =1}^{M} t_{i}^{a_{i}}$, $a_{i} \in \mathbb{Z}_{>0}$ will be
opposed to the Hodge filtration and will thus give us canonical
coordinates on the polydisc.  This affine structure corresponds to a
rational decoration of a Calabi-Yau variation of \nc-Hodge structures.

\

\medskip

\subsubsection {\bfseries Generalizations.} \ \label{ssec:generalize} 
Here we generalize the previous
discussion to the case of varieties with divisors. 

\

\smallskip

\noindent
{\bfseries (i)} \ Let $X$ be a $d$-dimensional smooth projective
variety over $\mathbb{C}$, and let $D \subset X$ be a normal
crossings anti-canonical divisor, i.e $\mathcal{O}_{X}(D) = K_{X}^{-1}
\in \op{Pic}(X)$. Typically such an $X$ will be a Fano or a
quasi-Fano. If $D$ is smooth, then by adjunction $D$ will be a
Calabi-Yau. Specifying such a divisor is equivalent to
specifying a logarithmic volume form on $X$. This is a unique up to
scale $n$-form $\bOmega_{X\log D} \in
\Gamma\left(X,\Omega_{X}^{d}(\log D)\right)$ on $X$ which has a first
order pole along $D$ and does not vanish anywhere on $X - D$.

Let $T_{X,D}$ be the subsheaf of $T_{X}$ of holomorphic vector fields
on $X$ which at the points of $D$ are tangent to $D$. This is a
locally free subsheaf of $T_{X}$ of rank $d$ which controls the
deformation theory of the pair $(X,D)$. The relevant dg
Batalin-Vilkovisky algebra $(A,d,\Delta)$ is an obvious generalization
of the one in the absolute case:
\[
\begin{split}
A & := \Gamma_{C^{\infty}}\left(X,\wedge^{\bullet}
T_{X,D}\otimes_{\mathcal{C}_{X}^{\infty}} 
\wedge^{\bullet}A^{0,1}_{X}\right) \\
d & := \bar{\partial} \\
\Delta & := \ddiv_{\bOmega_{X\log D}} = \iota_{\bOmega_{X\log D}}^{-1}\circ
\partial \circ \iota_{\bOmega_{X\log D}},
\end{split}
\]
where $\iota_{\bOmega_{X\log D}} : \wedge^{\bullet} T_{X,D} \to
\Omega^{d -\bullet}_{X}(\log D)$ is the isomorphism given by
contraction with $\bOmega_{X\log D}$. 

Again the map $\iota_{\bOmega_{X\log D}}$ identifies $(A,d,\Delta)$
with the logarithmic Dolbeault bicomplex
$\left(A^{\bullet,\bullet}(\log
D),\bar{\partial},\partial\right)$. In particular, for all $u \neq 0$
we get an identification of the cohomology of the complex
$(A,d+u\Delta)$ with the cohomology of the total complex of the double
complex $\left(\Omega_{X}^{\bullet,\bullet}(\log
D),\bar{\partial},\partial\right)$, which is equal 
\cite[Section~6.1]{voisin-v2} to the cohomology of the open variety
$X-D$. In other words for all $u \neq 0$ we have an isomorphism
\begin{equation} \label{eq:iso(i)}
H^{\bullet}(A,d+u\Delta) \cong H^{\bullet}_{dR}(X-D,\mathbb{C}).
\end{equation}
Now mixed Hodge theory implies the following

\begin{lemma} \label{lemma:deg.with.D(i)}
The logarithmic  dg Batalin-Vilkovisky algebra
$\left(A,d,\Delta\right)$ has the degeneration property. In
particular the formal moduli of the pair $(X,D)$ is smooth.
\end{lemma}

\

\noindent
We will return to the proof of this lemma in section~\ref{ssec:MHT}
but first we will discuss a couple of  variants of this geometric setup. 

\

\noindent
{\bfseries (ii)} \ Suppose $X$ is a smooth projective $d$-dimensional
Calabi-Yau manifold. Let as before $\bOmega_{X}$ be the holomorphic
volume form on $X$. Let $D \subset X$ be a normal crossings
divisor. Typically if $D$ is smooth, it will be a variety of general
type. 

Consider the dg
Batalin-Vilkovisky algebra $(A,d,\Delta)$ given by
\[
\begin{split}
A & := \Gamma_{C^{\infty}}\left(X,\wedge^{\bullet}
T_{X,D}\otimes_{\mathcal{C}_{X}^{\infty}} 
\wedge^{\bullet}A^{0,1}_{X}\right) \\
d & := \bar{\partial} \\
\Delta & := \ddiv_{\bOmega_{X}} = \iota_{\bOmega_{X}}^{-1}\circ
\partial \circ \iota_{\bOmega_{X}},
\end{split}
\]
The contraction $\iota_{\bOmega_{X}}$ identifies this algebra with the
dg Batalin-Vilkovisky algebra
\[
\left(\Gamma_{C^{\infty}}\left(X,\Omega^{\bullet}_{X}(\op{rel} D)
\otimes_{\mathcal{C}_{X}^{\infty}}
\wedge^{\bullet}A^{0,1}_{X}\right),\bar{\partial},\partial\right),
\]
where $\Omega^{k}_{X}(\op{rel} D) \subset \Omega^{k}_{X}$ denotes the
subsheaf of all holomorphic $k$-forms that restrict to $0 \in
\Omega^{k}_{D-\op{sing}( D)}$. The cohomology of the total complex
  associated with this double complex is the de Rham cohomology of the
  pair $(X,D)$, and so again we get an identification
\begin{equation} \label{eq:iso(ii)}
H^{\bullet}(A,d+u\Delta) \cong H^{\bullet}_{dR}(X,D;\mathbb{C})
\end{equation}
valid for all fixed $u \neq 0$. Again using this identification and
mixed Hodge theory one deduces the following

\begin{lemma} \label{lemma::deg.with.D(ii)}
The dg Batalin-Vilkovisky algebra $(A,d,\Delta)$ has the degeneration
property and hence the formal moduli space of the pair
$(X,D)$ is smooth.
\end{lemma}

\

\noindent
{\bfseries (iii)} \ The setups {\bfseries (i)} and {\bfseries (ii)}
have a natural common generalization. Fix a smooth projective complex
variety of dimension $d$, a normal crossings divisor $D =
\cup_{i \in I} D_{i} \subset X$, and a collection of weights
$\{a_{i}\}_{i \in I} \subset [0,1]\cap \mathbb{Q}$, so that
\[
\sum_{i\in I}\displaylimits a_{i}[D_{i}] = - K_{X} \in
\op{Pic}(X)\otimes \mathbb{Q}.
\]
Represent the $a_{i}$'s by reduced fractions, take $N \geq 1$ to be
the least common multiple of the denominators of these fractions and
such that 
\[
\sum_{i\in I} (Na_i) [D_i] = -N K_X \in \op{Pic}(X),
\]
and
set $n_{i} := a_{i}N$. In particular we have a unique up to scale
section $\widetilde{\bOmega}_{X} \in \Gamma\left(X,K_{X}^{\otimes
  (-N)}\right)$ whose divisor is $\sum_{i \in I} n_{i}D_{i}$.
In this situation we can again promote the Dolbeault dg Lie algebra 
which
computes the deformation theory of $(X,D)$ to a  dg Batalin-Vilkovisky algebra 
$(A,d,\Delta)$, where 
\[
\begin{split}
A & := \Gamma_{C^{\infty}}\left(X,\wedge^{\bullet}
T_{X,D}\otimes_{\mathcal{C}_{X}^{\infty}} 
\wedge^{\bullet}A^{0,1}_{X}\right) \\
d & := \bar{\partial} \\
\Delta & := \ddiv_{\widetilde{\bOmega}_{X}}.
\end{split}
\]
The divergence operator $\ddiv_{\widetilde{\bOmega}_{X}}$ is defined
as follows. Restricting the section $\widetilde{\bOmega}_{X}$ to $X -
D$ we get a nowhere vanishing section of $K_{X-D}^{\otimes (-N)}$,
i.e. a flat holomorphic connection on $K_{X-D}$. If $ U \subset X-D$
is a simply connected open, then we can choose $\bOmega_{U}$ a
holomorphic volume form on $U$ which is covariantly constant for this flat
connection, and define the associated divergence operator
$\ddiv_{\bOmega_{U}} := \iota_{\bOmega_{U}}^{-1}\circ \partial \circ
\iota_{\bOmega_{U}}$. But by the flatness of the connection it follows
that any other covariantly constant volume form
on $U$ will be proportional to $\bOmega_{U}$ with a constant
proportionality coefficient. Since by definition $\ddiv_{c\bOmega_{U}}
= \ddiv_{\bOmega_{U}}$ for any constant $c$ we get a well defined
divergence operator on $X-U$. Furthermore locally this divergence
operator is a given by a holomorphic volume form which is a branch of 
$\left(\widetilde{\bOmega}_{X}\right)^{-1/N}$ and so by continuity it gives a
well defined map of locally free sheaves 
$ \ddiv_{\widetilde{\bOmega}_{X}} : \wedge^{i}T_{X,D} \to
\wedge^{i-1}T_{X,D}$. 

Again we claim that 

\begin{lemma} \label{lemma::deg.with.D(iii)}
The dg Batalin-Vilkovisky algebra $(A,d,\Delta)$ has the degeneration
property and the formal moduli space of the pair
$(X,D)$ is smooth.
\end{lemma}
{\bfseries Proof.} \ The proof of this lemma again reduces to mixed
Hodge theory via a map similar to the isomorphisms \eqref{eq:iso(i)}
and \eqref{eq:iso(ii)}.  However constructing this map is a bit more
involved than the arguments we used to construct \eqref{eq:iso(i)}
and \eqref{eq:iso(ii)}.

Consider the root stack $Z = X\left\langle \left\{
\frac{D_{i}}{N}\right\}_{i \in I} \right\rangle$ as defined in
e.g. \cite{matsuki-olsson,iyer-simpson}. By construction $Z$ is a
smooth proper Deligne-Mumford stack, equipped with a finite and flat
morphism $\pi : Z \to X$.

Conceptually the best way to define the stack $Z$  is  as a moduli stack
classifying (special) log structures associated with $X$, the divisor
$D$ and the number $N$ (see
\cite{matsuki-olsson} for the details). 
Etale locally on $X$ the stack $Z$ can be described easily as a
quotient stack. Indeed choose etale locally an identification of $X$
with a neighborhood of zero in $\mathbb{A}^{d}$ with coordinates
$z_{1}, \ldots, z_{d}$, so that $D = D_{1}\cup \cdots \cup D_{r}$ and
$D_{i}$ is identified with the hyperplane $z_{i} = 0$. Then the
corresponding etale local patch in $Z$ is canonically isomorphic to
the stack quotient 
\[
\big[\mathbb{A}^{d}/\underbrace{\bmu_{N}\times \cdots \times
  \bmu_{N}}_{\text{$r$-times}}\big],
\] 
where $\bmu_{N} \subset \mathbb{C}^{\times}$ is the group of $N$-th
roots of unity, and $\left(\zeta_{1},\ldots,\zeta_{r}\right) \in
\bmu_{N}^{\times r}$ acts as $\left(z_{1}, \ldots, z_{r},z_{r+1},
\ldots,z_{d}\right) \mapsto \left(\zeta_{1}z_{1}, \ldots,
\zeta_{r}z_{r}, z_{r+1}, \ldots, z_{d}\right)$.

In particular, this description shows (see
\cite[Theorem~4.1]{matsuki-olsson}) that:
\begin{itemize}
\item The map $\pi$ is an isomorphism over
$X - D$ and in general  identifies $X$ with the
coarse moduli space of $Z$;
\item There is a strict normal crossings divisor $\widetilde{D} =
  \cup_{i \in I} \widetilde{D}_{i} \subset Z$, such that 
\[
\mathcal{O}_{Z}\left(-N\widetilde{D}_{i}\right) =
  \pi^{*}\mathcal{O}_{X}\left(-D_{i}\right)
\] 
as ideal subsheaves of
  $\mathcal{O}_{Z}$; 
\item For all $j$ we have the Hurwitz formula $\Omega_{Z}^{j}(\log
  \widetilde{D}) 
  = \pi^{*}\Omega^{j}_{X}(\log D)$.
\end{itemize}
In particular we have canonical isomorphisms
\[
\begin{split}
\pi^{*}K_{X} & \cong \mathcal{O}_{Z}\left(-\sum_{i \in I}
n_{i}\widetilde{D}_{i}\right) \\[+0.5pc]
\pi^{*}K_{X} & \cong  K_{Z}\otimes \mathcal{O}_{Z}\left((1-N)\sum_{i \in I}
\widetilde{D}_{i}\right) 
\end{split}
\]
the first given by the section
$\pi^{*}\widetilde{\bOmega}_{X}$ and the second coming from the
Hurwitz formula.

There is a natural complex local system of rank one on $X-D$ with
monodromy in $\bmu_N$ associated with the choices of $N$-th root of the
section $\widetilde{\bOmega}_{X}$.  It is easy to see that the
pullback of this local system admits a canonical extension (as a local
system) to $Z$, which we denote by $\bxi$.  Moreover, we have a
canonical meromorphic section $\bOmega_{Z}$ of
$K_{Z}\otimes_{\mathbb{C}} \bxi$ with divisor $\sum_{i \in I}
(N-1-n_{i})\widetilde{D}_{i}$. It is easy to check locally by using
the etale local description of $Z$ as a quotient stack the contraction
$\iota_{\bOmega_{Z}}$ gives a well defined isomorphism of locally free
sheaves:
\[
\xymatrix@1{\iota_{\bOmega_{Z}} : \hspace{-1.5pc} & \wedge^{j}
  T_{Z,\widetilde{D}} \ar[r]^-{\cong} & 
\Omega^{d-j}_{Z}\left(\log \widetilde{D}_{(1)}, \op{rel}
\widetilde{D}_{(0)}\right)\otimes_{\mathbb{C}} \bxi.
}
\]
Here
\[
\begin{split}
\widetilde{D}_{(0)} & := \cup_{i \in I_{0}} \widetilde{D}_{i}  \qquad
\qquad I_{0} = \{ i \in I | a_{i} = 0 \} \\
\widetilde{D}_{(1)} & := \cup_{i \in I_{1}} \widetilde{D}_{i}  \qquad
\qquad I_{1} = \{ i \in I | a_{i} = 1 \}.
\end{split}
\]
Now taking into account the Hurwitz isomorphism  $\wedge^{j}
  T_{Z,\widetilde{D}} \cong \pi^{*}\wedge^{j}
  T_{X,D}$ and using adjunction, we can view $\iota_{\bOmega_{Z}}$
as an isomorphism
\begin{equation} \label{eq:adjoint.contract}
\xymatrix@1{\wedge^{j} T_{X,D}\ar[r]^-{\cong} & 
 \left(\pi_{*}\Omega^{d-j}_{Z}\left(\log
\widetilde{D}_{(1)}, \op{rel} 
\widetilde{D}_{(0)}\right)\otimes_{\mathbb{C}} \bxi\right)
}
\end{equation}
It is immediate from the definition that the isomorphism
\eqref{eq:adjoint.contract} (taken for all $j$) identifies the dg
Batalin-Vilkovisky algebra $(A,d,\Delta)$ with the Dolbeault bicomplex
\[
\left(\Gamma_{C^{\infty}}\left(X,\left(
\pi_{*}\Omega^{\bullet}_{Z}\left(\log \widetilde{D}_{(1)}, \op{rel}
\widetilde{D}_{(0)}\right)\otimes_{\mathbb{C}} \bxi\right)
\otimes_{\mathcal{C}^{\infty}_{X}}A^{0,\bullet}_{X} 
\right),\bar{\partial},\partial\right).
\]
But the above complex
equipped with the differential $\partial + \bar{\partial}$ is the
Dolbeault resolution of the complex of sheaves
$\pi_{*}\left(\Omega^{\bullet}_{Z}\left(\log \widetilde{D}_{(1)},
\op{rel} \widetilde{D}_{(0)}\right)\otimes_{\mathbb{C}} \bxi,\partial\right)$ 
which is equal to the derived direct image 
$R\pi_{*}\left(\Omega^{\bullet}_{Z}\left(\log \widetilde{D}_{(1)},
\op{rel} \widetilde{D}_{(0)}\right)\otimes_{\mathbb{C}} \bxi,
\partial\right)$ since $\pi$ is
finite. Now combined with the Leray spectral sequence for $\pi$ this
gives, for all $u \neq 0$ an isomorphism 
\begin{equation} \label{eq:iso(iii)}
H^{\bullet}(A,d+u\Delta) \cong
H^{\bullet}_{dR}\left(Z-\widetilde{D}_{(1)},\widetilde{D}_{(0)} -
\widetilde{D}_{(1)}; \bxi\right),
\end{equation}
which specializes to  both isomorphisms \eqref{eq:iso(i)} and
\eqref{eq:iso(ii)}. 

Now the fact that $Z$ is a smooth and proper Deligne-Mumford stack and
mixed Hodge theory (see \ref{ssec:MHT}) for
$\left(Z-\widetilde{D}_{(1)},\widetilde{D}_{(0)} -
\widetilde{D}_{(1)}\right)$ endowed with local system
 $\bxi$ imply that $(A,d,\Delta)$ has the
degeneration property. \ \hfill $\Box$

\

\begin{rem} \label{rem:deg.mirror} The fact that the root stack in the
  previous proof can be viewed as the moduli stack of special log
  structures is very interesting. It suggests that the setup we just
  discussed may fit naturally in the recent approach of Gross-Siebert
  \cite{gross-siebert.i,gross-siebert.ii} to mirror symmetry and
  instanton corrections via log degenerations of toric Fano manifolds
  (see also \cite{ks-nonarchimedean,kontsevich.soibelman-tori}). The
  relationship between these two setups is certainly worth studying
  and we plan to return to it in the future.
\end{rem}

\

\noindent
{\bfseries (iv)} Yet another generalization of the previous picture arises
  when we take the variety $X$ to be  a normal-crossings Calabi-Yau. More
  precisely assume that $X$ is a strict normal crossings variety with
  irreducible components $X = \cup_{i \in I} X_{i}$ equipped with a
  holomorphic volume form $\bOmega_{X}$ on $X - X^{\op{sing}}$ such that 
  the restriction of  $\bOmega_{X}$ on each $X_{i}$ has a logarithmic
  pole along $X_{i}\cap \left( 
  \cup_{j\neq i} X_{j}\right)$ and the residues of these restricted
  forms cancel along each $X_{i}\cup X_{j}$. Taking a colimit along
  the projective system of all finite intersections of components of
  $X$ we get again a dg
  Batalin-Vilkovisky algebra $A_{\op{tot}}(X) = \op{{\sf
  colim}}_{J\subset I} A\left(\cap_{i \in J} X_{i}\right)$ and again
  by using mixed Hodge theory we can check that this algebra has the
  degeneration property. 

\

\medskip

\subsubsection {\bfseries Mixed Hodge theory in a nutshell.} \
\label{ssec:MHT} In this section we briefly recall the basic arguments from
Deligne's mixed Hodge theory \cite{deligne-h3} that are necessary for
proving the degeneration property of the dg Batalin-Vilkovisky algebras
in section~\ref{ssec:generalize} {\bfseries (i)-(iv)}. 

Suppose we are given:
\begin{itemize}
\item  a finite ordered collection $\left( X_{\alpha} \right)$ of
smooth complex projective varieties;
\item for every $\alpha$ a choice of a $\mathbb{Z}\times
  \mathbb{Z}$-graded complex of sheaves of differential forms which are either
  $C^{\infty}$ or are $C^{-\infty}$ (i.e. currents) and constrained so
  that their wave 
  front (singular support) is contained in a given conical Lagrangian in
  $T^{\vee}X_{\alpha}$ which is the conormal bundle to a normal
  crossings divisor in $X_{\alpha}$;
\item a collection of integers $n_{\alpha} \in \mathbb{Z}$.
\end{itemize}
Consider the complex $C^{\op{tot}} = \oplus_{\alpha}
C_{\alpha}^{\bullet}[n_{\alpha}]$  equipped with three differentials
$\partial$, $\bar{\partial}$, $\delta$, where $\delta = \sum_{\alpha < \beta}
\delta_{\alpha\beta}$, and the $\delta_{\alpha\beta}$ come from
pullbacks and pushforwards for some maps $X_{\beta} \hookrightarrow
X_{\alpha}$ or $X_{\alpha} \hookrightarrow X_{\beta}$. The statement
we need now can be formulated as follows:

\begin{claim} \label{claim:MHT} For every $k \geq 1$ the cohomology
\[
H^{\bullet}\left(C^{\op{tot}}[u]/(u^{k}), \bar{\partial} + \delta +
  u\partial\right)
\]
 is a free $\mathbb{C}[u]/(u^{k})$-module. 
\end{claim}
{\bfseries Proof.} \ If $X$ is smooth projective over $\mathbb{C}$ and
if $\left( A^{\bullet}(X), \bar{\partial} \right)$ is the
$\bar{\partial}$-complex of (either $C^{\infty}$ or $C^{-\infty}$)
differential forms on $X$, then the inclusion
\[
\left( \ker \partial, \bar{\partial} \right) \hookrightarrow \left(
A^{\bullet}(X), \bar{\partial} \right) 
\]
is a quasi-isomorphism. 

This implies that the horizontal arrows in the diagram of complexes
\[
\xymatrix@R-1pc{
\left( \ker \partial[u]/(u^{k}), \bar{\partial} + \delta + u \partial
\right) \ar@{=}[d] \ar[r] & \left(C^{\op{tot}}[u]/(u^{k}),
\bar{\partial} + \delta + 
  u\partial\right) \\
\left( \ker \partial[u]/(u^{k}), \bar{\partial} + \delta 
\right) \ar[r] & \left(C^{\op{tot}},
\bar{\partial} + \delta\right)[u]/(u^{k}),
}
\]
are quasi-isomorphisms. Indeed, this follows by noticing that there
are natural filtrations on both sides (by the powers of $u$ and the
index $\alpha$) which give rise to convergent spectral sequences and
induce the quasi-isomorphic inclusion $\left( \ker \partial,
\bar{\partial} \right) \hookrightarrow
\left(C^{\op{tot}},\bar{\partial}\right)$ on the associated
graded. This proves the claim. \ \hfill $\Box$

\

\begin{rem} \label{rem:MHT} $\bullet$ \ 
Note that the same reasoning implies
  that the natural map
\[
\left( \ker \partial,
\bar{\partial} + \delta \right) \twoheadrightarrow \left( \ker
\partial/\op{im} \partial, \bar{\partial}+\delta \right) = \left(
H^{\bullet}(X_{\alpha}), \delta \right),
\]
is also a quasi-isomorphism, which reduces the problem of computing
\linebreak 
$H^{\bullet}\left(C^{\op{tot}}[u]/(u^{k}), \bar{\partial} + \delta +
  u\partial\right)$ to a homological algebra question on a complex of
  finite dimensional vector spaces.

\

\noindent
$\bullet$ \ There is useful variant of the theory, also discussed in
 \cite{deligne-h3}: the previous discussion immediately generalizes to
 the case of cochain complexes of a collection of projective manifolds
 with coefficients in some unitary local systems.
\end{rem}

\

\noindent
Next we discuss a few examples and applications of the geometric setup from
section~\ref{ssec:generalize}.

\

\subsubsection {\bfseries The moduli stack of Fano varieties.} \ As a
consequence of section~\ref{ssec:generalize} {\bfseries (iii)} we get
a new proof and a refinement of the following result of Ran
\cite{ran,kawamata}:

\begin{theo} \label{theo:ran} Let $X$ be a complex Fano
manifold, that is let $X$ be a smooth proper $\mathbb{C}$-variety for
which $K_{X}^{-1}$ is ample. Then  
the versal deformations of $X$ are unobstructed. 
\end{theo}
{\bfseries Proof:} Choose $N > 1$ so that $K_{X}^{\otimes (-N)}$ is
very ample and all the higher cohomology groups $H^{k}\left(X,
K_X^{\otimes (-N)}\right)$ vanish for $k\geq 1$.  Choose a generic
section $\widetilde{\bOmega}_{X} \in H^{0}\left(X,K_{X}^{\otimes
(-N)}\right) = 0$ whose zero locus is a smooth and connected divisor
$D \subset X$.

Consider now $\mathfrak{g} = \bPi R\Gamma\left(X,\wedge^{\bullet}
T_{X,D}\right)$ with the Schouten bracket. By
Lemma~\ref{lemma::deg.with.D(iii)} this d$(\mathbb{Z}/2)$g Lie algebra
is homotopy abelian and so as in the proof of
Theorem~\ref{theo:Tian.Todorov} we conclude that $\mathfrak{g}^{(1)} =
R\Gamma\left(X,T_{X,D}\right)$ is homotopy abelian. Since this
d$(\mathbb{Z}/2)$g Lie algebra governs the deformation theory of
$(X,D)$ as a variety with a divisor, it follows that the formal germ
of the deformation space of the pair $(X,D)$ is smooth. Next we will
need the following simple

\begin{lemma} \label{lemma:volume.form}
Suppose $(X',D')$ is a small deformation of $(X,D)$ as a variety with
divisor. Then $X'$ is still a Fano with $K_{X'}^{\otimes (-N)}$ is
very ample and $D' \in \left|K_{X'}^{\otimes (-N)}\right|$.
\end{lemma}
{\bfseries Proof:} The condition of $K_{X}^{\otimes (-N)}$ being very
ample is open in the moduli of $X$. Furthermore by definition 
 $K_{X}^{\otimes (-N)}\otimes \mathcal{O}_{X}(-D) = \mathcal{O}_{X}$
and so by the small deformation hypothesis it follows that 
 $K_{X'}^{\otimes (-N)}\otimes \mathcal{O}_{X'}(-D)$ is in the
connected component of the identity of $\op{Pic}(X')$. But $X'$ is a
Fano and so $\op{Lie}(\op{Pic}^{0}(X')) = H^{1}(X',\mathcal{O}_{X'}) =
0$. Hence  $K_{X'}^{\otimes (-N)}\otimes \mathcal{O}_{X'}(-D) =
\mathcal{O}_{X}$ as \linebreak well. \ \hfill $\Box$

\

\noindent
The theorem now follows easily. The versal deformation space of smooth
connected $D$'s for a given $X$ is smooth and isomorphic to a domain
in $\mathbb{P}^{h^{0}\left(X,K_{X'}^{\otimes (-N)}\right) - 1}$. Since
the dimension of these projective spaces is locally constant in $X$ by
Riemann-Roch and vanishing of the higher cohomologies, it follows that
the map from the versal deformation space of the pairs $(X,D)$ to the
versal deformation stack of $X$ is smooth. In other words the versal
deformation stack of $X$ has a presentation in the smooth topology
with a smooth atlas - the versal deformation space for $(X,D)$. Hence
the versal deformations of $X$ are a smooth stack. \ \hfill $\Box$

\

\medskip

\subsubsection {\bfseries Algebras for the Landau-Ginzburg model.} \ 
Consider again the setup of a holomorphic Landau-Ginzburg
model. Suppose $Y$ is smooth and quasi-projective over $\mathbb{C}$
and of dimension $\dim Y = d$. Suppose there exists a nowhere vanishing
algebraic volume form $\bOmega_{Y} \in \Gamma(Y,K_{Y})$, and let
$\bw : Y \to \mathbb{A}^{1}$ be a regular function with compact
critical locus. 

This data gives a dg Batalin-Vilkovisky algebra 
$(A,d,\Delta)$ where
\[
\begin{split}
A & := \Gamma_{C^{\infty}}\left(Y,\wedge^{\bullet}
T_{Y}^{1,0}\otimes_{\mathcal{C}_{Y}^{\infty}} 
\wedge^{\bullet}A^{0,1}_{Y}\right) \\
d & := \bar{\partial} + \iota_{d\bw}\\
\Delta & := \ddiv_{\widetilde{\bOmega}_{Y}}.
\end{split}
\]
Again the contraction $\iota_{\bOmega_{Y}}$ identifies $(A,d,\Delta)$
with the twisted Dolbeault bicomplex \linebreak 
$\left(A^{\bullet}(Y),\bar{\partial} + d\bw\wedge,
\partial\right)$. The latter satisfies the degeneration property by 
the work of Barannikov and the second author, Sabbah
\cite{sabbah-twisted}, or 
Ogus-Vologodsky \cite{ogus.vologodsky}

\

\begin{rem} It will be interesting to combine
the previous discussion with the discussion in
section~\ref{ssec:generalize} {\bfseries (iii)} or with the broken
Calabi-Yau geometry from section~\ref{ssec:generalize} {\bfseries
(iv)}.  Suppose we have a quasi-projective smooth complex $Y$, a
regular function $\bw : Y \to \mathbb{A}^{1}$ with compact critical
locus, and suppose we are given a normal crossings divisor $D =
\cup_{i \in I} D_{i}$ and a system of weights $\{ a_{i} \}_{i \in I}$
as in section~\ref{ssec:generalize} {\bfseries (iii)}. Then we can
write the $\bw$-twisted version of the dg Batalin-Vilkovisky algebra
for $(Y,D)$ which by general nonsense will compute the deformation
theory of the data $(Y,D,\bw)$. Similarly we can add a potential to a
$Y$ which itself is a normal-crossings Calabi-Yau, as in
section~\ref{ssec:generalize} {\bfseries (iv)}.
We expect
that the resulting algebras will again have the degeneration property but we
have not investigated this question.
\end{rem}

\

\subsection{\bfseries Categorical framework: spherical functors}
\label{ssec:spherical}

In this section we briefly discuss some algebraic aspects of the
deformation theory of \nc-spaces (see
section~\ref{sssec:categorical}). For simplicity we will discuss the
$\mathbb{Z}$-graded case but in fact all definitions and statements readily 
generalize to the $\mathbb{Z}/2$ case.

\

\subsubsection {\bfseries Calabi-Yau \nc-spaces.} \ 
\label{sssec:CYnc}
Suppose $X = \ncSpec(A)$
is a graded \nc-affine \nc-space over $\mathbb{C}$. 
If $X$ is smooth, then
$A \in \Perf_{X\times X^{\op{op}}} = \Perf\left(A\otimes
A^{\op{op}}-\text{{\sf mod}}\right)$ and we define the {\em\bfseries smooth
dual} of $A$ to be $A^{!} := \op{Hom}_{A\otimes
A^{\op{op}}}\left(A,A\otimes A\right)$. 
Similarly if $X$ is compact,
then $A \in \Perf_{\op{pt}}$ and we define the {\em\bfseries compact
  dual} of $A$ to be $A^{*} := \op{Hom}_{\mathbb{C}}(A,\mathbb{C}) \in
\left(A\otimes
A^{\op{op}}-\text{{\sf mod}}\right)$. 

If $X$ is both a smooth and compact \nc-space, then we have isomorphisms 
\[
A^{!}\otimes_{A}
A^{*} \cong  A^{*}\otimes_{A}
A^{!} \cong A
\]
in the category $\left(A\otimes
A^{\op{op}}-\text{{\sf mod}}\right)$. The endofunctor $S_{X} : C_{X}
\to C_{X}$ given by the $A$-bimodule $A^{*}$ is called {\em\bfseries
  the Serre functor} of $X$. It is an autoequivalence of $C_{X}$ which
is central (i.e. commutes with all autoequivalences). Moreover for any two
objects $\mathcal{E}, \mathcal{F} \in \Perf_{X}$ there is a functorial
identification 
\[
\op{Hom}_{X}(\mathcal{E},\mathcal{F})^{\vee} \cong
\op{Hom}_{X}(\mathcal{F},S_{X}\mathcal{E}).
\]
With this notation we have the following definition (see also
\cite{ks-ncgeometry}):

\begin{defi} \label{defi:ncCY} We say that a smooth graded
\nc-affine \nc-space $X =
  \ncSpec(A)$ is a {\bfseries Calabi-Yau of dimension $d \in
\mathbb{Z}$} if $A^{!} \cong A[-d] \quad \text{ in } \quad
\left(A\otimes A^{\op{op}}-\text{{\sf mod}}\right)$.  We say that a
compact \nc-affine \nc-space $X = \ncSpec(A)$ is a {\bfseries
Calabi-Yau of dimension $d \in \mathbb{Z}$} if $A^{*} \cong A[d] \quad
\text{ in } \quad \left(A\otimes A^{\op{op}}-\text{{\sf mod}}\right)$.
\end{defi}

\

\noindent
The definition works also in the $\mathbb{Z}/2$-graded case, where the
dimension $d$ is 
understood as an element of $\mathbb{Z}/2$.

\

\noindent
For a \nc-space which is both smooth and compact the two conditions
are equivalent and are equivalent to having an isomorphism of
endofunctors $S_{X} \cong [d]$. 

\

\begin{rem} \label{rem:CYcondition} This definition of a Calabi-Yau
  structure on a smooth compact \nc-space is somewhat simplistic and
  should be taken with a grain of salt. The true definition (see
  \cite{ks-ncgeometry}) implies the isomorphism of functors $S_{X}
  \cong [d]$ but also involves higher homotopical data which is
  encoded in a cyclic category structure on $C_{X}$. We will suppress
  the cyclic structure here in order to simplify the discussion.
\end{rem}

\

\noindent
We are interested in \nc-space analogues of the Tian-Todorov
theorem. The unobstructedness of graded smooth and
compact \nc-Calabi-Yau spaces was recently  analyzed by Pandit
\cite{pranav} via the $T^{1}$-lifting property of Ran \cite{ran} and
Kawamata \cite{kawamata}. Here we formulate the following general 

\begin{theo} \label{theo:ncCY-TT} Suppose that $X$ is a smooth
  and compact \nc-Calabi-Yau space of dimension $d \in
  \mathbb{Z}$ (or of dimension $d \in \mathbb{Z}/2$ in the
  $\mathbb{Z}/2$-graded case). Assume that $X$ satisfies the degeneration
  conjecture (see section~\ref{sssec:degeneration}). Then:
\begin{itemize}
\item the Hochschild cochain algebra $C^{\bullet}(X)$ of $X$ is a
homotopy abelian $L_{\infty}$ algebra; 
\item the formal moduli space $\MM_{X}$ of $X$ is a formal
  supermanifold, i.e. 
\[
\MM_{X} := \MM_{C^{\bullet}(A,A)} \cong
  \op{Spf} \mathbb{C}[[x_{1}, \ldots, x_{N},\xi_{1}, \ldots,
  \xi_{M}]];
\] 
\item the negative cyclic homology of the universal family over
  $\MM_{X}$ gives a vector bundle $H \to \MM_{X}\times \mathbb{D}$
  which is equipped with a flat meromorphic connection $\nabla$ so
  that $\nabla_{u \partial/\partial x_{i}}$, $\nabla_{u
  \partial/\partial \xi_{j}}$, and  $\nabla_{u^{2} \partial/\partial
  u}$ are regular;
\item  $(H,\nabla)$ is the de Rham part of a
Calabi-Yau variation of \nc-Hodge structures.
\end{itemize}
\end{theo}

\

\noindent
We will only sketch some of the highlights of the proof of this
theorem here since going into full details will take us too far
afield. The proof is based on a mildly generalized version of
Deligne's conjecture (see e.g. \cite{ks-deligne,tamarkin}) which
states that the Hochschild cochain complex of an affine
\nc-space is also an algebra over the operad of chains of the little
discs operad. The first step is to show that under the Calabi-Yau
assumption the Hochschild cochain complex $C^{\bullet}(X)$  is also
naturally an algebra over the cyclic operad of chains of the framed
little discs operad (i.e. the operad of little discs with a marked
point point on the boundary). Next one shows that  the validity of the
degeneration conjecture for $X$ implies that the induced
$S^{1}$-action on the cochain complex, is homotopically
trivial. Finally by a topological argument one deduces from this the fact
that all the higher $L_{\infty}$ operations on $C^{\bullet}(X)$ must
vanish.

\

\begin{rem} It seems certain\footnote{We borrowed this delightful
    expression from \cite{misha}.} that from deformation quantization
    it follows that if $X$ is a smooth and projective Calabi-Yau
    variety, then the data described in the above theorem is
    canonically isomorphic to the formal completion of the variation
    of \nc-Hodge structures described in section~\ref{sssec:CYcan}. 

For a general smooth and compact \nc-Calabi-Yau space we expect that
the formal variation of \nc-de Rham data in theorem \ref{theo:ncCY-TT}
converges to give an analytic de Rham data which 
contains a compatible \nc-Betti data $\mycal{E}_{B}$ and so extends to
an honest variation of \nc-Hodge structures. 
\end{rem}

\

\medskip

\subsubsection {\bfseries Spherical functors.} \ In this section 
 we introduce a special version of the general notion of a spherical
functor \cite{anno} which is tailored to the Calabi-Yau condition. We
begin with a definition:

\begin{defi} \label{defi:ncmap} Let
$X$ and $Y$ be two graded \nc-spaces. A {\bfseries morphism $f : X
    \to Y$}  is a triple of functors 
\[
\xymatrix{
C_{X} \ar[d]|-{f_{*}} \\
C_{Y} \ar@<-1pc>[u]|-{f^{!}} \ar@<1pc>[u]|-{f^{*}}
}
\]
so that $(f^{*},f_{*})$ and $(f_{*},f^{!})$ are
$(\text{left},\text{right})$ pairs of adjoint functors. 
\end{defi}

\

\noindent
Suppose now $X$, $Y$ are smooth and compact graded \nc-spaces and let
$Y$ be a \nc-Calabi-Yau of dimension $d$.

\begin{defi} \label{defi-spherical} A morphism 
$f : X \to Y$ is {\bfseries spherical} if:
\begin{itemize}
\item[{\em\bfseries (a)}]  the cone of the
natural adjunction morphism $\op{id}_{C_{X}} \to f^{!}\circ f_{*}$ is
isomorphic to the shifted Serre functor of $X$:
$\op{cone}\left( \op{id}_{C_{X}} \to f^{!}\circ f_{*} \right) \cong
S_{X}[1-d]$,
\item[{\em\bfseries (b)}] the natural map 
$f^{!} \to S_{X}[1-d]\circ f^{*}$, 
induced from the isomorphism in {\em\bfseries (a)} and the adjunction 
$f^{!} \to f^{!}\circ f_{*}\circ f^{*}$ is an isomorphism of functors.
\end{itemize}
\end{defi}

\

\begin{rem} \label{rem:spherical}
{\bfseries (a)} \ If $f$ is spherical, then the associated
{\em\bfseries reflection functor} \linebreak 
$\mathcal{R}_{f} := \op{cone}\left( f_{*}\circ f^{!} \to
\op{id}_{C_{Y}}\right)$
is an auto-equivalence of $C_{Y}$ \cite{anno}. 

\noindent
{\bf (b)} \ Similarly to the definition of a Calabi-Yau structure the
above notion of a spherical functor should be viewed as a weak
preliminary version of a stronger more refined notion which has to
involve higher homotopical data and has yet to be defined carefully.
\end{rem}

\

\medskip

\begin{ex} \label{ex:ncpairs} 
{\bf (i)} \ Let $X = \op{pt}$, and let $Y$ be a $d$-dimensional
  smooth and compact \nc-Calabi-Yau and let $\mathcal{E} \in
  C_{Y}$ be a spherical object, i.e. an object for which the complex
  of $\mathbb{C}$-vector spaces
  $\op{Hom}_{Y}(\mathcal{E},\mathcal{E})$ is quasi isomorphic to
  $(H^{\bullet}(S^{d},\mathbb{C}),0)$. The morphism of
  \nc-spaces $f : \op{pt} \to Y$ given by $f_{*}(V)
  = \mathcal{E}\otimes V$, for any $V \in C_{\op{pt}} =
  \left(\text{\sf Vect}_{\mathbb{C}}\right)$ is spherical. 

\

\noindent
{\bf (ii)} \ Let $X$ be smooth and projective of dimension $d+1$, and
let $i : Y \hookrightarrow X$ be a smooth anti-canonical divisor in
$X$. The $Y$ is a $d$-dimensional Calabi-Yau and we have a natural spherical
\nc-morphism $f : X \to Y$ given by $f_{*} := i^{*}$, $f^{!} :=
i_{*}$, etc.

\

\noindent
{\bf (iii)} \ Let $Y$ be a smooth projective $d$-dimensional
Calabi-Yau. Let $i : X \hookrightarrow Y$  be a smooth hypersurface. 
Then we have a natural spherical \nc-morphism $f : X \to Y$ given by
$f_{*} = i_{*}$, $f^{!} = i^{!}$, and $f^{*} = i^{*}$. 
\end{ex}

\

\begin{rem}
The geometry of Example~\ref{ex:ncpairs} {\bfseries (ii)}, where $X$
is taken to be a smooth projective Fano, and $i : Y \hookrightarrow X$
is a smooth anti-canonical divisor, can be encoded algebraically in
the categories $C_{X} = D(\text{\sf Qcoh}(X))$, $C_{Y} = D(\text{\sf
Qcoh}(Y))$, the functor $f_{*} = i^{*}$, and another natural triple of
categories:
\begin{itemize}
\item the compact category 
$D_{\substack{\text{compact} \\ \text{support}}}(\text{\sf Qcoh}(X-Y)) =
  \ker(f_{*})$,
\item the compact category $D_{\op{supp}\, Y}(\text{\sf Qcoh}(X)) = $
 the subcategory in $D(\text{\sf Qcoh}(X))$ generated by
 $i_{*}D(\text{\sf Qcoh}(Y))$,
\item the smooth category $D(\text{\sf Qcoh}(X-Y)) =$ the quotient 
$D(\text{\sf Qcoh}(X))/D_{\op{supp}\, Y}(\text{\sf Qcoh}(X))$.
\end{itemize}
There is a similar triple of categories for the setup in
Example~\ref{ex:ncpairs} {\bfseries (iii)}. It will be very
interesting to describe the categorical data that encodes
anti-canonical divisors
with normal crossings or more generally the
fractional anti-canonical divisor setup from
section~\ref{ssec:generalize} {\bfseries (iii)}. It seems
likely that in this situation one gets a system of nested categories
and functors with a ``spherical'' condition imposed on the whole
system rather than on individual functors. This is a very interesting
question that we plan to investigate in the future.
\end{rem}

\

\begin{rem} It is clear from the examples above that spherical
  functors give a unifying framework for handling different type of
  geometric pairs. 

Suppose that $X$ and $Y$ are smooth and compact \nc-spaces, $Y$ is a
  \nc-Calabi-Yau, $f : X \to Y$ is a spherical map, and the
  degeneration conjecture holds for both $X$ and $Y$. In this
  situation we expect that the deformation theory of $f : X \to Y$ is
  controlled by a homotopy abelian d$(\mathbb{Z}/2)$g Lie algebra
  which is $L_{\infty}$-quasi-isomorphic to
\begin{equation} \label{eq:Linfcomplex}
\op{cone}\left(\xymatrix@1@C-1pc{C_{\bullet}(Y) \ar[r]^-{f^{!}} &
  C_{\bullet}(X)}\right)[1-d].
\end{equation}
Moreover, using $f_{*}$ (or $f^{!}$) we can build a new \nc-space $Z$
by taking $C_{Z}$ to be the semi-orthogonal extension $C_{Z} =
\left\langle C_{X}, C_{Y} \right\rangle$, where we set
\[
\begin{split}
\op{Hom}_{Z}\left(C_{Y},C_{X}\right) & := 0 \\[0.5pc]
\op{Hom}_{Z}\left(\mathcal{E},\mathcal{F}\right) & :=
\op{Hom}_{Y}\left(f_{*}\mathcal{E},\mathcal{F}\right) \quad \text{for
  all } \mathcal{E} \in C_{X}, \mathcal{F} \in C_{Y}.
\end{split}
\]
We expect that the deformation theory of $f : X \to Y$ is equivalent
to the deformation theory of $Z$ and in particular that the
$L_{\infty}$ algebra $C^{\bullet}(Z)$ is quasi-isomorphic to the
algebra \eqref{eq:Linfcomplex}.
\end{rem}

\

\begin{rem} We should point out that even though
  deformation quantization provides a conceptual bridge between the
  categorical framework and the geometric framework of the previous
  section, the actual connection between the two frameworks is tenuous
  at best. The source of the problem lies in the fact that the
  deformation quantization of general Poisson maps can be obstructed
  \cite{willwacher}.
\end{rem}

\subsection{\bfseries $A$-model framework: symplectic Landau-Ginzburg models}
\label{ssec:symplecticLG}

We already noted in Examples~\ref{ex:extended.symplectic} and
\ref{ex:can.basic} that there are natural canonical coordinates and a
Calabi-Yau variation of \nc-Hodge structures that one can attach to
the $A$-model on a compact symplectic manifold. An interesting open
problem is to find an algebraic description of these coordinates and
variation in terms of some $d(\mathbb{Z}/2)g$ Batalin-Vilkovisky algebra
that is naturally attached to the Fukaya category. This question is
hard and we will not study it directly here. Instead we will look at
the question of finding canonical coordinates and variation in another
symplectic context, i.e. for symplectic Landau-Ginzburg models, and
try to get an insight into a possible algebraic formulation in that
case. It will be interesting to compare our formalism with
the recent work of Fan-Jarvis-Ruan \cite{fan-jarvis-ruan} on the
symplectic geometry of 
quasi-homogeneous Landau-Ginzburg potentials with isolated
singularities but at the moment we do not see a direct relationship.

\subsubsection {\bfseries Symplectic geometry with potentials.} \

The objects we would like to understand are triples $(Y,\bw,\omega)$,
where
\begin{itemize}
\item $Y$ is a $C^{\infty}$-manifold and $\omega$ is a
  $C^{\infty}$-symplectic form on $Y$.
\item $\bw : Y \to \mathbb{C}$ is a proper $C^{\infty}$-map such that
  there exists an $R > 0$ so that over \linebreak 
$\{ z \in \mathbb{C} | \, |z|
  \geq R \}$ the map $\bw$ is a smooth fibration with fibers symplectic
  submanifolds in $(Y,\omega)$. 
\end{itemize} 
Similarly to the case of compact symplectic manifolds one can associate
Gromov-Witten invariants to such a geometry. Specifically, if we fix
$n \geq 1$, $g \geq 0$, and $\beta \in H_{2}(Y,\mathbb{Z})$, then we
can use stable pseudo-holomorphic pointed curves in $Y$ to 
define a natural linear (correlator) map
\[
\xymatrix@1{
I_{g,\beta,n-1}^{(1)} : & \hspace{-1pc}  H^{\bullet}(Y,\mathbb{Q})^{\otimes
  (n-1)}\otimes H^{\bullet}\left(\overline{\mathcal{M}}_{g,n},\mathbb{Q}\right)
  \ar[r] &  H^{\bullet}(Y,\mathbb{Q}). 
}
\] 
Indeed, note that Poincar\'{e} duality gives an identification
\[
H^{\bullet}(Y,\mathbb{Q}) \cong H_{\bullet}(Y,Y_{R};\mathbb{Q})[-\dim Y],
\]
where $R > 0$ is as above and $Y_{R} = \bw^{-1}\left( \{ z \in
  \mathbb{C} | \, |z| 
  \geq R \} \right) \subset Y$. Combining  this identification  with
  the isomorphism $(H^{\bullet})^{\vee} = H_{\bullet}$ we see that
  $I_{g,\beta,n-1}^{(1)}$ will be given by a class in 
$H_{\bullet}(Y,\mathbb{Q})^{\otimes (n-1)}\otimes
H_{\bullet}(Y,Y_{R};\mathbb{Q})\otimes
H_{\bullet}(\overline{\mathcal{M}}_{g,n},\mathbb{Q})$.

Next consider the usual moduli stack
$\overline{\mathcal{M}}_{g,n}(Y,\beta)$ of stable pseudo-holomorphic
maps. Here it will be convenient to assume that an almost-complex structure
 on $Y$ tamed by $\omega$ is chosen in such a way that $\bw_{|Y_{R}}$ is
holomorphic. The stack $\overline{\mathcal{M}}_{g,n}(Y,\beta)$ is non compact
but near infinity it parameterizes only pseudo-holomorphic maps
$\varphi : C \to Y$ such that $\bw\circ \varphi : C \to \mathbb{C}$ is
constant and $\bw\circ \varphi(C) \in \mathbb{C}$ is close to
infinity. Thus the virtual fundamental class of
$\overline{\mathcal{M}}_{g,n}(Y,\beta)$ is well defined as a class in
the relative homology 
\[
\begin{split}
\left[\overline{\mathcal{M}}_{g,n}(Y,\beta)\right]_{\op{vir}} & \in 
H_{\bullet}\left(Y^{n}\times \overline{\mathcal{M}}_{g,n}, Y^{n-1}\times
Y_{R} \times \overline{\mathcal{M}}_{g,n}; \mathbb{Q}\right) \\ 
& = 
H_{\bullet}(Y,\mathbb{Q})^{\otimes (n-1)}\otimes
H_{\bullet}(Y,Y_{R};\mathbb{Q})\otimes
H_{\bullet}(\overline{\mathcal{M}}_{g,n},\mathbb{Q}).
\end{split}
\]
We define $I_{g,\beta,n-1}^{(1)}$ to be the map given by the relative
virtual fundamental class
$\left[\overline{\mathcal{M}}_{g,n}(Y,\beta)\right]_{\op{vir}}$.

This collection of correlates satisfies analogues of the usual axioms
of a cohomological field theory \cite{kontsevich-manin} but we will
not discuss them here. Consider now a cohomology class
\[
x = (x_{2},x_{\neq 2}) \in H^{\bullet}(Y,\mathbb{C}) =
H^{2}(Y,\mathbb{C})\oplus H^{\neq 2}(Y,\mathbb{C}),
\] 
where $H^{\bullet}(Y,\mathbb{C})$ is viewed as a supermanifold over
$\mathbb{C}$. 

Now for every such $x$ we define a quantum product
\[
\xymatrix@1{\bullet *_{x} \bullet : & \hspace{-1pc}
  H^{\bullet}(Y,\mathbb{C})\otimes H^{\bullet}(Y,\mathbb{C}) \ar[r] &
  H^{\bullet}(Y,\mathbb{C})
}
\]
by the formula \label{page:qproduct}
\[
\alpha_{1} *_{x} \alpha_{2} := \sum_{m \geq 0}\displaylimits 
\sum_{\beta \in H_{2}(Y,\mathbb{Z})}\displaylimits 
\exp \left( \langle \beta, x_{2} \rangle \right)\cdot  \frac{1}{m!} 
I_{g,\beta,m+1}^{(1)} \big( (\alpha_{1}\otimes \alpha_{2} \otimes
\underbrace{x_{\neq 2} \otimes \cdots \otimes x_{\neq 2}}_{m \text{
    times }})\otimes 1_{\overline{\mathcal{M}}_{0,m+1}} \big).
\]
Now this quantum multiplication together with the usual formulas (see
Examples~\ref{ex:extended.symplectic} and \ref{ex:can.basic}) can be
used to define a decorated variation of \nc-Hodge structures over the
(conjecturally non-empty) domain in $H^{\bullet}(Y,\mathbb{C})$ where the
series defining $*_{x}$ is absolutely convergent.

\

\begin{rem} There are some interesting variants of this
  construction. For instance we can take a symplectic manifold
  $(Y,\omega)$ with no potential and a pseudo-convex boundary. In this
  situation $\overline{\mycal{M}}_{g,n}(Y,\beta)$ is already compact,
  as long as $\beta \neq 0$. Also in a symplectic Landau-Ginzburg
  model $(Y,\omega,\bw)$ we can allow for $\bw$ to be non-proper and
  instead require that its fibers have pseudo-convex boundary. Finally
  one can consider a symplectic $Y$ equipped with a proper map $Y \to
  \mathbb{C}^{k}$, holomorphic at infinity and with $k \geq 2$. 
\end{rem}

\

\subsubsection {\bfseries Categories of branes} \ Let $(Y,\omega,\bw)$
be a symplectic geometry with a proper potential. There are two
natural categories that we can attach to this geometry: the Fukaya
category of the general fiber of $\bw$, and the Fukaya-Seidel category
of $\bw$. Understanding the structure properties of these categories
or even defining them properly is a difficult task which requires a
lot of effort and hard work. We will not explain any of these
intricate details but will rather use the Fukaya and Fukaya-Seidel
categories as conceptual entities. For details of the definitions and
a rigorous development of the theory we refer the reader to the main
sources \cite{fooo,fukaya-ono},
\cite{seidel-book,seidel-directed}. The categories that we are
interested in are:

\

\begin{figure}[!ht]
\begin{center}
\psfrag{L1}[c][c][1][0]{{$L_{1}$}}
\psfrag{L2}[c][c][1][0]{$L_{2}$}
\epsfig{file=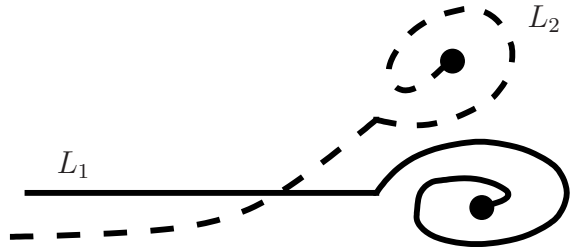,width=3in} 
\end{center}
\caption{Tadpole-like $\bw$-images of two Lagrangian
submanifolds.}\label{fig:wiggle}
\end{figure}

\noindent
{\bfseries (1)} \  The Fukaya-Seidel category $\FS(Y,\omega,\bw)$ of the
potential $\bw$ has objects which are unitary local systems $\mathbb{V}$ on
(graded)  $\omega$-Lagrangian submanifolds $L
\subset Y$ such that:
\begin{itemize}
\item $\bw(L) \subset (\text{compact})\cup \mathbb{R}_{\leq 0}$;
\item The restriction of $L$ over the ray $\mathbb{R}_{\leq 0}$ is a
  fibration on $\mathbb{R}_{\leq -R}$ and when $z \in \mathbb{R}_{\leq
  0}$, and $z \to -\infty$, we have that the fiber $L_{z} \subset
  Y_{z}$ is a Lagrangian submanifold in the symplectic manifold
  $\left(Y_{z},\omega_{|Y_{z}}\right)$. 
\end{itemize}
The morphisms between two objects $(L_{1},\mathbb{V}_{1})$ and
$(L_{2},\mathbb{V}_{2})$ are defined as homomorphisms between the
fibers of the local systems at the intersection points of the two
Lagrangians. As usual to make this work one has to perturb one of the
Lagrangians, say $L_{2}$ by a Hamiltonian isotopy to ensure
transversality of the intersection. A new feature of this setup
(compared to the situation of symplectic manifolds with no potential)
is that the allowable isotopies are tightly controlled - they
correspond to small wiggling, see Figure~\ref{fig:wiggle}, of the tail
of the tadpole-like image of the Lagrangian in $\mathbb{C}$ and a lift
of this wiggling to $Y$ given by a non-linear symplectic connection
identifying the fibers of $Y$.

\noindent
The compositions of morphisms are given by correlators counting
pseudo-holomorphic discs whose boundary is contained in the given
Lagrangian submanifolds.

\

\smallskip

\noindent
{\bfseries (2)} \ The Fukaya category $\Fuk(Y_{z})$ of a fiber
$\left(Y_{z},\omega_{|Y_{z}}\right)$ over a point  $z \in \mathbb{C}$
which is not a critical value for $\bw$. The objects in this category
are again pairs consisting of (graded) Lagrangian submanifolds in
$Y_{z}$ equipped with unitary local systems, and morphisms and
compositions are defined again by maps between the fibers of the local
systems at the intersection points and by counting discs. 
The parallel transport w.r.t. a non-linear symplectic connection on
$\bw : Y  \to \mathbb{C}$ identifies symplectically all fibers
$\left(Y_{z},\omega_{|Y_{z}}\right)$ over points $z \in
\mathbb{R}_{\leq 0}$  when $z \to -\infty$. We will denote any one
such fiber as $(Y_{-\infty},\omega_{-\infty})$.

\

\smallskip

\noindent
Now observe that by intersecting a Lagrangian $L \subset Y$ with the
fiber $Y_{-\infty}$ we get an assignment $L \mapsto L_{-\infty} :=
L\cap Y_{-\infty}$. We expect that this assignment can be promoted to
a spherical functor (see also \cite{seidel-directed} for a similar 
discussion)
\[
\xymatrix@1{ F : & \hspace{-1.5pc}  \FS(Y,\bw,\omega) \ar[r] &
  \Fuk(Y_{-\infty},\omega_{-\infty})
}
\]
so that the associated spherical twist $\mathcal{R}_{F} :
\Fuk(Y_{-\infty},\omega_{-\infty}) \to
\Fuk(Y_{-\infty},\omega_{-\infty})$ categorifies the monodromy around
the circle $\{ z \in \mathbb{C} | \; |z| = R \}$. 

In this situation one can also define relative Gromov-Witten invariants 
\[
\xymatrix@1{
J_{g,\beta,n-2}^{(1)} : & \hspace{-1pc}  H^{\bullet}(Y,Y_{-\infty};
\mathbb{Q}) \otimes H^{\bullet}(Y, \mathbb{Q})^{\otimes
  (n-2)}\otimes H^{\bullet}\left(\overline{\mathcal{M}}_{g,n},\mathbb{Q}\right)
  \ar[r] &  H^{\bullet}(Y,Y_{-\infty};\mathbb{Q}). 
}
\]
For we again use the duality $(H^{\bullet})^{\vee} \cong H_{\bullet}$ and the
Poincar\'{e} duality $H^{\bullet}(Y,Y_{-\infty};\mathbb{Q}) \cong
H_{\bullet}(Y,Y_{+\infty};\mathbb{Q})$ to rewrite 
$J_{g,\beta,n-2}^{(1)}$ as a class in 
\[
H_{\bullet}(Y,Y_{-\infty};\mathbb{Q})\otimes
H_{\bullet}(Y,Y_{+\infty};\mathbb{Q}) \otimes H_{\bullet}(Y,
\mathbb{Q})^{\otimes (n-2)}\otimes
H_{\bullet}\left(\overline{\mathcal{M}}_{g,n},\mathbb{Q}\right). 
\]
This class can again be defined as a virtual fundamental class 
space $\overline{\mathcal{M}}_{g,n}(Y,\beta)$ 
of stable pseudo-holomorphic maps. Again we can interpret the virtual
class as a relative homology class: 
\[
\begin{split}
\left[\overline{\mathcal{M}}_{g,n}(Y,\beta)\right]_{\op{vir}} & \in
H_{\bullet}\left(Y^{n}\times \overline{\mathcal{M}}_{g,n},\,
Y^{n-2}\times \left((Y_{R,\varepsilon}^{-}\times Y) \cup (Y \times
Y_{R,\varepsilon}^{+})\right) \times
\overline{\mathcal{M}}_{g,n}; \mathbb{Q}\right) \\
& = H_{\bullet}(Y,Y_{R,\varepsilon}^{-};\mathbb{Q})\otimes
H_{\bullet}(Y,Y_{R,\varepsilon}^{+};\mathbb{Q}) \otimes H_{\bullet}(Y,
\mathbb{Q})^{\otimes (n-2)}\otimes
H_{\bullet}\left(\overline{\mathcal{M}}_{g,n},\mathbb{Q}\right)
\\
& = H_{\bullet}(Y,Y_{-\infty};\mathbb{Q})\otimes
H_{\bullet}(Y,Y_{+\infty};\mathbb{Q}) \otimes H_{\bullet}(Y,
\mathbb{Q})^{\otimes (n-2)}\otimes
H_{\bullet}\left(\overline{\mathcal{M}}_{g,n},\mathbb{Q}\right),
\end{split}
\]
and so it gives the desired map $J_{g,\beta,n-2}^{(1)}$.

\begin{figure}[!ht]
\begin{center}
\psfrag{0}[c][c][1][0]{{$0$}}
\psfrag{D}[c][c][1][0]{{
$\boldsymbol{\mathfrak{D}}_{\boldsymbol{\varepsilon}}$}}
\epsfig{file=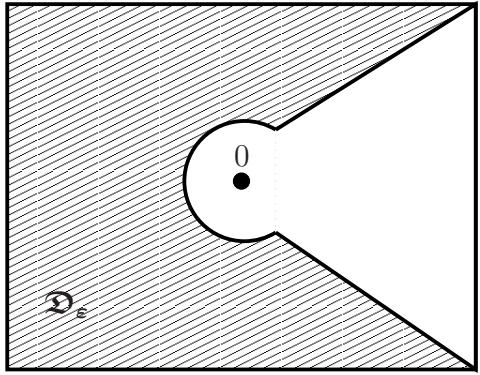,width=2.5in} 
\end{center}
\caption{The domain $\mathfrak{D}_{\varepsilon}$.}\label{fig:epsilon}
\end{figure}

Here $1 \gg \varepsilon > 0$, and $Y_{R,\varepsilon}^{\pm} =
\bw^{-1}(\pm \mathfrak{D}_{\varepsilon})$, where
$\mathfrak{D}_{\varepsilon} \subset \mathbb{C}$ is 
the domain given by (see Figure~\ref{fig:epsilon})
\[
\mathfrak{D}_{\varepsilon} := \left\{ z \in \mathbb{C} \; \left| \;
|z| \geq R \text{ and } \op{Arg} z \in \left( \frac{\pi}{2} -
\varepsilon, \frac{3\pi}{2} + \varepsilon \right) \right.\right\}.
\]
Again the relative invariants $J_{g,\beta,n-2}^{(1)}$ give rise to a
quantum multiplication and through the usual formulas from
Examples~\ref{ex:extended.symplectic} and \ref{ex:can.basic} we again
get a decorated variation of \nc-Hodge structures over a
(conjecturally non-empty) domain in $H^{\bullet}(Y,\mathbb{C})$ with
fiber $H^{\bullet}(Y,Y_{-\infty})$.

\

\subsubsection {\bfseries Mirror symmetry.} \ In conclusion we
systematize all the objects introduced above in a mirror 
table (see also \cite{auroux-complement}) 
describing the corresponding $A$ and $B$-model entities in
parallel:

\

\medskip

\begin{center}
\begin{tabular}{|c||c|c|}
\hline
invariants &
\begin{minipage}[c]{2.5in} \addtolength{\baselineskip}{-4pt}
\begin{center}
$A$-model
\end{center}
\end{minipage} &
\begin{minipage}[c]{2.3in} \addtolength{\baselineskip}{-4pt}
\begin{center}
$B$-model
\end{center}
\end{minipage} \\
\hline
\hline
geometry
&
\begin{minipage}[c]{2.5in} \addtolength{\baselineskip}{-4pt}
\

a triple $(Y,\bw,\omega)$ where: \\

$\bw : Y \to \mathbb{C}$  is a proper \linebreak $C^{\infty}$-map
$(Y,\omega)$ is symplectic with $c_{1}(T_{Y}) = 0$

\

\end{minipage} &
\begin{minipage}[c]{2.3in} \addtolength{\baselineskip}{-4pt}
\

a pair $Z\subset X$  where:\\

$X$ is smooth projective, and \linebreak
$Z \subset X$ is a smooth anticanonical
  divisor 

\

\end{minipage} \\
\hline
\hline
cohomology &
\begin{minipage}[c]{2.5in} \addtolength{\baselineskip}{-4pt}

\

$\left.
\text{
\begin{minipage}[c]{0.65in}\addtolength{\baselineskip}{-4pt}
$H^{\bullet}(Y,\mathbb{C})$ \\
$H^{\bullet}(Y,Y_{-\infty};\mathbb{C})$ \\
$H^{\bullet}(Y_{-\infty},\mathbb{C})$
\end{minipage}
} \qquad \right\}$
\ \begin{minipage}[c]{0.9in}\addtolength{\baselineskip}{-4pt}
variations of \nc HS over a domain in 
$H^{\bullet}(Y,\mathbb{C})$
\end{minipage}

\

\end{minipage} & 
\begin{minipage}[c]{2.5in} \addtolength{\baselineskip}{-4pt}

\

$\left.
\text{
\begin{minipage}[c]{1in} \addtolength{\baselineskip}{-4pt}
$H^{\bullet}(X-Z,\mathbb{C})$ \\
$H^{\bullet}(X,\mathbb{C})$ \\
$H^{\bullet}(Z,\mathbb{C})$
\end{minipage}
}  \right\}$ \hspace{0.5pc}
\begin{minipage}[c]{0.95in} \addtolength{\baselineskip}{-4pt}
variations of \nc HS over a domain in 
$H^{\bullet}(X-Z,\mathbb{C})$
\end{minipage}

\

\end{minipage}\\
\hline
\hline
\multirow{4}{*}{categories} &
\begin{minipage}[c]{2.5in} \addtolength{\baselineskip}{-4pt}

\

\begin{minipage}[c]{0.6in} \addtolength{\baselineskip}{-4pt}
$\xymatrix{
\Fuk(Y_{-\infty})  \\
\FS(Y) \ar[u]^-{F}
}$
\end{minipage}
\hspace{0.5pc} :
\begin{minipage}[c]{1.2in} \addtolength{\baselineskip}{-4pt}
$\Fuk(Y_{-\infty})$ is a CY category and  $F$ is a spherical functor
\end{minipage}

\

\end{minipage} &
\begin{minipage}[c]{2.5in} 

\

\begin{minipage}[c]{0.6in}
$\xymatrix{
D(Z)  \\
D(X) \ar[u]^-{F}
}$
\end{minipage}
\hspace{1pc} :
\begin{minipage}[c]{1.2in} \addtolength{\baselineskip}{-4pt}
$D(Z)$ is a CY category and  $F$ is a spherical functor
\end{minipage}

\

\end{minipage} \\
\cline{2-3}
&
\begin{minipage}[c]{2.5in} \addtolength{\baselineskip}{-4pt}

\

The part of $\FS(Y)$ consisting of  Lagrangians fibered over
\epsfig{file=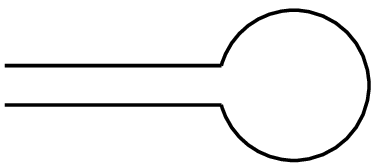,width=0.5in}
where the circle is of radius $R \gg 0$. 

\

\end{minipage}
& 
\begin{minipage}[c]{2.5in} \addtolength{\baselineskip}{-4pt}

\

$D_{\op{supp} Z}(X)$ : 
\begin{minipage}[c]{1.4in}\addtolength{\baselineskip}{-4pt}
a full compact (non smooth) subcategory in $D(X)$
\end{minipage}

\

\end{minipage}
\\
\cline{2-3}
& 
\begin{minipage}[c]{2.5in} \addtolength{\baselineskip}{-4pt}
\

The part of $\FS(Y)$ consisting \linebreak 
of compact Lagrangian sumanifolds in
$Y$ 

\

\end{minipage}
& 
\begin{minipage}[c]{2.5in} \addtolength{\baselineskip}{-4pt}
\

$D_{\substack{\text{compact} \\ \text{support}}}(X-Z)$ : 
\begin{minipage}[c]{1.2in}\addtolength{\baselineskip}{-4pt}
a full compact (non smooth) subcategory in $D(X)$
\end{minipage}

\

\end{minipage}
\\
\cline{2-3}
& 
\begin{minipage}[c]{2.5in} \addtolength{\baselineskip}{-4pt}
\

The wrapped $\FS$ category: the Hom space
between $(L_{1},\mathbb{V}_{1})$ and $(L_{2},\mathbb{V}_{2})$ is the
sum of  $\op{Hom}(\mathbb{V}_{1},\mathbb{V}_{2})_{x}$, $x \in L_{1}
\cap L_{2}$, and $L_{2}$ is deformed so that $\bw(L_{2})$ becomes a spiral:
\begin{center}
\epsfig{file=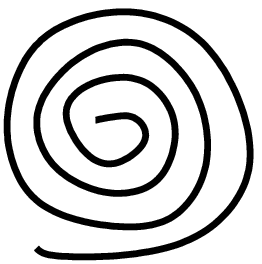,width=0.5in} 
\end{center}
wrapped infinitely many times

\

\

\end{minipage}
& 
\begin{minipage}[c]{2.5in} \addtolength{\baselineskip}{-4pt}
\

$D(X-Z)$ : 
\begin{minipage}[c]{1.4in}\addtolength{\baselineskip}{-4pt}
a smooth (non compact) category
\end{minipage}

\

\end{minipage}

\\
\hline
\end{tabular}
\end{center}

\bibliographystyle{halpha} 
\bibliography{mirror,catinvariants}

\

\bigskip

\noindent
Ludmil Katzarkov, {\sc University of Miami}, 
l.katzarkov@math.miami.edu

\smallskip

\noindent
Maxim Kontsevich, {\sc IHES, and University of Miami}, maxim@ihes.fr

\smallskip

\noindent
Tony Pantev, {\sc University of Pennsylvania}, tpantev@math.upenn.edu

\end{document}